
\input amstex
\documentstyle{amsppt}
\loadbold

\magnification=\magstep1
\hsize=6.5truein
\vsize=9.1truein

\document

\baselineskip=14pt

\font\smallit=cmti10 at 9pt
\font\smallsl=cmsl10 at 9pt
\font\sc=cmcsc10

\def \loongrightarrow {\relbar\joinrel\relbar\joinrel\rightarrow}
\def \llongrightarrow
{\relbar\joinrel\relbar\joinrel\relbar\joinrel\rightarrow}

\def \gerg {\frak g}
\def \gerk {\frak k}

\def \gerh {\frak h}
\def \germ {\frak m}

\def \id {\hbox{\rm id}}
\def \und1 {\underline{1}}
\def \undepsilon {\underline{\epsilon}}
\def \u {{\bold u}}

\def \L {\Cal{L}}
\def \U {\Cal{U}}
\def \calA {\Cal{A}}
\def \calB {\Cal{B}}
\def \calC {\Cal{C}}

\def \HA {\Cal{HA}}
\def \A {\Bbb{A}}

\def \N {\Bbb{N}}
\def \Z {\Bbb{Z}}

\def \hyp {{\text{\sl Hyp}\,}}
\def \Char {\hbox{\sl Char}\,}


\topmatter

\title
  The Crystal Duality Principle: from Hopf Algebras
to Geometrical Symmetries
\endtitle

\author
       Fabio Gavarini
\endauthor

\rightheadtext{ The crystal duality principle }

\affil
  Universit\`a degli Studi di Roma ``Tor Vergata'' ---
Dipartimento di Matematica  \\
  Via della Ricerca Scientifica 1, I-00133 Roma --- ITALY  \\
\endaffil

\address\hskip-\parindent
  Fabio Gavarini  \newline
     \indent   Universit\`a degli Studi di Roma ``Tor Vergata''  
---   Dipartimento di Matematica  \newline
     \indent   Via della Ricerca Scientifica 1, I-00133 Roma,
ITALY ---   e-mail: gavarini\@{}mat.uniroma2.it
\endaddress

\abstract
   We give functorial recipes to get, out of any Hopf algebra over
a field, two pairs of Hopf algebras with some geometrical content.
If the ground field has characteristic zero, the first pair is made
by a function algebra  $ F[G_+] $  over a connected Poisson group and
a universal enveloping algebra  $ U(\gerg_-) $  over a Lie bialgebra
$ \gerg_- \, $.  In addition, the Poisson group as a variety is an
affine space, and the Lie bialgebra as a Lie algebra is graded. 
Forgetting these last details, the second pair is of the same type,
namely  $ \, \big( F[K_+], U(\gerk_-) \big) \, $  for some Poisson
group  $ K_+ $  and some Lie bialgebra  $ \gerk_- \, $.  When the
Hopf algebra  $ H $  we start from is already of geometric type the
result involves Poisson duality.  The first Lie bialgebra associated
to  $ \, H = F[G] \, $  is  $ \, \gerg^* \, $  (with  $ \, \gerg
:= \text{\it Lie}\,(G) \, $),  \, and the first Poisson group
associated to  $ \, H = U(\gerg) \, $  is of type  $ G^* $,  i.e.{}
it has  $ \gerg $  as cotangent Lie bialgebra.  If the ground field
has positive characteristic, the same recipes give similar results, but
the Poisson groups obtained have dimension 0 and height 1, and  {\sl
restricted\/}  universal enveloping algebras are obtained.  We show
how these geometrical Hopf algebras are linked to the initial one via
1-parameter deformations, and explain how these results follow from
quantum group theory.  
   \hbox{We examine in detail the case of group algebras.}   
\endabstract

\endtopmatter

%
%
\footnote""{\noindent   Keywords:  {\sl filtered, graded Hopf algebras;
Poisson algebraic groups; Lie bialgebras; quantum groups}.}

\footnote""{ 2000 {\it 
%
%
 MSC:} \ 
Primary 16W30, 16W70, 16W50; Secondary 14L17, 16S30, 17B62--63, 81R50. }

\vskip-9pt

\hfill  \hbox{\vbox{ \baselineskip=10pt
 \hbox{\smallit \hskip5pt ``Yet these crystals are to Hopf algebras }
 \hbox{\smallit \hskip1pt but as is the body to the Children of Rees: }
 \hbox{\smallit   the house of its inner fire, that is within it }
 \hbox{\smallit \hskip4pt  and yet in all parts of it,
and is its life'' }
      \vskip7pt
  \hbox{\smallsl  \hskip37pt  N.~Barbecue, ``Scholia'' } }
\hskip1truecm }

\vskip19pt

\centerline {\bf Introduction }

\vskip10pt

   Among all Hopf algebras over a field  $ \Bbbk $,  there are two
special families which are of relevant interest for their geometrical
meaning.  The function algebras  $ F[G] $  of algebraic groups  $ G $
and the universal enveloping algebras  $ U(\gerg) $  of Lie algebras
$ \gerg $,  if  $ \, \Char(\Bbbk) = 0 \, $,  or the restricted
universal enveloping algebras  $ \u(\gerg) $  of restricted Lie
algebras  $ \gerg \, $, if  $ \, \Char(\Bbbk) > 0 \, $.  For brevity,
we call both  $ U(\gerg) $  and  $ \u(\gerg) $  ``enveloping algebras''
and denote them by  $ \U(\gerg) \, $.  Similarly by ``restricted Lie
algebra'', when  $ \, \Char(\Bbbk) = 0 \, $,  we shall simply mean
``Lie algebra''.  Function algebras are exactly those Hopf algebras
which are commutative, and enveloping algebras those which are
connected,
      \hbox{cocommutative and generated by their primitives.}
                                         \par
   In this paper we give functorial recipes to get out from any Hopf
algebra two pairs of Hopf algebras of geometrical type, say  $ \big(
F[G_+], \U(\gerg_-) \big) $  and  $ \big( F[K_+], \U(\gerk_-) \big) $. 
In addition, the algebraic groups obtained in this way are connected  {\sl
Poisson\/}  groups,  and the (restricted) Lie algebras are (restricted) 
{\sl Lie bialgebras}.  Therefore, to each Hopf algebra   --- which encodes
a general notion of ``symmetry'' ---   we can associate in a functorial
way some symmetries of geometrical type, where the geometry involved is in
fact Poisson geometry.  Moreover, if  $ \, \Char(\Bbbk) > 0 \, $  these
Poisson groups have dimension 0 and height 1, which makes them very
interesting for arithmetic geometry, hence for number theory too.   
                                            \par
   The construction of the pair  $ (G_+,\gerg_-) $  uses pretty
classical (as opposite to ``quantum'') methods: in fact, it might
be part of the content of any basic textbook on Hopf algebras
(and, surprisingly enough, it is not!).  Instead, in order to
obtain the pair  $ (K_+,\gerk_-) $  one relies on the construction
of the first pair, and uses the theory of quantum groups.
                                            \par
   Let's describe our results in detail.  Let  $ \, J := \text{\sl
Ker}\,(\epsilon_{\scriptscriptstyle H}) \, $  be the augmentation
ideal of the Hopf algebra  $ H $  (where  $ \epsilon_{\scriptscriptstyle
H} $  is the counit of  $ H $),  and let  $ \, \underline{J} :=
{\big\{ J^n \big\}}_{n \in \N} \, $  be the associated  $ J $--adic
filtration,  $ \, \widehat{H} := G_{\underline{J}}(H) \, $  the
associated graded vector space and  $ \, H^\vee := H \Big/
\bigcap_{n \in \N} J^n \, $.  One proves that  $ \underline{J} $
is a Hopf algebra filtration, hence  $ \widehat{H} $  is a graded
Hopf algebra: the latter happens to be connected, cocommutative and
generated by its primitives, so  $ \, \widehat{H} \cong \U(\gerg_-)
\, $  for some restricted Lie algebra  $ \gerg_- \, $.  In addition,
since  $ \widehat{H} $  is graded also  $ \gerg_- $  itself is graded
as a restricted Lie algebra.  The fact that  $ \widehat{H} $  be
cocommutative allows to define a Poisson cobracket on it (from the
natural Poisson cobracket  $ \, \nabla := \Delta - \Delta^{\text{op}}
\, $  on  $ H $),  which makes  $ \widehat{H} $  a graded  {\sl
co-Poisson\/}  Hopf algebra.  Eventually, this implies that
$ \gerg_- $  is a  {\sl Lie bialgebra}.  So the right-hand
side half of the first pair of ``Poisson geometrical''
     \hbox{Hopf algebras is just  $ \widehat{H} $.}   
                                            \par
   On the other hand, one considers a second filtration   ---
increasing, whereas  $ \underline{J} $  is decreasing ---
namely  $ \, \underline{D} \, $  which is defined in a dual
manner to  $ \underline{J} \, $.  For each  $ \, n \in \N \, $, 
let  $ \delta_n $  be the composition of the  $ n $--fold  iterated
coproduct followed by the projection onto  $ J^{\otimes n} $  (note
that  $ \, H = \Bbbk \! \cdot \! 1_{\scriptscriptstyle H} \oplus J
\, $);  then  $ \, \underline{D} := {\big\{ D_n := \text{\sl Ker}\,
(\delta_{n+1}) \big\}}_{n \in \N} \, $.  Let now  $ \, \widetilde{H}
:= G_{\underline{D}}(H) \, $  be the associated graded vector space
and  $ \, H' := \bigcup_{n \in \N} D_n \, $.  Again, one shows
that  $ \underline{D} $  is a Hopf algebra filtration, hence
$ \widetilde{H} $  is a graded Hopf algebra: moreover, the latter
is commutative, so  $ \, \widetilde{H} = F[G_+] \, $  for some
algebraic group  $ G_+ $.  One proves also that  $ \, \widetilde{H}
= F[G_+] \, $  has no non-trivial idempotents, thus  $ G_+ $  is
connected; in addition, since  $ \widetilde{H} $  is graded,  $ G_+ $
as a variety is just an affine space.  A deeper analysis shows that
in the positive characteristic case  $ G_+ $  has dimension 0 and
height 1.  The fact that  $ \widetilde{H} $  be commutative
allows to define on it a Poisson bracket (from the natural
Poisson bracket on  $ H $  given by the commutator) which makes
$ \widetilde{H} $  a graded  {\sl Poisson\/}  Hopf algebra.  This
means that  $ G_+ $  is an algebraic  {\sl Poisson group}.  So the
   \hbox{left-hand side half of the first pair of ``Poisson
geometrical'' Hopf algebras is just  $ \widetilde{H} $.}
                                            \par
   The relationship among  $ H $  and the ``geometrical'' Hopf
algebras  $ \widehat{H} $  and  $ \widetilde{H} $  can be expressed
in terms of ``reduction steps'' and regular 1-parameter deformations,
namely
  $$  \widetilde{H}  \hskip1pt
\underset{ {\Cal{R}^t_{\underline{D}}(H)}}  \to
{\overset{0 \,\leftarrow\, t \,\rightarrow\, 1} \to
{\longleftarrow\joinrel\relbar\joinrel\relbar\joinrel%
\relbar\joinrel\llongrightarrow}}  \hskip1pt  H'  \hskip1pt
\,{\lhook\joinrel\relbar\joinrel\relbar\joinrel\rightarrow}\,
\hskip1pt  H  \relbar\joinrel\relbar\joinrel\twoheadrightarrow
H^\vee  \hskip0pt  \underset{ {\Cal{R}^t_{\underline{J}}(H^\vee)}}
\to  {\overset{1 \,\leftarrow\, t \,\rightarrow\, 0} \to
{\longleftarrow\joinrel\relbar\joinrel\relbar\joinrel%
\relbar\joinrel\relbar\joinrel\relbar\joinrel\longrightarrow}}
\hskip1pt   \widehat{H}   \eqno (\bigstar)  $$
Here the one-way arrows are Hopf algebra morphisms and
the two-ways arrows are regular 1-parameter deformations
of Hopf algebras, realized through the  {\sl Rees\/}
Hopf algebras  $ {\Cal{R}^t_{\underline{D}}(H)} $  and
$ {\Cal{R}^t_{\underline{J}} (H^\vee)} $  associated to
the filtration  $ \underline{D} $  of  $ H $  and to
the filtration  $ \underline{J} $  of  $ H^\vee $.
                                            \par
   The construction of the pair  $ (K_+,\gerk_-) $
uses quantum group theory, the basic ingredients
being  $ {\Cal{R}^t_{\underline{D}}(H)} $  and
$ {\Cal{R}^t_{\underline{J}} (H^\vee)} $.  In
the present framework, by quantum group we mean,
loosely speaking, a Hopf  $ \Bbbk[t] $--algebra
($ t $  an indeterminate)  $ H_t $  such that either \,
{\it (a)}  $ \, H_t \Big/ t \, H_t \cong F[G] \, $  for some
connected Poisson group  $ G $   --- then we say  $ H_t $  is a
QFA ---   or \,  {\it (b)}  $ \, H_t \Big/ t \, H_t \cong \U(\gerg)
\, $,  \, for some restricted Lie bialgebra  $ \gerg $   --- then we
say  $ H_t $  is a QrUEA.  Formula  $ (\bigstar) $  says that  $ \;
H'_t := {\Cal{R}^t_{\underline{D}}(H)} \; $  is a QFA, with  $ \;
H'_t \Big/ t \, H'_t \, \cong \, \widetilde{H} = F[G_+] \, $,  \;
and that  $ \; H^\vee_t := {\Cal{R}^t_{\underline{J}}(H)} \; $
is a QrUEA, with  $ \; H^\vee_t \Big/ t \, H^\vee_t \, \cong \,
\widehat{H} \, = \U(\gerg_-) \, $.  Now, a general result   ---
the ``Global Quantum Duality Principle'', in short GQDP, see [Ga1--2]
---  teaches us how to construct from the QFA  $ H'_t $  a QrUEA, call
it  $ \big( H'_t \big)^\vee $,  and how to build out of the QrUEA
$ H^\vee_t $  a QFA, say  $ \big( H^\vee_t \big)' $.  Then  $ \,
\big(H'_t\big)^\vee \Big/ t \, \big(H'_t\big)^\vee = \U(\gerk_-) \, $
for some restricted Lie bialgebra  $ \gerk_- \, $,  and  $ \, \big(
H^\vee_t \big)' \Big/ t \, \big( H^\vee_t \big)' = F[K_+] \, $  for
some connected Poisson group  $ K_+ \, $.  This gives the pair
$ (K_+,\gerk_-) \, $.  The very construction implies that  $ \big(
H'_t \big)^\vee $  and  $ \big( H^\vee_t \big)' $  yield another
frame of regular 1-parameter deformations for  $ H' $  and
$ H^\vee $,  namely
  $$  \U(\gerk_-)  \hskip1pt
\underset{ (H'_t)^\vee}  \to
{\overset{0 \,\leftarrow\, t \,\rightarrow\, 1}
\to {\longleftarrow\joinrel\relbar\joinrel\relbar%
\joinrel\relbar\joinrel\llongrightarrow}}  \hskip1pt  H'
\hskip1pt  \,{\lhook\joinrel\relbar\joinrel\rightarrow}\,
\hskip1pt  H  \relbar\joinrel\relbar\joinrel\twoheadrightarrow
H^\vee  \hskip0pt  \underset{(H^\vee_t)'}
\to  {\overset{1 \,\leftarrow\, t \,\rightarrow\, 0} \to
{\longleftarrow\joinrel\relbar\joinrel\relbar\joinrel\relbar%
\joinrel\relbar\joinrel\relbar\joinrel\longrightarrow}}
\hskip1pt   F[K_+]   \eqno (\maltese)  $$
which is the analogue of  $ (\bigstar) $.  In addition, when
$ \, \Char(\Bbbk) = 0 \, $  the GQDP also claims that the two
pairs  $ (G_+,\gerg_-) $  and  $ (K_+,\gerk_-) $  are related
by Poisson duality: namely,  $ \gerk_- $  {\sl is the cotangent
Lie bialgebra of\/}  $ G_+ \, $,  and  $ \gerg_- $  {\sl is the
cotangent Lie bialgebra\/}  of  $ K_+ \, $;  \, in short, we write
$ \, \gerk_- = \gerg^\times \, $  and  $ \, K_+ = G_-^\star \, $.
Therefore the four ``Poisson symmetries''  $ G_+ $,  $ \gerg_- $,
$ K_+ $  and  $ \gerk_- $,  attached to  $ H $  are actually encoded
simply by the pair  $ (G_+,K_+) $.
                                            \par
   In particular, when  $ \, H^\vee \! = H = H' \, $  from
$ (\bigstar) $  and  $ (\maltese) $  together we find
%
%
%
  $$  \hskip9pt   \U(\gerg_-)  \hskip2pt
\underset{H^\vee_t}  \to
{\overset{0 \,\leftarrow\, t \,\rightarrow\, 1}
\to{\longleftarrow\joinrel\relbar\joinrel%
\relbar\joinrel\relbar\joinrel\relbar\joinrel\llongrightarrow}}
\hskip2pt {}  H^\vee  \hskip1pt
\underset{\;(H^\vee_t)'}  \to
{\overset{1 \,\leftarrow\, t \,\rightarrow\, 0}
\to{\longleftarrow\joinrel\relbar\joinrel%
\relbar\joinrel\relbar\joinrel\relbar\joinrel\llongrightarrow}}
\hskip2pt  F[K_+]
\hskip9pt  \Big( = F\big[G_-^\star\big]  \,\text{\ if \ }
\!\text{\it Char}\,(\Bbbk) = 0 \,\Big)  $$
 \vskip-25pt
  $$  \hskip-135pt   ||  $$
 \vskip-21pt
  $$  \hskip-129pt  H_{\phantom{\displaystyle I}}  $$
 \vskip-23pt
  $$  \hskip-135pt   ||  $$
 \vskip-29pt
  $$  \hskip5pt   F[G_+]  \hskip3pt
\underset{\;H'_t}  \to
{\overset{0 \,\leftarrow\, t \,\rightarrow\, 1}
\to{\longleftarrow\joinrel\relbar\joinrel%
\relbar\joinrel\relbar\joinrel\relbar\joinrel\llongrightarrow}}
\hskip2pt {}  H'  \hskip1pt
\underset{\;(H'_t)^\vee}  \to
{\overset{1 \,\leftarrow\, t \,\rightarrow\, 0}
\to{\longleftarrow\joinrel\relbar\joinrel%
\relbar\joinrel\relbar\joinrel\relbar\joinrel\llongrightarrow}}
\hskip3pt  \U(\gerk_-)
\hskip11pt  \Big( = U\big(\gerg_+^{\,\times}\big)  \,\text{\ if \ }
\text{\it Char}\,(\Bbbk) = 0 \,\Big)  $$
\noindent   This gives  {\sl four\/}  different regular 1-parameter
deformations from  $ H $  to Hopf algebras encoding geometrical objects
of Poisson type, i.e.~Lie bialgebras or Poisson algebraic groups.
                                            \par
   When the Hopf algebra  $ H $  we start from is already of geometric
type, the result involves Poisson duality.  Namely, if  $ \, \Char
(\Bbbk) = 0 \, $  and  $ \, H = F[G] \, $,  then  $ \, \gerg_- = \gerg^*
\, $  (where  $ \, \gerg := \text{\sl Lie}\,(G) \, $),  \, and if
$ \, H = \U(\gerg) = U(\gerg) \, $,  then  $ \, \text{\sl Lie}\,(G_+)
= \gerg^* $,  i.e.~$ G_+ $  has  $ \gerg $  as cotangent Lie bialgebra.
If instead  $ \, \Char(\Bbbk) > 0 \, $,  we have only a slight variation
on this result.
                                            \par
   The construction of  $ \widehat{H} $  and  $ \widetilde{H} $
needs only ``half the notion'' of a Hopf algebra.  In fact, we
construct  $ \widehat{A} $  for  any augmented algebra  $ A $
(roughly, an algebra with an augmentation, or counit, i.e.~a
character), and  $ \widetilde{C} $  for any coaugmented coalgebra
$ C $  (a coalgebra with a coaugmentation, or unit, i.e.~a coalgebra
morphism from  $ \Bbbk $  to  $ C $).  In particular this applies to
bialgebras, for which both  $ \widehat{B} $  and  $ \widetilde{B} $  are
(graded) Hopf algebras.  We can also perform a second construction using
$ \big( B'_t \big)^\vee $  and  $ \big( B^\vee_t \big)' $  (via a stronger
version of the GQDP), and get from these a second pair of bialgebras
$ \Big( \big( B'_t \big)^\vee{\Big|}_{t=0} , \big( B^\vee_t \big)'
{\Big|}_{t=0} \Big) \, $.  Then again  $ \, \big( B'_t \big)^\vee
{\Big|}_{t=0} \! \cong \U(\gerk_-) \, $  for some restricted Lie
bialgebra  $ \gerk_- \, $,  \, while  $ \big( B^\vee_t \big)'
{\Big|}_{t=0} $  is commutative with no non-trivial idempotents, but
it's not, in general, a Hopf algebra.  So the spectrum of  $ \big(
B^\vee_t \big)'{\Big|}_{t=0} $  is a connected algebraic Poisson
monoid, but not necessarily a Poisson group.
                                            \par
   It is worthwhile pointing out that everything in fact follows from
the GQDP, which   --- in the stronger formulation ---   deals with
augmented algebras and coaugmented coalgebras over 1-dimensional
domains.  The content of this paper can in fact be obtained as
a corollary of the GQDP as follows.  Pick any augmented algebra
or coaugmented coalgebra over  $ \Bbbk $,  and take its scalar
extension from  $ \Bbbk $  to  $ \Bbbk[t] \, $:  the latter ring
is a 1-dimensional domain, hence we can apply the GQDP, and
(almost) every result in the present paper will follow.
                                            \par
   In the last section we apply these results to the case of group
algebras and their duals.  Another interesting application   ---
based on a non-commutative version of the function algebra of the
group of formal diffeomorphism on the line, also called ``Nottingham
group'' ---   is illustrated in detail in a separate paper (see [Ga3]).  

\vskip15pt

\centerline{ \sc acknowledgements }

\vskip1pt

  The author thanks N.~Andruskiewitsch, G.~Carnovale, L.~Foissy,
C.~Gasbarri, C.~Menini, H.-J.~Schneider, C.~M.~Scoppola,
E.~Strickland and E.~Taft for many helpful conversations.

\vskip1,3truecm

\centerline {\bf \S \; 1 \ Notation and terminology }

\vskip10pt

  {\bf 1.1 Algebras, coalgebras, and further structures.} \,
Let  $ \Bbbk $  be a field, which will stand fixed throughout, with 
$ \, p := \text{\sl Char}\,(\Bbbk) \, $.  In this paper we deal with
(unital associative)  $ \Bbbk $--algebras  and (counital coassociative) 
$ \Bbbk $--coalgebras  in the standard sense, cf.~[Sw] or [Ab].  In
particular we use notations as in [Ab].  For any (counital coassociative) 
$ \Bbbk $--coalgebra  $ C $  we denote by  $ \text{\sl coRad}\,(C) $ 
its coradical, and by  $ \, G(C) := \big\{\, c \in C \,\big\vert\,
\Delta(c) = c \otimes c \,\big\} \, $  its set of group-like elements. 
We say  $ C $  is  {\sl monic\/}  if  $ \, \big| G(C) \big| = 1 \, $; 
we say it is  {\sl connected\/}  if  $ \text{\sl coRad} \,(C) $  is
one-dimensional:  \, of course ``connected'' implies ``monic''. 
Any unital associative  $ \Bbbk $--alge\-bra  $ A $  is said  {\sl
idempotent-free\/}  (in short,  {\sl i.p.-free\/}) iff it has
no non-trivial idempotents.   
                                            \par
  We call  {\sl augmented algebra\/}  any unital associative
$ \Bbbk $--algebra  $ A $  {\sl together with\/}  a special unital
algebra morphism  $ \, \undepsilon : \, A \longrightarrow \Bbbk \, $ 
(so the unit  $ \, u : \, \Bbbk \longrightarrow A \, $  is a section
of  $ \undepsilon \, $):  \, these form a category in the obvious way.
We call  {\sl indecomposable elements\/}  of an augmented algebra
$ A $  the elements of the set  $ \, Q(A) := J_{\scriptscriptstyle A}
\big/ {J_{\scriptscriptstyle A}}^{\hskip-3pt 2} \, $  with  $ \,
J_{\scriptscriptstyle A} := \text{\sl Ker}\,\big(\, \undepsilon \,
\colon A \longrightarrow \Bbbk \,) \, $.  We denote by  $ \calA^+ $ 
the category of all augmented  $ \Bbbk $--algebras.  
                                            \par
   We call  {\sl coaugmented coalgebra\/}  any counital coassociative 
$ \Bbbk $--coalgebra  $ C $  {\sl together with\/}  a special counital
coalgebra morphism  $ \, \underline{u} : \, \Bbbk \longrightarrow C
\, $  (so  $ \underline{u} $  is a section of  $ \, \epsilon : \, C
\longrightarrow \Bbbk \, $),  and let  $ \, \und1 := \underline{u}(1)
\, $,  \, a group-like element in  $ C \, $:  \, these form a category
in the obvious way.  For such a  $ C $  we said  {\sl primitive\/}  the
elements of the set  $ \, P(C) := \big\{ c \in C \,|\; \Delta(c) = c
\otimes \und1 + \und1 \otimes c \,\big\} \, $.  We denote by  $ \calC^+ $ 
the category of all coaugmented   $ \Bbbk $--coalgebras.
                                            \par
   We denote by  $ \calB $  the category of all  $ \Bbbk $--bialgebras. 
Clearly each bialgebra  $ B $  can be seen both as an augmented
algebra,  w.r.t.~$ \, \undepsilon = \epsilon \equiv \epsilon_{\!
\scriptscriptstyle B} \, $  (the counit of  $ B \, $)  and as a
coaugmented coalgebra, w.r.t.~$ \, \underline{u} = u \equiv
u_{\! \scriptscriptstyle B} \, $  (the unit map of  $ B \, $),  so
that  $ \, \und1 = 1 = 1_{\! \scriptscriptstyle B} \, $:  then
$ Q(B) $  is naturally a Lie coalgebra and  $ P(B) $  a Lie algebra
over  $ R \, $.  In the following we'll use such an interpretation
throughout, looking at objects of  $ \calB $  as objects of
$ \calA^+ $  and of  $ \calC^+ $.  We call  $ \HA $  the category
of all Hopf  $ \Bbbk $--algebras;  this naturally identifies with
a subcategory of  $ \calB $.
                                              \par
   We call  {\sl Poisson algebra\/}  any (unital) commutative algebra
$ A $  endowed with a Lie bracket  $ \; \{\ ,\ \} : A \otimes A
\longrightarrow A \; $  (i.e.,  $ \, \big( A, \{\ ,\ \} \big) \, $
is a Lie algebra) such that the  {\sl Leibnitz identities}
  $$  \{a \, b, c\} = \{a, c\} \, b + a \, \{b, c\} \, ,
\quad  \{a, b \, c\} = \{a, b\} \, c + b \, \{a, c\}  $$
hold (for all  $ \, a $,  $ b $,  $ c \in A \, $).  We call  {\sl
Poisson bialgebra},  or  {\sl Poisson Hopf algebra},  any bialgebra,
or Hopf algebra, say  $ H $,  which is also a Poisson algebra
(w.r.t.~the same product) enjoying
  $$  \Delta\big(\{a,b\}\big) = \big\{\Delta(a),\Delta(b)\big\}
\, ,  \qquad  \epsilon\big(\{a,b\}\big) = 0 \, ,  \qquad
S\big(\{a,b\}\big) = \big\{S(b),S(a)\big\}  $$
(for all  $ \, a $,  $ b $,  $ c \in H \, $)   --- the condition on
the antipode   $ S $  being required in the Hopf algebra case ---
where the (Poisson) bracket on  $ \, H \otimes H \, $  is defined
by  $ \, \{a \otimes b, c \otimes d\,\} := \{a,b\} \otimes c \, d +
a \, b \otimes \{c,d\,\} \, $  (for all  $ \, a $,  $ b $,  $ c $,
$ d \in H \, $).
                                            \par
   We call  {\sl co-Poisson coalgebra\/}  any (counital) cocommutative
coalgebra  $ C $  with a Lie cobracket  $ \, \delta \,\colon\, C
\longrightarrow C \otimes C \, $  (i.e.~$ \, \big( C, \delta \big) \, $
is a Lie
       \hbox{coalgebra) such that the  {\sl co-Leibnitz identity}}
  $$  \big( id \otimes \Delta \big) \circ \big(\delta(a)\big) \; = \;
{\textstyle \sum_{(a)}} \, \Big( \delta\big(a_{(1)}\big) \otimes
a_{(2)} + \sigma_{1,2} \Big(a_{(1)} \otimes \delta\big(a_{(2)}\big)
\Big) \Big)  $$
holds for all  $ \, a \in C \, $,  \, where  $ \, \sigma_{1,2}:
\! C^{\otimes 3} \longrightarrow C^{\otimes 3} \, $  is given by  $ \,
x_1 \otimes x_2 \otimes x_3 \mapsto x_2 \otimes x_1 \otimes x_3 \, $.
We call  {\sl co-Poisson bialgebra},  or  {\sl co-Poisson Hopf
algebra},  any bialgebra, or Hopf algebra, say  $ H $,  which is
also a co-Poisson algebra (w.r.t.~the same co-product) enjoying
  $$  \delta(a \, b) = \delta(a) \, \Delta(b) + \Delta(a) \,
\delta(b) \, ,  \qquad  (\epsilon \otimes \epsilon)
\big(\delta(a)\big) = 0 \, ,  \qquad  \delta\big(S(a)\big)
= \tau \Big(\! \big(S \otimes S\big) \big(\delta(a)\big) \!\Big)  $$
(where  $ \tau $  is the flip) for all  $ \, a $,  $ b \in H \, $, 
the condition on the antipode  $ S $  being required in the Hopf
algebra case.  Finally, we call  {\sl bi-Poisson bialgebra},  or 
{\sl bi-Poisson Hopf algebra},  any bialgebra, or Hopf algebra, say 
$ H $,  which is simultaneously a Poisson and co-Poisson bialgebra,
or Hopf algebra,     
     \hbox{for some Poisson bracket and cobracket enjoying,
for all  $ \, a $,  $ b \in H \, $,}   
  $$  \delta\big(\{a,b\}\big) = \big\{ \delta(a), \Delta(b) \big\}
+ \big\{ \Delta(a), \delta(b) \big\} \; .  $$
See [CP] and [Tu1], [Tu2] for further details on the above notions.
                                             \par
   A  {\sl graded algebra\/}  is an algebra  $ A $  which is
$ \Z $--graded  as a vector space and whose structure maps
$ \, m $,  $ u $  and  $ S $  are morphisms of degree zero in
the category of graded vector spaces, where  $ \, A \otimes A \, $
has the standard grading inherited from  $ A $  and  $ \Bbbk $  has
the trivial grading.  Similarly we define the graded versions of
coalgebras, bialgebras and Hopf algebras, and also the graded
versions of Poisson algebras, co-Poisson coalgebras,
Poisson/co-Poisson/bi-Poisson bialgebras, and
Poisson/co-Poisson/bi-Poisson Hopf algebras, but for
the fact that the Poisson bracket, resp.~cobracket, must
be a morphism (of graded spaces) of degree  $ -1 $,  resp.~$ +1 $.
We write  $ \, V = \oplus_{z \in \Z} V_z \, $  for the degree
splitting of any graded vector space  $ V $.

\vskip7pt

  {\bf 1.2 Function algebras.} \, According to standard theory, the
category of commutative Hopf algebras  is antiequivalent to the
category of algebraic groups (over  $ \Bbbk $).  Then we call
$ \hbox{\it Spec}\,(H) $  ({\sl spectrum of  $ H \, $})  the image
of a Hopf algebra  $ H $  in this antiequivalence, and conversely
we call  {\sl function algebra\/}  or  {\sl algebra of regular
functions\/}  the preimage  $ F[G] $  of an algebraic group  $ G $.
Note that we do  {\sl not\/}  require algebraic groups to be reduced
(i.e.~$ F[G] $  to have trivial nilradical) and we do  {\sl not\/}
make any restrictions on dimensions: in particular we deal with
{\sl pro-affine\/}  as well as  {\sl affine\/}  algebraic groups.
We say that  $ G $  is  {\sl connected\/}  if  $ F[G] $  is
i.p.-free; this is equivalent to the classical topological
notion when  $ dim(G) $  is finite.
                                           \par
   Given an algebraic group  $ G $,  let  $ \, J_G := \text{\it Ker}
\,(\epsilon_{\scriptscriptstyle F[G]}) \; $;  \, the  {\sl cotangent
space of  $ G $\/}  (at its unity) is  $ \, \gerg^\times := J_G \Big/
{J_G}^{\hskip-3pt 2} = Q\big(F[G]\big) \, $,  \, endowed with its weak
topology.  The  {\sl tangent space of  $ G $}  (at its unity) is the
topological dual  $ \, \gerg := {\big( \gerg^\times \big)}^\star \, $
of  $ \gerg^\times \, $:  \, this is a Lie algebra, the  {\sl tangent
Lie algebra of  $ G $}.  If  $ \, p > 0 \, $,  \, then  $ \gerg $  is a
{\sl restricted\/}  Lie algebra (also called  ``$ p $--Lie  algebra'').
We say that  $ G $  is an algebraic  {\sl Poisson\/}  group if  $ \,
F[G] \, $  is a Poisson Hopf algebra.  Then the tangent Lie algebra
$ \gerg $  of  $ G $  is a  {\sl Lie bialgebra},  and the same holds
for  $ \gerg^\times $.  If  $ \, p > 0 \, $,  \, then  $ \gerg $
and  $ \gerg^\times $  are  {\sl restricted Lie bialgebras},  the
$ p $--operation  on  $ \gerg^\times $  being trivial.

\vskip7pt

  {\bf 1.3 Enveloping algebras and symmetric algebras.}  \, Given
a Lie algebra  $ \gerg \, $,  we denote by  $ U(\gerg) $  its  {\sl
universal enveloping algebra\/}.  If  $ \, p > 0 \, $  and  $ \gerg $
is a  {\sl restricted \/}  Lie algebra, we denote  $ \; \u(\gerg) =
U(\gerg) \Big/ \big( \big\{\, x^p - x^{[p\hskip0,5pt]} \,\big|\, x \in
\gerg \,\big\} \big) \, $  its  {\sl restricted universal enveloping
algebra\/}.  If  $ \, p = 0 \, $,  \, then  $ \, P \big( U(\gerg) \big)
= \gerg \, $;  \, if instead  $ \, p > 0 \, $,  \, then  $ \, P \big(
U(\gerg) \big) = \gerg^\infty := \text{\it Span}\,\Big(\! \big\{ x^{p^n}
\big\}_{n \in \N} \Big) \, $,  \, the latter carrying a natural structure
of restricted Lie algebra with  $ \, X^{[p\,]} := X^p \, $.
                                           \par
   Note that  $ \, U(\gerg) = \u(\gerg^\infty) \, $  for any Lie
algebra  $ \gerg \, $,  \, so any universal enveloping algebra can
be thought of as a restricted universal enveloping algebra.  Both
$ U(\gerg) $  and  $ \u(\gerg) $  are cocommutative connected Hopf
algebras, generated by  $ \gerg $  itself.  Conversely, if  $ \,
p = 0 \, $  then each cocommutative connected Hopf algebra is the
universal enveloping algebra of some Lie algebra, and if  $ \, p >
0 \, $  then each cocommutative connected Hopf algebra  $ H $  which
is generated by  $ P(H) $  is the restricted universal enveloping
algebra of some restricted Lie algebra (cf.~[Mo], Theorem 5.6.5, and
references therein).  Thus, in order to unify termino\-logy and notations,
we call both universal enveloping algebras (when  $ \, p = 0 \, $)  and
restricted universal enveloping algebras ``enveloping algebras'', and
denote them by  $ \U(\gerg) $;  similarly, we talk of ``restricted Lie
algebra'' even when  $ \, p = 0 \, $  simply meaning ``Lie algebra''.
                                            \par
   If a cocommutative connected Hopf algebra generated by its
primitive elements is also {\sl co-Poisson},  then the restricted
Lie algebra  $ \gerg $ such that  $ \, H = \U(\gerg) \, $  is indeed
a (restricted) Lie bialgebra.  Conversely, if a (restricted) Lie
algebra  $ \gerg $  is also a Lie bialgebra then  $ \U(\gerg) $
is a cocommutative connected  {\sl co-Poisson\/}  Hopf algebra
(cf.~[CP]).
                                      \par
   Let  $ V $  be a vector space: then the symmetric algebra  $ S(V) $
has a natural structure of Hopf algebra, given by  $ \, \Delta(x) =
x \otimes 1 + 1 \otimes x \, $,  $ \, \epsilon(x) = 0 \, $  and  $ \,
S(x) = -x \, $  for all  $ \, x \in V \, $.  If  $ \gerg $  is a Lie
algebra, then  $ S(\gerg) $  is also a Poisson Hopf algebra w.r.t.~the
Poisson bracket given by  $ \, {\{x,y\}}_{S(\gerg)} = {[x,y]}_\gerg
\, $  for all  $ \, x $,  $ y \in \gerg \, $.  If  $ \gerg $  is a
Lie coalgebra, then  $ S(\gerg) $  is also a co-Poisson Hopf algebra
w.r.t.~the Poisson cobracket determined by  $ \, \delta_{S(\gerg)}(x)
= \delta_\gerg(x) \, $  for all  $ \, x \in \gerg \, $.  Finally, if
$ \gerg $  is a Lie bialgebra, then  $ S(\gerg) $  is a bi-Poisson Hopf
algebra with respect to the previous Poisson bracket and cobracket
(cf.~[Tu1] and [Tu2] for details).

\vskip7pt

  {\bf 1.4 Filtrations.} \, Let  $ \; {\big\{F_z\big\}}_{z \in \Z}
=: \underline{F} : \big( \{0\} \subseteq \big) \cdots \subseteq F_{-1}
\subseteq F_0 \subseteq F_1 \subseteq \cdots \big( \subseteq V \big)
\; $  be a filtration of a vector space  $ V $.  We denote by  $ \,
G_{\underline{F}}(V) := \bigoplus_{z \in \Z} F_z \big/ F_{z-1} \, $
the associated graded vector space.  We say that  $ \, \underline{F}
\, $  is  {\sl exhaustive\/}  if  $ \; V^{\underline{F}} := \bigcup_{z
\in \Z} F_z = V \; $;  \; we say it is  {\sl separating\/}  if  $ \;
V_\downarrow := \bigcap_{z \in \Z} F_z = \{0\} \; $.  We say that
a filtered vector space is  {\sl exhausted\/}  if the filtration is
exhaustive; we say that it is  {\sl separated\/}  if the filtration
is separating.
%
%
                                            \par
  A filtration  $ \, \underline{F} = {\big\{ F_z
\big\}}_{z \in \Z} \, $  in an algebra  $ A $  is
said to be an  {\sl algebra
              filtration\/}  iff\break
                    $ \; m\big(F_\ell \otimes F_m\big) \subseteq
F_{\ell+m} \; $  for all  $ \, \ell $,  $ m $,  $ n \in \Z \, $.
Similarly, a filtration  $ \, \underline{F} = {\big\{ F_z \big\}}_{z
\in \Z} \, $  in a coalgebra  $ C $  is said to be a  {\sl coalgebra
filtration\/}  iff  $ \; \Delta\big(F_z\big) \subseteq \sum_{r+s=z}
F_r \otimes F_s \; $  for all  $ \, z \in \Z \, $.  Finally, a
filtration  $ \, \underline{F} = {\big\{ F_z \big\}}_{z \in \Z} \, $
in a bialgebra, or in a Hopf algebra,  $ H $  is said to be a  {\sl
bialgebra filtration\/},  or a  {\sl Hopf (algebra) filtration\/},
iff it is both an algebra and a coalgebra filtration and   --- in the
Hopf case ---   in addition  $ \; S\big(F_z\big) \subseteq F_z \; $
for all  $ \, z \in \Z \, $.  The notions of  {\sl exhausted\/}  and
{\sl separated\/}  for filtered algebras, coalgebras, bialgebras and
Hopf algebras are defined like for vector spaces, with respect to
the proper type of filtrations.

\vskip7pt

\proclaim{Lemma 1.5} \; Let  $ \, \underline{F} \, $  be an algebra
filtration of an algebra  $ A $.  Then  $ \, G_{\underline{F}}\,(A) \, $
\, is a graded algebra; if, in addition, it is commutative, then it is
a commutative graded  {\sl Poisson}  algebra.  If  $ E $  is another
algebra with algebra filtration  $ \, \underline{\varPhi} \, $  and
$ \, \phi : A \longrightarrow E \, $  is a morphism of algebras such
that  $ \, \phi(F_z) \subseteq \varPhi_z \, $  for all  $ \, z \in \Z
\, $,  \, then the morphism  $ \, G(\phi) : G_{\underline{F}}\,(A)
\longrightarrow G_{\underline{\varPhi}}\,(E) \, $  associated to
$ \phi $  is a morphism of graded algebras.  In addition, if  $ \,
G_{\underline{F}}\,(A) \, $  and  $ \, G_{\underline{\varPhi}}\,(E)
\, $  are commutative, then  $ \, G(\phi) \, $  is a morphism of
graded commutative  {\sl Poisson}  algebras.
                             \hfill\break
   \indent   The analogous statement holds replacing ``algebra'' with
``coalgebra'', ``commutative'' with ``cocommutative'' and ``Poisson''
with ``co-Poisson''.  In addition, if we start from bialgebras, or
Hopf algebras, with bialgebra filtrations, or Hopf filtrations, then
we end up with graded commutative Poisson bialgebras, or Poisson
Hopf algebras, and graded cocommutative co-Poisson bialgebras,
or co-Poisson Hopf algebras, respectively.
\endproclaim

\demo{Proof}  The only non-trivial part concerns the Poisson structure
on  $ \, G_{\underline{F}}(A) \, $  and the co-Poisson structure on
$ \, G_{\underline{F}}(C) \, $  (for a coalgebra  $ C \, $).  Indeed,
let  $ \, \underline{F} := {\big\{ F_z \big\}}_{z \in \Z} \, $  be
an algebra filtration of  $ A $.  The Poisson bracket on  $ \,
G_{\underline{F}}\,(A) \, $  is defined as follows.  For any  $ \,
\overline{x} \in F_z \big/ F_{z-1} \, $,  $ \, \overline{y} \in
F_\zeta \big/ F_{\zeta-1} \, $  ($ \, z $,  $ \zeta \in \Z \, $),
let  $ \, x \in F_z \, $,  \, resp.~$ \, y \in F_\zeta \, $,  \,
be a lift of  $ \overline{x} $,  resp.~of  $ \overline{y} \, $. 
Then  $ \, [x,y] := (xy-yx) \in F_{z + \zeta - 1} \, $  because
$ \, G_{\underline{F}}\,(A) \, $  is commutative; thus we set
$ \; \big\{ \overline{x}, \overline{y} \,\big\} := \overline{[x,y]}
\equiv [x,y] \mod F_{z + \zeta - 2} \, \in \, F_{z + \zeta - 1}
\big/ F_{z + \zeta - 2} \, $,
\; which is easily seen to define a Poisson bracket on  $ \,
G_{\underline{F}}\,(A) \, $  which makes it into a graded commutative
{\sl Poisson\/}  algebra.  Similarly, if  $ \, \underline{F} := {\big\{
F_z \big\}}_{z \in \Z} \, $  is a coalgebra filtration of  $ C $,  we
define a co-Poisson bracket on  $ \, G_{\underline{F}}\,(C) \, $  as
follows.  For any  $ \, \overline{x} \in F_z \big/ F_{z-1} \, $  ($ \,
z \in \Z \, $),  let  $ \, x \in F_z \, $  be a lift of  $ \overline{x}
\; $.  Then  $ \, \nabla(x) := \big( \Delta(x) - \Delta^{\text{op}}(x)
\big) \in \sum_{r+s=z-1} F_r \otimes F_s \, $  because  $ \,
G_{\underline{F}}\,(C) \, $  is cocommutative; thus  $ \;
\delta(\overline{x}\,) := \overline{\nabla(x)} \equiv \nabla(x)
\mod \hskip-3pt {\textstyle \sum\limits_{r+s=z-2}} \hskip-2pt
F_r \otimes F_s \, \in \hskip-4pt {\textstyle \sum\limits_{r+s=z-1}}
\hskip-1pt \big( F_r \big/ F_{r-1} \big) \otimes \big( F_s \big/
F_{s-1} \big) \; $  defines a Poisson cobracket which makes 
$ \, G_{\underline{F}}\,(C) \, $  into a graded cocommutative 
{\sl co-Poisson\/}  algebra.   \qed
\enddemo

\vskip7pt

\proclaim{Lemma 1.6} \, Let  $ C $  be a coalgebra.  If  $ \,
\underline{F} := {\big\{F_z \big\}}_{z \in \Z} \, $  is a coalgebra
filtration, then  $ \; C^{\underline{F}} := \bigcup\limits_{z \in \Z}
\! F_z \; $  is a coalgebra, which injects into  $ C $,  and  $ \;
C_{\underline{F}} := C \Big/ \hskip-4pt \bigcap\limits_{z \in \Z}
\hskip-3pt F_z \; $  is a coalgebra, which  $ C $  surjects onto.
The same holds for algebras, bialgebras, Hopf algebras, with algebra,
bialgebra, Hopf algebra filtrations respectively.
\endproclaim

\demo{Proof} The claim for  $ C^{\underline{F}} $  is trivial, while
for  $ C_{\underline{F}} \, $ we must prove  $ \, F_\downarrow := \!
\bigcap\limits_{z \in \Z} \! F_z \, $  is a coideal.
                                           \par
   Fix a basis  $ \, B_\downarrow \, $  of  $ F_\downarrow $  and a
basis  $ \, B_+ \, $  of any chosen complement of  $ F_\downarrow $,
so that  $ \, B := B_\downarrow \cup B_+ \, $  be a basis of  $ C $.
In addition, we can choose  $ B_+ $  to have the following property:
the span of  $ \, B_+ \cap \big( F_z \setminus F_{z-1} \big) \, $
has trivial intersection with  $ F_{z-1} \, $,  \, for all  $ \,
z \in \Z \, $.  Then  $ \, C \otimes F_\downarrow + F_\downarrow
\otimes C \, $  has basis  $ \, \big( B_\downarrow \otimes B_+ \big)
\cup \big( B_+ \otimes B_\downarrow \big) \cup \big( B_\downarrow
\otimes B_\downarrow \big) \, $.  Moreover, for each  $ \, b_+ \in
B_+ \, $  there is a unique  $ \, z(b_+) \in \Z \, $  such that
$ \, b_+ \in F_{z(b_+)} \setminus F_{z(b_+)-1} \; $.  Now pick
$ \, f \in F_\downarrow \, $,  \, and let  $ \, \Delta(f) =
\sum_{b, b' \in B} c_{b,b'} \cdot b \otimes b' \, $  be the
expansion of  $ \Delta(f) $  w.r.t.~the basis  $ \, B \otimes
B \, $  of  $ \, C \otimes C \, $:  \, then  $ \, \Delta(f)
\subseteq C \otimes F_\downarrow + F_\downarrow \otimes C \, $
if and only if  $ \, \big\{ (b,b') \in B_+ \times B_+ \;\big|\;
c_{b,b'} \not= 0 \,\big\} = \emptyset \, $.  So assume the latter
set is non empty, and let  $ \, \nu := \min \big\{\, z(b) + z(b')
\;\big|\; b, \, b' \in B_+ \, , \, c_{b,b'} \not= 0 \,\big\} \, $.
Then  $ \, \Delta(f) \not\in \! \sum\limits_{r+s=\nu-1} \hskip-5pt
F_r \otimes F_s \, $, \, which contradicts  $ \; \Delta(f) \in
\Delta \big( F_\downarrow \big) \subseteq \Delta \big( F_{\nu-1}
\big) \subseteq \! \sum\limits_{r+s=\nu-1} \hskip-5pt
F_r \otimes F_s \; $.
                                           \par
   As for algebras, the definition of algebra filtration implies that
all terms of  $ \underline{F} $  are ideals, so  $ F_\downarrow $  is
an ideal too, and we conclude.  Finally, the bialgebra case follows
from the previous two cases, and the Hopf algebra case follows too
once we note that in addition each term of  $ \underline{F} $  is
$ S $--stable  (by assumption), hence the same is true for
$ F_\downarrow $  as well.   \qed
\enddemo

 \vskip1,3truecm

\centerline {\bf \S \; 2 \  Connecting functors on
(co)augmented (co)algebras }

\vskip10pt

  {\bf 2.1 The  $ \undepsilon $--{\,}filtration  $ \, \underline{J}
\; $  on augmented algebras.}  \, Let  $ A $  be an augmented algebra
(cf.~\S 1.1).  Let  $ \, J := \hbox{\it Ker}\,(\undepsilon\,) \, $:
\, then  $ \; \underline{J} := {\big\{ J_{-n} := J^n , \, J_n :=
A \big\}}_{n \in \N} \; $  is clearly an algebra filtration of
$ A $,  which we call  {\sl the  $ \undepsilon $--filtration\/}
of  $ A \, $.  To simplify notation we shall usually forget
its (fixed) positive part.  We say that  $ A $  is  {\sl
$ \undepsilon $--separated\/}  if  $ \, \underline{J} \, $
is separating, that is  $ \, J^\infty := \bigcap_{n \in \N}
J^n = \{0\} \, $.  Next lemma points out some properties of
the  $ \undepsilon $--filtration  $ \underline{J} \, $:

\vskip7pt

\proclaim{Lemma 2.2}
                                   \hfill\break
   \indent   (a)  $ \; \underline{J} \, $  is an algebra filtration
of  $ \, A $,  \, which contains the radical filtration of  $ \, A $,
\, that is  $ \, J^n \supseteq {\hbox{\it Rad}\,(A)}^n \, $  for all
$ \, n \in \N \, $  where  $ \, \hbox{\it Rad}\,(A)\, $  is the
(Jacobson) radical of  $ \, A \, $.
                                   \hfill\break
   \indent   (b) \, If  $ \, A \, $  is  $ \undepsilon $--separated,
then it is i.p.-free.
                                   \hfill\break
   \indent   (c)  $ \; A^\vee := A \Big/ \bigcap_{n \in \N} J^n
\; $  is a quotient augmented algebra of  $ A \, $,  \, which is
$ \, \undepsilon $--separated.
\endproclaim

\demo{Proof} {\it (a)} \, By definition  $ \, \hbox{\it Rad}\,(A) \, $ 
is the intersection of all maximal left (or right) ideals of  $ A $, 
and $ J $  is one of them: so  $ \, J \supseteq \hbox{\it Rad}\,(A)
\, $,  \, whence the claim follows at once.
                                                  \par
   {\it (b)} \, Let  $ \, e \in A \, $  be idempotent, and
let  $ \, e_0 := \undepsilon(e) \, $,  $ \, e_+ := e - e_0
\cdot 1 \, $:  \, then  $ \, e_+ \in J \, $,  \, and  $ \,
{e_0}^{\!2} = \undepsilon\big(e^2\big) = \undepsilon\,(e) =
e_0 \, $,  \, i.e.~$ \, e_0 \in \Bbbk \, $  is idempotent,
whence  $ \, e_0 \in \{0,1\} \, $.  If  $ \, e_0 = 0 \, $ 
then  $ \, e_+ = e = e^n = {e_+}^{\!n} \in J^n \, $  for
all  $ \, n \in \N \, $  so  $ \, e_+ \in J^\infty :=
\bigcap_{n \in \N} \, $.  If  $ \, e_0 = 1 \, $  then  $ \,
{e_+}^{\!2} = e^2 - 2 \, e + 1 = -(e-1) = -e_+ \, $  whence 
$ \, e_+ = {(-1)}^{n+1} {e_+}^{n} \in J^n \, $  for all  $ \,
n \in \N \, $,  \, so  $ \, e_+ \in J^\infty \, $.  Therefore,
if  $ A $  is  $ \undepsilon $--separated  it is also i.p.-free.   
                                                     \par
   {\it (c)}  \, Lemma 1.6 proves that  $ \, A_{\underline{J}} =
A^\vee \, $  is a quotient algebra of  $ A \, $.  The augmentation
of  $ A $  induces an augmentation on  $ A^\vee $  too, and the
latter is  $ \undepsilon $--separated  by construction.   \qed  
\enddemo

\vskip7pt

\proclaim{Proposition 2.3} \, Mapping  $ \; A \mapsto A^\vee
:= A \big/ {\displaystyle \cap_{n \in \N}} J^n \; $  gives a
well-defined functor from the category  $ \calA^+ $  to the
subcategory of all  $ \undepsilon $--separated augmented algebras. 
Moreover, the augmented algebras  $ A $  of the latter subcategory
are characterized by  $ \, A^\vee = A \, $.   
\endproclaim

\demo{Proof}  By  Lemma 2.2{\it (c)},  the functor is well-defined
on objects, and we are left with defining the functor on morphisms.
The last part of the claim will be immediate.
                                       \par
   Let  $ \, \varphi : \, A \longrightarrow E \, $  be
a morphism of augmented algebras, i.e.~such that  $ \,
\undepsilon_{\scriptscriptstyle E} \circ \varphi =
\undepsilon_{\scriptscriptstyle A} \, $.  Then  $ \;
J_{\scriptscriptstyle A} = {\undepsilon_{\scriptscriptstyle
A}}^{\!-1} (0) = \varphi^{-1} \big( {\undepsilon_{\scriptscriptstyle
E}}^{\!-1} (0) \big) =  \varphi^{-1} \big( J_{\scriptscriptstyle E}
\big) \, $,  \, so  $ \, \varphi \big( {J_{\scriptscriptstyle
A}}^{\!n} \big) \subseteq {J_{\scriptscriptstyle E}}^{\!n}
\, $  for all  $ n \in \N \, $:  \, thus  $ \varphi \big(
{J_{\scriptscriptstyle A}}^{\!\infty} \big) \subseteq
{J_{\scriptscriptstyle E}}^{\!\infty} \, $,  \,
whence  $ \, \varphi \, $  induces a morphism of
$ \undepsilon $--separated  augmented algebras.   \qed
\enddemo

 \vskip3pt

   {\bf Remark 2.4:} \; It is worthwhile mentioning a special example
of  $ \undepsilon $--separated  augmented algebras, namely the  {\sl
graded\/}  ones (i.e.~those augmented algebras  $ A $  with an algebra
grading such that  $ \undepsilon_{\scriptscriptstyle A} $  is a morphism
of graded algebras w.r.t.~the trivial grading on the ground field 
$ \Bbbk $)  which are also  {\sl connected},  i.e.~their zero-degree
subspace is the  $ \Bbbk $--span  of  $ 1_{\scriptscriptstyle \! A}
\, $.  Then   
                                               \par
   {\it Every  {\sl graded connected}  augmented algebra  $ A $  is
{\sl  $ \epsilon $--separated},  or equivalently  $ \, A = A^\vee \, $.}
                                               \par
\noindent   
 Indeed, by definition, each non-zero homogeneous element in  $ \, J
:= J_{\scriptstyle A} \, $  has positive degree: thus any non-zero
homogeneous element of  $ J^n $  has degree at least  $ n $,  \,
so any non-zero homogeneous element of  $ J^\infty $  should have
degree at least  $ n $  {\sl for any}  $ \, n \in \N \, $.
Then  $ \, J^\infty = \{0\} \, $.

\vskip7pt

  {\bf 2.5 Drinfeld's  $ \delta_\bullet $--{\,}maps.} \, Let  $ C $
be a coaugmented coalgebra (cf.~\S 1.1).  For every  $ \, n \in \N
\, $,  \, define  $ \; \Delta^n \colon H \longrightarrow H^{\otimes n}
\; $  by  $ \, \Delta^0 := \epsilon \, $,  $ \, \Delta^1 :=
\id_{\scriptscriptstyle C} $,  \, and  $ \, \Delta^n := \big( \Delta
\otimes \id_{\scriptscriptstyle C}^{\,\otimes (n-2)} \big) \circ
\Delta^{n-1} \, $  if  $ \, n > 2 \, $.  For any ordered subset
$ \, \Phi = \{i_1, \dots, i_k\} \subseteq \{1, \dots, n\} \, $
with  $ \, i_1 < \dots < i_k \, $,  \, define the linear map
$ \, j_{\scriptscriptstyle \Phi} : H^{\otimes k} \longrightarrow
H^{\otimes n} \, $  by  $ \, j_{\scriptscriptstyle \Phi} (a_1 \otimes
\cdots \otimes a_k) := b_1 \otimes \cdots \otimes b_n \, $  with  $ \,
b_i := \und1 \, $  if  $ \, i \notin \Phi \, $  and  $ \, b_{i_m}
:= a_m \, $  for  $ \, 1 \leq m \leq k \, $.  Then set  $ \;
\Delta_\Phi := j_{\scriptscriptstyle \Phi} \circ \Delta^k \, $,
$ \, \Delta_\emptyset := \Delta^0 \, $,  and  $ \; \delta_\Phi
:= \sum_{\Psi \subset \Phi} {(-1)}^{n-|\Psi|} \Delta_\Psi \, $,
$ \; \delta_\emptyset := \epsilon \, $.  By the inclusion-exclusion
principle, the inverse for\-mula
$ \, \Delta_\Phi = \sum_{\Psi \subseteq \Phi} \delta_\Psi \, $
holds.
     \hbox{To be short we'll write  $ \, \delta_0 \! :=
\delta_\emptyset \, $  and  $ \, \delta_n := \! \delta_{\{1, 2,
\dots, n\}} \, $  too.}

\vskip7pt

\proclaim{Lemma 2.6}   
                                     \hfill\break
   \indent  (a) \quad \;\;  $ \delta_n \, = \, {\big(
\id_{\scriptscriptstyle C} - \underline{u} \circ \epsilon
\,\big)}^{\otimes n} \circ \Delta^n $  \qquad  for all
$ \; n \in \N_+ \; $;
                                     \hfill\break
   \indent  (b) \; The maps  $ \, \delta_n \, $  (and similarly the
$ \delta_\Phi $'s,  for all finite  $ \, \Phi \subseteq \N \, $)
are coassociative, i.e.
  $$  {} \hskip25pt  \big( \id_{\scriptscriptstyle C}^{\,\otimes s}
\otimes \delta_\ell \otimes \id_{\scriptscriptstyle C}^{\,\otimes
(n-1-s)} \big) \circ \delta_n \, = \, \delta_{n+\ell-1}   \eqno
\hbox{for all}  \quad  n, \ell, s \in \N \, ,  \;\, 0 \leq s
\leq n-1 \; .  $$
\endproclaim

\demo{Proof}  Claim  {\it (a)\/}  is proved by an easy induction,
and claim  {\it (b)\/}  follows from  {\it (a)\/}.   \qed
\enddemo

\vskip7pt

  {\bf 2.7 The  $ \delta_\bullet $--\,filtration  $ \, \underline{D}
\; $  on coaugmented coalgebras.}  \, Let  $ C $  be as above, and
take notations of \S 2.5.  For all  $ \, n \in \N \, $,  \, let
$ \; D_n := \text{\it Ker}\,(\delta_{n+1}) \, $:  \, then  $ \;
\underline{D} := {\big\{ D_{-n} := \{0\} , \, D_n \big\}}_{n \in \N}
\; $  is clearly a filtration of  $ C $,  which we call  {\sl the
$ \delta_\bullet $--filtration\/}  of  $ C \, $;  to simplify
notation we shall usually forget its (fixed) negative part.  We
say that  $ C $  is  {\sl  $ \delta_\bullet $--exhausted\/}
   \hbox{if  $ \, \underline{D} \, $  is exhaustive,
i.e.~$ \, \bigcup_{n \in \N} D_n = C \, $.}
                                           \par
   The  $ \delta_\bullet $--filtration  has several remarkable
properties: next lemma highlights some of them, in particular shows
that  $ \, \underline{D} \, $  is sort of a refinement of the
coradical filtration of  $ C $.  We make use of the notion of ``wedge''
product, namely  $ \, X \bigwedge Y := \Delta^{-1} \big( C \otimes Y
+ X \otimes C \big) \, $  for all subspaces  $ X $,  $ Y $  of  $ C $,
with  $ \, \bigwedge^1 \! X := X \, $  and  $ \, \bigwedge^{n+1} X
:= \big( \bigwedge^n X \big) \bigwedge X  \, $  for all  $ \, n \in
\N_+ \, $.

\vskip7pt

\proclaim{Lemma 2.8}
                                   \hfill\break
   \indent   (a)  $ \; D_0 = \Bbbk \cdot \und1 \, ,  \quad
D_n = \Delta^{-1} \big( C \otimes D_{n-1} + D_0 \otimes C \,\big)
= \bigwedge^{n+1} \! D_0  \quad  \text{for all}  \;\; n \in \N \, $.
                                   \hfill\break
   \indent   (b) $ \; \underline{D} \, $  is a coalgebra
filtration of  $ \, C $,  \, which is contained in the coradical
filtration of  $ \, C $,  \, that is  $ \, D_n \subseteq C_n
\, $  if  $ \; \underline{C} \, := {\big\{C_n\big\}}_{n \in \N}
\; $  is the coradical filtration of  $ \, C \, $.
                                   \hfill\break
   \indent   (c)  $ \; C \, $  is  $ \delta_\bullet $--exhausted
$ \quad \Longleftrightarrow \quad C \, $  is connected  $ \quad
\Longleftrightarrow \quad \underline{D} \, = \, \underline{C} \; $.
                                   \hfill\break
   \indent   (d)  $ \; C' := \bigcup_{n \in \N} D_n \; $  is
a subcoalgebra of  $ \, C \, $:  \, more precisely, it is
the irreducible (hence connected) component of  $ \, C $
containing  $ \und1 \, $.
\endproclaim

\demo{Proof}
   {\it (a)} \, Definitions give  $ \; D_0 := \text{\it Ker}\,
(\delta_1) = \Bbbk \cdot 1 \, $  and, using induction and
coassociativity,  $ \; D_n := \text{\it Ker}\,(\delta_{n+1}) =
%
%
\text{\it Ker}\,\big( (\delta_1 \otimes \delta_n) \circ
\Delta \big) =
%
%
\Delta^{-1} \big( C \otimes \text{\it Ker}\,
(\delta_n) + \text{\it Ker}\,(\delta_1) \otimes C \,\big) =
\Delta^{-1} \big( C \otimes D_{n-1} + D_0 \otimes C \,\big) \; $
for all  $ \, n \in \N \, $.  Since  $ \, D_0 := \bigwedge^1 \! D_0
\, $  we have by induction  $ \, D_n = \Delta^{-1} \big( C \otimes
D_{n-1} + D_0 \otimes C \,\big) = D_{n-1} \bigwedge D_0 = \big(
\bigwedge^n D_0 \big) \bigwedge D_0 = \bigwedge^{n+1} D_0 \, $
for all  $ \, n \in \N \, $.
                                                  \par
   {\it (b)} \, (cf.~[Ab], Theorem 2.4.1) From  {\it (a)\/}  we
have that  $ D_0 $  is a subcoalgebra, and then the  $ D_n $  are
subcoalgebras too, because the wedge of two subcoalgebras is again a
subcoalgebra.  In addition,  $ \, D_n = \bigwedge^{n+1} D_0 = \big(
\bigwedge^i D_0 \big) \bigwedge \big( \bigwedge^{n+1-i} D_0 \big)
= D_{i-1} \bigwedge D_{n-i} \, $  so  $ \, \Delta(D_n) \subseteq
C \otimes D_{n-i} + D_{i-1} \otimes C \, $  for  $ \, 1 \leq i \leq
n \, $,  \, and also for  $ \, i \in \{0,n+1\} \, $  because  $ D_n $
is a subcoalgebra.  Then  $ \, \Delta(D_n) \subseteq \sum_{r+s=n}
D_r \otimes D_s \, $  (for all  $ n \, $)  follows from the fact
(cf.~[Sw], \S 9.1.5) that for any vector space  $ V $  and any
filtration  $ \, {\big\{ V_{-\ell} := 0, V_\ell \big\}}_{\ell \in
\N} \, $  of  $ V $  one has  $ \; \bigcap_{i=0}^{\ell+1} \big( V
\otimes V_{\ell+1-i} + V_i \otimes V \,\big) \, = \, \sum_{i=1}^{\ell
+ 1} \, V_i \otimes V_{\ell+2-1} \, $  (set  $ \, V = C \, $  and
$ \, V_\ell = D_{\ell-1} \, $  for all  $ \ell \, $).
                                                     \par
   Second we prove that  $ \, \underline{D} \, $  is contained in
$ \, \underline{C} \, $.  Recall that  $ \, \underline{C} \, $  may
be defined inductively by  $ \; C_0 := \hbox{\sl coRad}\,(C) \, $,
$ \, C_n := \Delta^{-1} \big( C \otimes C_{n-1} + C_0 \otimes C
\,\big) = C_0 \bigwedge C_{n-1} = \bigwedge^{n+1} C_0 \, $;  \;
then we get  $ \, D_n \subseteq C_n \, $  for all  $ \, n \in
\N \, $  by induction using  $ \, C_0 \supseteq \Bbbk \cdot 1
= D_0 \, $,  \, by  {\it (a)}.
                                                     \par
   {\it (c)}  $ \, C \, $  is connected iff  $ \, C_0 = \Bbbk \cdot
1 \, $.  But in this case  $ \, C_0 = D_0 \, $,  \, hence  $ \; C_n
= \bigwedge^{n+1} C_0 = \bigwedge^{n+1} D_0 = D_n \; $  for all
$ n \, $  (by  {\it (a)\/}),  that is  $ \, \underline{D} =
\underline{C} \, $.  Conversely, if  $ \, \underline{D} =
\underline{C} \, $  then  $ \, \text{\sl coRad}\,(C) = C_0 =
D_0 = \Bbbk \cdot 1 \, $  thus  $ C $  is connected.  Further,
if  $ \, \underline{D} = \underline{C} \, $  then  $ C $  is
$ \delta_\bullet $--exhausted,  because the coradical filtration
$ \underline{C} $  is always exhaustive.  Conversely, assume
$ C $  is  $ \delta_\bullet $--exhausted:  then by a standard
result (cf.~[Ab], Theorem 2.3.9(ii), or [Mo], Lemma 5.3.4) we
have  $ \, D_0 \supseteq C_0 \, $,  \, so  $ \, D_0 = C_0 \, $
by  {\it (a)},  and we conclude by induction as above that
$ \, \underline{D} = \underline{C} \, $.
                                                     \par
   {\it (d)}  $ \, C' := \bigcup_{n \in \N} D_n \, $  is a
subcoalgebra of  $ \, C \, $  because  $ \, \underline{D} \, $
is a coalgebra filtration, so  Lemma 1.6 applies.  Moreover, by
construction  $ \, \und1 \in D_0 \subseteq C' \, $,  \, and
$ C' $  is  $ \delta_\bullet $--exhausted  (w.r.t.~the same
group-like element  $ \und1 \, $),  hence it is connected by
{\it (c)}.  Thus  $ C' $  is the irreducible (connected)
component of  $ C $  containing  $ \und1 \, $.   \qed
\enddemo  

\vskip7pt

\proclaim{Proposition 2.9} \, Mapping  $ \; C \mapsto C' :=
\bigcup_{n \in \N} D_n \; $  gives a well-defined functor from 
$ \calC^+ $  to the subcategory of all  $ \delta_\bullet $--exhausted 
(=connected) coaugmented coalgebras.  Moreover, the coaugmented
coalgebras  $ C $  of the latter subcategory are characterized
by  $ \, C' = C \, $.   
\endproclaim

\demo{Proof}  By Lemma 2.8{\it (d)},  if  $ C $  is
a coaugmented coalgebra, then  $ \, C' \, $  with its
$ \delta_\bullet $--filtration  (w.r.t.~the same group-like
element  $ \und1 \, $)  is an exhausted coaugmented coalgebra:
in other words it is  $ \delta_\bullet $--exhausted,  so by
Lemma 2.8{\it (c)\/}  it is connected too.  The last part of
the claim being immediate, we are left to define the functor
on morphisms.
                                       \par
   Let  $ \, \varphi : \, C \longrightarrow K \, $  be a morphism of
coaugmented coalgebras (mapping  $ \, \und1 \in G(C) \, $  to  $ \,
\und1 \in G(K) \, $).  Then  $ \; \delta_m \circ \varphi = {\big(
\id - \underline{u} \circ \epsilon \,\big)}^{\otimes m} \! \circ
\Delta^m \circ \varphi = \varphi^{\otimes m} \circ {\big( \id -
\underline{u} \circ \epsilon \,\big)}^{\otimes m} \! \circ \Delta^m
= \varphi^{\otimes m} \circ \delta_m \; $  for all  $ m \in \N \, $, 
\, which yields  $ \, \varphi\big(D_n(C)\big) = \varphi \big( \hbox{\it
Ker}\,\big(\delta_{n+1}{\big|}_C\big) \big) \subseteq \hbox{\it Ker}\,
\big(\delta_{n+1}{\big|}_K\big) = D_n(K) \, $.  Thus  $ \, \varphi \, $ 
maps  $ \, \underline{D}(C) \, $  to  $ \, \underline{D}(K) \, $  and
so  $ \, C' \, $  to  $ \, K' \, $,  \, so the restriction  $ \,
\varphi{\big\vert}_{C'} \, $  of  $ \, \varphi \, $  to  $ \, C' \, $ 
is a morphism of exhausted coaugmented coalgebras which we take as 
$ \, \varphi' \in \hbox{\it Mor}\,\big(C',K'\big) \, $.   \qed  
\enddemo 

\vskip7pt

   {\bf Remark 2.10:} \; Lemma 2.8{\it (c)\/}  yields alternative
proofs of two well-known facts.   
                                               \par
   {\it (a) \, Every  {\sl graded connected}  (in the graded sense,
i.e.~$ \dim(C_0) = 1 $)  coaugmented coalgebra  $ C $  is  {\sl
connected}  (in the mere coalgebra sense).}  \; Indeed, given 
$ \, c \in C_{\partial(c)} \, $  homogeneous of degree  $ \,
\partial(c) \, $,  \, one has  $ \, \delta_n(c) \in {C_+}^{\hskip-3pt
\otimes n} \bigcap \big( \sum_{r+s=\partial(c)} C_r \otimes C_s \big)
\, $  for all  $ \, n \, $.  Thus  $ \, c \in \hbox{\it Ker}\,
(\delta_{\partial(c)+1}) =: D_{\partial(c)} \, $,  \, hence  $ \,
C = \bigoplus_{n \in \N} C_n \subseteq \bigcup_{n \in \N} D_n =
C' \subseteq C \, $  so  $ \, C = C' \, $.   
                                               \par
   {\it (b) \, Every  {\sl connected}  coaugmented coalgebra is  {\sl
monic}.}  \; By  Lemma 2.8{\it (c)--(d)\/}  it is enough to prove  $ \,
G(C') = \{ \und1 \,\} \, $.  As  $ \, G(C') = G(C) \cap C' \, $,  \,
let $ \, g \in G(C) \, $:  \, then  $ \, \delta_n(g) = {( g - \und1
\,)}^{\otimes n} \, $  for all  $ \, n \in \N \, $,  \, so  $ \, g
\in C' := \bigcup_{n \in \N} D_n = \bigcup_{n \in \N} \hbox{\it
Ker}\,(\delta_{n+1}) \, $  if and only if  $ \, g = \und1 \, $.

\vskip7pt

   {\bf 2.11 Connecting functors and Hopf duality.} \, Let's start
from an augmented algebra  $ A $,  with  $ \, J := \hbox{\it Ker}\,
(\undepsilon\,) \, $.  Its  $ \undepsilon $--filtration  has
$ \; J^n \, := \,\text{\sl Im}\Big( J^{\otimes n}
\,{\lhook\joinrel\longrightarrow}\, H^{\otimes n} \,{\buildrel
{\;\mu^n} \over {\relbar\joinrel\relbar\joinrel\twoheadrightarrow}}\,
H \Big) \; $  as  $ n $--th  term,  the left hand side arrow being the
natural embedding induced by  $ \, J \,{\lhook\joinrel\longrightarrow}\,
H \, $  and  $ \, {\mu^n} \, $  being the  $ n $--fold  iterated
multiplication of  $ H $.  Similarly, let  $ C $  be a
coaugmented coalgebra.  The  $ s $--th  term of its
$ \delta_\bullet $--filtration  is  $ \; D_s \, := \,
\text{\sl Ker}\Big( H \,{\buildrel {\;\Delta^{s+1}} \over
{\lhook\joinrel\relbar\joinrel\relbar\joinrel\longrightarrow}}\,
H^{\otimes (s+1)} \,{\relbar\joinrel\twoheadrightarrow}\,
J^{\otimes (s+1)} \Big) \; $  where the right hand side arrow is the
natural projection induced by  $ \, \big( \id_{\scriptscriptstyle C}
- \underline{u} \circ \epsilon \big) \, \colon \, H
\,{\relbar\joinrel\relbar\joinrel\twoheadrightarrow}\,
J \, $  and  $ \, \Delta^{s+1} \, $  is the  $ (s+1) $--fold
iterated comultiplication of  $ C $.  This clearly means that
 \vskip3pt
   {\it (a)}  \,  {\sl The notions of  $ \undepsilon $--filtration 
and of  $ \delta_\bullet $--filtration  are dual to each other.}
 \vskip3pt
\noindent
   Similarly, as taking intersection and taking unions are mutually
dual operations, and taking submodules and quotient modules are
mutually dual too, we have the following two facts:
 \vskip3pt
   {\it (b)}  \,  {\sl The notions of  $ \; A^\vee := A \Big/
\bigcap_{n \in \N} J^n \; $  and of  $ \; C' := \bigcup_{n \in \N}
D_n \; $  are dual to each other;}
 \vskip3pt
   {\it (c)}  \,  {\sl The notions of  $ \undepsilon $--separated 
(for an augmented algebra) and  $ \delta_\bullet $--exhausted 
(for a coaugmented coalgebra) are dual to each other.}

 \vskip7pt

   {\bf Remark 2.12:}  {\it (a)  $ (\ )^\vee $  and  $ (\ )' $  as
``{\sl connecting}  functors''.} \, We now explain in what sense
both  $ A^\vee $  and  $ C' $  are ``connected'' objects.  Indeed,
$ C' $  is truly connected, in the sense of coalgebra theory
(cf.~Lemma 2.8{\it (d)\/}).  By duality, we might expect that
$ A^\vee $  be (or correspond to) a ``connected'' object too.
In fact, when  $ A $  is commutative, it is the algebra of
regular functions  $ F[\Cal{V}] $  of an algebraic variety
$ \Cal{V} \, $;  the augmentation on  $ A $  is a character,
hence corresponds to the choice of a point  $ \, P_0 \in \Cal{V}
\, $:  thus  $ A $  does correspond to the  {\sl pointed variety}
$ \, (\Cal{V},P_0) \, $.  Then  $ \, A^\vee = F[\Cal{V}_0] \, $
where  $ \Cal{V}_0 $  is the connected component of  $ \Cal{V} $
containing  $ P_0 \, $.  This follows at once because  $ A^\vee $
is i.p.-free (cf.~Lemma 2.2).  More in general, for any  $ \, A
\in \calA^+ $,  if  $ A^\vee $  is commutative then its spectrum
is a  {\sl connected\/}  algebraic variety.  Therefore we shall call
``connecting functors'' both functors  $ \, A \mapsto A^\vee \, $
and  $ \, C \mapsto C' \, $.
                                           \par
   {\it (b) Asymmetry of connecting functors on bialgebras.} \,
Let  $ B $  be a bialgebra.  As the notion of  $ B^\vee $  is dual
to that of  $ B' $,  and since  $ \, B = B' \, $  implies that  $ B $
is a  {\sl Hopf\/}  algebra  (cf.~Corollary 3.4{\it (b)\/}),  one might
dually conjecture that  $ \, B = B^\vee \, $  imply that  $ B $  is a
{\sl Hopf\/}  algebra.  Actually,  {\sl this is false},  the bialgebra
$ \, B := F \big[ \text{\it Mat}\,(n,\Bbbk) \big] \, $  yielding a
counter\-example:  $ \, F \big[ \text{\it Mat}\,(n,\Bbbk) \big] =
{F \big[ \text{\it Mat}\,(n,\Bbbk) \big]}^\vee $,  \, and yet
$ \, F \big[ \text{\it Mat}\,(n,\Bbbk) \big] \, $  is  {\sl
not\/}  a Hopf algebra.
                                           \par
   {\it (c) Hopf duality and augmented pairings.} \,  The most
precise description of the relationship between connecting functors
of the two types uses the notion of ``augmented pairing'':

\vskip7pt

\proclaim{Definition 2.13} \, Let  $ \, A \in \calA^+ \, $,  $ \,
C \in \calC^+ \, $.  We call  {\sl augmented pairing}  between
$ A $  and  $ C $  any bilinear mapping  $ \; \langle \,\ ,
\ \rangle \, \colon \, A \times C \loongrightarrow \Bbbk \; $
such that, for all  $ \, x, x_1, x_2 \in A $  and  $ \, y \in C $,
$ \; \big\langle x_1 \cdot x_2 \, , y \big\rangle = \big\langle x_1
\otimes x_2 \, , \Delta(y) \big\rangle := \sum_{(y)} \big\langle x_1
\, , y_{(1)} \big\rangle \cdot \big\langle x_2 \, , y_{(2)} \big\rangle
\, $,  $ \, \langle 1 \, , y \rangle \, = \, \epsilon(y) \, $,  $ \,
\langle x \, , \und1 \, \rangle \, = \, \undepsilon(x) \, $.  For any
$ \, B, P \in \calB \, $  we call  {\sl bialgebra pairing}  between
$ B $  and  $ P $  any augmented pairing between the latters (see \S
1.1) such that we have also, symmetrically,  $ \; \big\langle x \, ,
y_1 \cdot y_2 \big\rangle \, = \, \big\langle \Delta(x) \, , y_1
\otimes y_2 \big\rangle \, := \, \sum_{(x)} \big\langle x_{(1)}
\, , y_1 \big\rangle \cdot \big\langle x_{(2)} \, , y_2 \big\rangle
\, $,  \, for all  $ \, x \in B $,  $ \, y_1, y_2 \in P \, $.  For
any  $ \, H, K \in \HA \, $  we call  {\sl Hopf algebra pairing} 
(or  {\sl Hopf pairing})  between  $ H $  and  $ K $  any bialgebra
pairing such that, in addition,  $ \; \big\langle S(x), y \big\rangle
\, = \, \big\langle x, S(y) \big\rangle \; $  for all  $ \, x \in H
\, $,  $ \, y \in K \, $.
                                       \hfill\break
   \indent   We say that a pairing as above is  {\sl perfect on
the left (right)}  if its left (right) kernel is trivial; we say
it is  {\sl perfect}  if it is both left and right perfect.
\endproclaim

\vskip7pt

\proclaim{Theorem 2.14} \, Let  $ \, A \in \calA^+ $,  $ C \in
\calC^+ $  and let  $ \, \pi \colon \, A \times C \loongrightarrow
\Bbbk \, $  be an augmented pairing.  Then  $ \pi $  induces a  {\sl
filtered}  augmented pairing  $ \, \pi_f \, \colon A^\vee \times
C' \loongrightarrow \Bbbk \, $  and a  {\sl graded}  augmented
pairing  $ \, \pi_{{}_G} \, \colon \, G_{\underline{J}}(A) \times
G_{\underline{D}}(C) \loongrightarrow \Bbbk \, $  (notation of \S
1.4), both perfect on the right.  If in addition  $ \, \pi $  is
perfect then  $ \, \pi_f $  and  $ \, \pi_{{}_G} $  are perfect
as well.
\endproclaim

\demo{Proof} \, Consider the filtrations  $ \, \underline{J} =
{\big\{ {J_{\!{}_A}}^{\! n} \big\}}_{n \in \N} \, $  and  $ \,
\underline{D} = {\big\{ D^{\scriptscriptstyle C}_n \big\}}_{n
\in \N} \, $.  The key fact is that
  $$  D^{\scriptscriptstyle C}_n = \big( {J_{\!{}_A}}^{\! n+1}
\big)^\perp  \qquad  \hbox{and}  \qquad  {J_{\!{}_A}}^{\! n+1}
\subseteq \big( D^{\scriptscriptstyle C}_n \big)^\perp  \qquad
\qquad  \hbox{for all} \;\; n \in \N \, .   \eqno (2.1)  $$
   \indent   Indeed, let  $ E $  be an algebra,  $ K $  a coalgebra, and
assume there is a bilinear pairing  $ \, E \times K \!\longrightarrow
\Bbbk \, $  enjoying  $ \; \big\langle x_1 \cdot x_2, y \big\rangle \,
= \, \big\langle x_1 \otimes x_2, \Delta(y) \big\rangle \; $  (for
all  $ \, x_1, x_2 \in E $  and  $ \, y \in K \, $).  If  $ X $  is
a subspace of  $ K $,  then  $ \, \bigwedge^n X = {\big( {(X^\perp)}^n
\big)}^\perp \, $  (cf.~\S 2.7) for all  $ \, n \in \N \, $,  \, where
the superscript  $ \perp $  means ``orthogonal subspace'' (either in
$ E $  or in  $ K \, $)  w.r.t.~the pairing under exam (cf.~[Ab] or
[Mo]).  Now,  Lemma 2.8{\it (a)}  gives  $ \, D^{\scriptscriptstyle
C}_n = \bigwedge^{n+1} (\Bbbk \!\cdot\! \und1 {}_{\,{}_C}) \, $,  \,
thus  $ \, D^{\scriptscriptstyle C}_n = \bigwedge^{n+1} (\Bbbk
\!\cdot\! \und1 {}_{\,{}_C}) = \Big( \big( {(\Bbbk \!\cdot\! \und1
{}_{\,{}_C})}^\perp \big)^{n+1} \Big)^\perp = {\big( {J_{{}_A}}^{\!n+1}
\big)}^\perp \, $  because  $ \, {(\Bbbk \!\cdot\! \und1
{}_{\,{}_C})}^\perp = J_{\!{}_A} \, $  (w.r.t.~the pairing
$ \pi $  above).  Therefore  $ \, D^{\scriptscriptstyle C}_n
= \big( {J_{\!{}_A}}^{\! n+1} \big)^\perp \, $,  \, and this
also implies  $ \, {J_{\!{}_A}}^{\! n+1} \subseteq \big(
D^{\scriptscriptstyle C}_n \big)^\perp \, $.
                                         \par
   Now  $ \, C' := \bigcup_{n \in \N} D^{\scriptscriptstyle C}_n =
\bigcup_{n \in \N} \left( {J_{\!{}_A}}^{\! n+1} \right)^\perp \! =
\left( \, \bigcap_{n \in \N} {J_{\!{}_A}}^{\! n+1} \right)^\perp \!
= \big( {J_{\!{}_A}}^{\!\infty} \big)^\perp \, $.  So  $ \pi $
induces a Hopf pairing  $ \, \pi_f \, \colon \, A^\vee \times C'
\loongrightarrow \Bbbk \, $  as required, and by (2.1) this respects
the filtrations on either side.  Then by general theory  $ \pi_f $
induces a graded Hopf pairing  $ \pi_{{}_G} $  as required: in
particular  $ \pi_{{}_G} $  is well-defined because  $ \,
D^{\scriptscriptstyle C}_n \subseteq \big( {J_{\!{}_A}}^{\! n+1}
\big)^\perp \, $  and  $ \, {J_{\!{}_A}}^{\! n+1} \subseteq \big(
D^{\scriptscriptstyle C}_n \big)^\perp \, $  ($ \, n \in \N_+ \, $)
by (2.1).  Moreover, both  $ \pi_f $  and  $ \pi_{{}_G} $  are perfect
on the right because all the inclusions  $ \, D^{\scriptscriptstyle C}_n
\subseteq \big( {J_{\!{}_A}}^{\! n+1} \big)^\perp \, $  happen to be
identities.  Finally if  $ \pi_f $  is perfect it is easy to see that
$ \pi_{{}_G} $  is perfect as well; note that this improves (2.1), for
then  $ \, {J_{\!{}_A}}^{\! n+1} = \big( D^{\scriptscriptstyle C}_n
\big)^\perp \, $  for all  $ \, n \in \N \, $.   \qed
\enddemo

\vskip7pt

{\bf 2.15 The crystal functors.} \, It is clear from the very
construction that mapping  $ \, A \mapsto \widehat{A} :=
G_{\underline{J}}(A) \, $  (for all  $ \, A \in \calA^+ \, $)
defines a functor from  $ \calA^+ $  to the category of graded,
augmented  $ \Bbbk $--algebras,  which factors through the functor 
$ \, A \mapsto A^\vee $.  Similarly, mapping  $ \, C \mapsto
\widetilde{C} := G_{\underline{D}}(C) \, $  (for all  $ \, C \in
\calC^+ \, $)  defines a functor from  $ \calC^+ $  to the category
of graded, coaugmented  $ \Bbbk $--coalgebras,  which factors through
the functor  $ \, C \mapsto C' \, $.   
                                                      \par
   In \S 4 we'll show that each of these functors can be seen as
a ``crystallization process'', in the sense, loosely speaking, of
Kashiwara's terminology of ``crystal basis'' for quantum groups:
we move from one fiber to another, very special one, within a
1-parameter deformation.  Therefore, we call them  {\sl ``crystal
functors''}.

 \vskip1,3truecm

\centerline {\bf \S \, 3 Connecting and crystal functors
on bialgebras and Hopf algebras }

\vskip10pt

{\bf 3.1 The program.} \, We apply now the connecting functors to
bialgebras and Hopf algebras.  Then we look at the graded objects
associated to the filtrations  $ \underline{J} $  and  $ \underline{D} $ 
in a bialgebra: this leads to the crystal functors on  $ \Cal{B} $  and 
$ \HA \, $,  the main topic of this section.
                                          \par
   From now on, every bialgebra  $ B $  will be considered as
a coaugmented coalgebra w.r.t.~its unit map, hence w.r.t~the
group-like element  $ 1 $  (the unit of  $ B \, $),  and the
corresponding maps  $ \, \delta_n \, $  ($ \, n \in \N \, $)
and  $ \delta_\bullet $--filtration  $ \, \underline{D} \, $
will be taken into account.  Similarly,  $ B $  will be
considered as an augmented algebra w.r.t.~the special
algebra morphism  $ \, \undepsilon = \epsilon \, $  (the counit
of  $ B \, $),  and the corresponding  $ \undepsilon $--filtration
(also called  $ \epsilon $--filtration)  $ \, \underline{J} \, $
will be considered.
                                          \par
   We begin with a technical result about the ``multiplicative''
properties of the maps  $ \, \delta_n \, $.

\vskip7pt

\proclaim{Lemma 3.2} ([KT], Lemma 3.2)  Let  $ \, B \in \calB \, $,
$ \, a $, $ b \in B $,  and  $ \, \Phi \! \subseteq \! \N \, $,  with
$ \Phi $  finite.  Then
\vskip2pt
   \indent  (a) \quad \;\;  $ \delta_\Phi(ab) \; = \; \sum_{\Lambda
\cup Y = \Phi} \delta_\Lambda(a) \, \delta_Y(b) \;\; $;
                                     \hfill\break
\vskip-11pt
   \indent  (b) \quad  if  $ \, \Phi \not= \emptyset \, $,  \;
then  \quad  $ \delta_\Phi(ab - ba) \; = \; \sum_{\Sb  \Lambda
\cup Y = \Phi  \\
\Lambda \cap Y \not= \emptyset  \endSb} \,
\big( \delta_\Lambda(a) \, \delta_Y(b) -
\delta_Y(b) \, \delta_\Lambda(a) \big) \;\; $.  \qed
\endproclaim

\vskip7pt

\proclaim{Lemma 3.3} \, Let  $ B $  be a bialgebra.  Then
$ \, \underline{J} \, $  and  $ \, \underline{D} \, $  are
bialgebra filtrations.  If  $ \, B $  is also a Hopf algebra,
then  $ \, \underline{J} \, $  and  $ \, \underline{D} \, $
are Hopf algebra filtrations.
\endproclaim

\demo{Proof} \, By  Lemma 2.2{\it (a)},  $ \underline{J} $  is an
algebra filtration.  In addition, since  $ \, J := \hbox{\it Ker}\,
(\epsilon) \, $  is a biideal we have  $ \; \Delta \big( J^n \big) =
{\Delta(J)}^n \subseteq {\big( B \otimes J + J \otimes B \,\big)}^n
\subseteq \sum_{r+s=n} J^r \otimes J^s \; $  for all  $ \, n \in \N
\, $,  \, so that  $ \underline{J} $  is a bialgebra filtration as
well.  Moreover, in the Hopf algebra case  $ J $  is a Hopf ideal,
hence  $ \, S\big(J^n\big) \subseteq J^n \, $  for all  $ \, n \in
\N \, $,  \, which means  $ \underline{J} $  is a Hopf algebra
filtration.
                                              \par
   As for  $ \underline{D} \, $,  by  Lemma 2.8{\it (b)\/}  it
is a coalgebra filtration.  The fact that it is also a bialgebra
filtration is a direct consequence of  Lemma 3.2{\it (a)}.  The
claim in the Hopf case follows noting that  $ \, \delta_n
\circ S = S^{\otimes n} \circ \delta_n \, $  so  $ \, S(D_n) =
S\big(\hbox{\it Ker}\,(\delta_{n+1})\big) \subseteq \hbox{\it
Ker}\,(\delta_{n+1}) = D_n \, $  for all  $ \, n \in \N \, $.
Otherwise, the claim for  $ \underline{D} $  follows from [Ab],
Theorem 2.4.1, or [Mo], 5.2.8.   \qed
\enddemo

\vskip7pt

\proclaim{Corollary 3.4} \, Let  $ B $  be a bialgebra.  Then
                                    \hfill\break
   \indent   (a)  $ \, B^\vee \! := B  \hskip-0.3pt  \Big/
\hskip-0.3pt  \bigcap_{n \in \N} \! J^n \, $  is an
$ \epsilon $--separated  (i.p.-free) bialgebra,
                      \hbox{which  $ B $  surjects onto.}
   \indent   (b)  $ \, B' := \bigcup_{n \in \N} D_n \, $  is
a  $ \delta_\bullet $--exhausted  (connected) Hopf algebra, which
injects into  $ B \, $:  \, more precisely, it is the irreducible
(actually, connected) component of  $ \, B $  containing\/  $ 1 $.
                                    \hfill\break
   \indent   (c) \, If in addition  $ \, B = H \, $  is a Hopf
algebra, then  $ H^\vee $  is a Hopf algebra quotient of  $ \, H $
and  $ H' $  is a Hopf subalgebra of  $ \, H $.
\endproclaim

\demo{Proof}  Claim  {\it (a)\/}  follows from the results in \S\S 2--3.
As for  {\it (b)},  the results of \S\S 2--3 imply that  $ B' $  is a
$ \delta_\bullet $--exhausted  (connected) bialgebra, embedding into
$ B $,  namely the irreducible component of  $ B $  containing  $ 1 $. 
In addition, a general result ensures that any connected bialgebra is
a Hopf algebra (cf.~[Ab], Theorem 2.4.24), whence this holds true
for  $ B' $.  Finally, claim  {\it (c)\/}  follows from {\it (a)},
{\it (b)},  and Lemma 3.3.   \qed
\enddemo

\vskip7pt

\proclaim{Theorem 3.5} \, Let  $ B $  be a bialgebra,
$ \underline{J} \, $,  $ \underline{D} $  its
$ \epsilon $--filtration  and
$ \delta_\bullet $--filtration  respectively.
                                    \hfill\break
   \indent   (a)  $ \; \widehat{B} := G_{\underline{J}}\,(B) \, $
\, is a graded
%
%
%
cocommutative
co-Poisson Hopf algebra generated by  $ \, P\big(G_{\underline{J}}\,
(B)\big) \, $,  \, the set of its primitive elements.  Therefore
$ \, \widehat{B} \cong \U(\gerg_-) \, $  as graded co-Poisson Hopf
algebras, for some restricted Lie bialgebra  $ \gerg_- $  which is
graded as a Lie algebra.  In particular, if  $ \; p = 0 \, $  and
$ \, \dim(B) \!\in \N \, $  then  $ \, \widehat{B} = \Bbbk \!
\cdot \! 1 \, $  and  $ \, \gerg_- = \{0\} \, $.
                                    \hfill\break
   \indent   (b)  $ \; \widetilde{B} := G_{\underline{D}}\,(B) \, $
\, is a graded
%
%
%
%
commutative Poisson Hopf algebra.  Therefore,  $ \, \widetilde{B}
\cong F[G_+] \, $  for some connected algebraic Poisson group  $ G_+ $
which, as a variety, is a (pro)affine space.  If  $ \; p = 0 \, $
then  $ \, \widetilde{B} \cong F[G_+] \, $  is a polynomial algebra,
i.e.~$ \, F[G_+] = \Bbbk \big[ {\{x_i \}}_{i \in \Cal{I}} \big] \, $
(for some set  $ \Cal{I} $);  in particular, if  $ \, \dim(B) \in \N
\, $  then  $ \, \widetilde{B} = \Bbbk \cdot 1 \, $  and  $ \, G_+ =
\{1\} \, $.  If  $ \, p > 0 \, $  then  $ G_+ $  has dimension 0 and
height 1, and if\/  $ \Bbbk $  is perfect then  $ \, \widetilde{B} \cong
F[G_+] \, $  is a truncated polynomial algebra,  i.e.~$ \, F[G_+] =
\Bbbk \big[ {\{x_i\}}_{i \in \Cal{I}} \big] \Big/ \big( {\{ x_i^{\,p}
\}}_{i \in \Cal{I}} \big) \, $  (for some set  $ \Cal{I} $).
\endproclaim

\demo{Proof} {\it (a)} \, Thanks to  Corollary 3.4{\it (a)},
I must only prove that  $ \, G_{\underline{J}}\,(B) \, $  is
cocommutative, and then apply Lemma 1.5 to have the co-Poisson
structure.  In fact, even for this, it is enough to show that 
$ \, G_{\underline{J}}\,(B) \, $  is generated by  $ \, P \big(
G_{\underline{J}}\,(B) \big) \, $.  Moreover, I must show
that  $ G_{\underline{J}}\,(B) $  is in fact a (graded
$ \epsilon $--separated cocommutative co-Poisson)
{\sl Hopf algebra},  not only a bialgebra!  The
last part of the statement then will follow from
\S 1.3.
                                               \par
   Since  $ \, G_{\underline{J}}\,(B) \, $  is generated by
$ \, J \big/ J^2 \, $,  \, it is enough to show that  $ \, J
\big/ J^2 \subseteq P\big(G_{\underline{J}}\,(B)\big) \, $.
Let  $ \, 0 \not= \overline{\eta} \in J \big/ J^2 \, $,  \, and
let  $ \, \eta \in J \, $  be a lift of  $ \, \overline{\eta} \, $:
\, then  $ \, \Delta(\eta) = \epsilon(\eta) \cdot 1 \otimes 1 +
\delta_1(\eta) \otimes 1 + 1 \otimes \delta_1(\eta) + \delta_2(\eta)
= \eta \otimes 1 + 1 \otimes \eta + \delta_2(\eta) \, $  (by the
very definitions and by  $ \, \eta \in J \, $).  Now,  $ \,
\delta_2(\eta) \in J \otimes J \, $,  \, hence  $ \,
\overline{\delta_2(\eta)} = \overline{0} \in {\big(
G_{\underline{J}}\,(B) \otimes G_{\underline{J}}\,(B)
\big)}_1 = \big( B \big/ J \big) \otimes J \big/ J^2
+ J \big/ J^2 \otimes \big( B \big/ J \big) \, $.  So
$ \, \Delta(\overline{\eta}\,) := \overline{\Delta(\eta)} =
\overline{\eta \otimes 1 + 1 \otimes \eta} = \overline{\eta}
\otimes 1 + 1 \otimes \overline{\eta} \, $,  \, which gives
$ \, J \big/ J^2 \subseteq P\big(G_{\underline{J}}\,(B)\big)
\, $.
                                               \par
   Finally,  $ \, G_{\underline{J}}\,(B) \, $  is connected by
Remark 2.10{\it (a)},  so is monic.  As it's also cocommutative,
by a standard result
    \hbox{(cf.~[Ab], Corollary 2.4.29) we deduce that it is
a Hopf algebra.}
                                               \par
   {\it (b)} \, Thanks to  Corollary 3.4{\it (b)},  the
sole non-trivial thing to prove is that  $ \, G_{\underline{D}}\,(B)
\, $  is commutative, for then Lemma 1.5 applies.  Indeed, let  $ \,
0 \not= \overline{a} \in D_r \big/ D_{r-1} \, $,  $ \, 0 \not=
\overline{b} \in D_s \big/ D_{s-1} \, $  (for  $ \, r $,  $ s \in
\N \, $,  \, with  $ \, D_{-1} := \{0\} \, $),  \, and let  $ \, a
\, $  and  $ \, b \, $  be a lift of  $ \, \overline{a} \, $  and
$ \, \overline{b} \, $  in  $ \, D_r \, $  and in  $ \, D_s \, $
respectively: then  $ \, \delta_{r+1}(a) = 0 \, $,  $ \, \delta_r(a)
\not= 0 \, $,  \, and  $ \, \delta_{s+1}(b) = 0 \, $,  $ \, \delta_s(b)
\not= 0 \, $.  Now  $ \, \big( \overline{a} \, \overline{b} - \overline{b}
\, \overline{a} \,\big) \in D_{r+s} \big/ D_{r+s-1} \, $  by construction,
but
  $$  \delta_{r+s} \big( a \, b - b \, a \big)  \hskip3pt  =
\hskip3pt  {\textstyle \sum_{\hskip-10pt  \Sb   \Lambda \cap Y
\not= \emptyset  \\   \Lambda \cup Y = \{1,\dots,r+s\}  \endSb}}
\hskip-10pt  \big( \delta_\Lambda(a) \, \delta_Y(b) - \delta_Y(b)
\, \delta_\Lambda(a) \big)  $$
by  Lemma 3.2{\it (b)}.  In the formula above one has 
$ \, \Lambda \cup Y = \{1,\dots,r+s\} \, $  {\sl and\/}  $ \,
\Lambda \cap Y \not= \emptyset \, $  if and only if either  $ \,
|\Lambda| > r \, $,  \, whence  $ \, \delta_\Lambda(a) = 0 \, $, 
\, or  $ \, |Y| > s \, $,  \, whence  $ \, \delta_Y(b) = 0 \, $, 
\, thus in any case  $ \; \delta_\Lambda(a) \, \delta_Y(b) = 0 \, $,
\, and similarly  $ \; \delta_Y(b) \, \delta_\Lambda(a) = 0 \, $. 
The outset is  $ \, \delta_{r+s} \big( a \, b - b \, a \big) =
0 \, $,  \, whence  $ \, (a \, b - b \, a) \in D_{r+s-1} \, $  so
$ \, \big( \overline{a} \, \overline{b} - \overline{b} \, \overline{a}
\,\big) = \overline{0} \in D_{r+s} \big/ D_{r+s-1} \, $,  \, which
implies commutativity.
                                          \par
   The ``geometrical part'' of the statement is clear by \S 1.2.
In particular, note that the Poisson group  $ G_+ $  is connected
because  $ \, F[G_+] \cong \widetilde{B} \, $  is i.p.-free because
it is  $ \epsilon $--separated,  by Remark 2.4, and then  Lemma
2.2{\it (b)\/}  applies.  In addition, since  $ \, \widetilde{B}
\cong F[G_+] \, $  is graded, when  $ \, p = 0 \, $  the (pro)affine
variety  $ G_+ $  is a  {\sl cone\/}:  but it is also smooth, since
it is an algebraic group, hence it has no vertex.  So it is a (pro)affine
space, say  $ \, G_+ \cong \Bbb{A}_{\,\Bbbk}^{\times \Cal{I}} =
\Bbbk^{\Cal{I}} \, $  for some index set  $ \Cal{I} \, $.  Then  $ \,
\widetilde{B} \cong F[G_+] \cong \Bbbk \big[ {\{x_i\}}_{i \in \Cal{I}}
\,\big] \, $  is a polynomial algebra, as claimed.   
                                             \par
   Finally, when  $ \, p > 0 \, $  the group  $ G_+ $  has dimension
0 and height 1.  Indeed, we must show that  $ \, \bar{\eta}^p = 0
\, $  for each  $ \, \eta \in \widetilde{B} \setminus (\Bbbk \cdot 1)
\, $,  \, which may be taken homogeneous.  Letting  $ \, \eta \in B'
\, $  be any lift of  $ \bar{\eta} $,  we have  $ \, \eta \in D_\ell
\setminus D_{\ell-1} \, $  for a unique  $ \, \ell \in \N \, $,  hence
$ \, \delta_{\ell+1}(x) = 0 \, $.  From  $ \, \Delta^{\ell+1}(\eta) =
\sum_{\Lambda \subseteq \{1, \dots, \ell+1\}} \delta_\Lambda(\eta) \, $
(cf.~\S 2.5) and the multiplicativity of  $ \Delta^{\ell+1} $  we have
 \vskip-14pt
  $$  \displaylines{
   \quad   \Delta^{\ell+1}(\eta^p) \; = \; {\big( \Delta^{\ell+1}(\eta)
\big)}^p \; = \; {\Big( {\textstyle \sum_{\Lambda \subseteq \{1,
\dots, \ell+1\}}} \, \delta_\Lambda(\eta) \Big)}^p  \; \in   \hfill  \cr
   \in \; {\textstyle \sum_{\Lambda \subseteq \{1,\dots,\ell+1\}}} \,
{\delta_\Lambda(\eta)}^p  \; + \;  {\textstyle \sum_{\Sb \! e_1, \dots,
e_p < p \\   \! e_1 + \cdots
+ e_p = p \endSb}}  {\textstyle \Big(\! {p \atop {e_1, \dots, e_p}}}
\!\Big) \, {\textstyle \sum\limits_{\Lambda_1, \dots, \Lambda_p
\subseteq \{1,\dots,\ell+1\}}} \; {\textstyle \prod\limits_{k=1}^p}
\, {\delta_{\Lambda_k}(\eta)}^{e_k}  \; +  \cr
   {} \hfill   +  \; {\textstyle \sum_{k=0}^{\ell}} \,
{\textstyle \sum_{\Sb \Psi \subseteq \{1,\dots,\ell+1\}  \\
|\Psi|=k  \endSb}}  j_\Psi \big( {J_{\!{}_{B'}}}^{\!\!\otimes k} \big)
\;  +  \; {\big( \hbox{ad}_{[\ ,\ ]}(D_{(\ell)})\big)}^{p-1} \big(
D_{(\ell)} \big)  \cr }  $$
 \vskip-6pt
\noindent   (since  $ \, \delta_\Lambda(\eta) \in j_\Lambda
\Big( {J_{\! {}_{B'}}}^{\!\!\otimes |\Lambda|} \Big) \, $
   \hbox{for all  $ \, \Lambda \subseteq \{1, \dots, \ell + 1\} \, $)
with  $ \; D_{(\ell)} := \sum\limits_{\sum_k s_k = \, \ell} \hskip-3pt
\otimes_{k=1}^{\ell+1} D_{s_k} \; $  and}
$ \; {\big(\hbox{ad}_{[\ ,\ ]}(D_{(\ell)})\big)}^{p-1} \! \big(
D_{(\ell)} \big) \! := \! \big[ D_{(\ell)}, \! \big[ D_{(\ell)}, \dots,
\! \big[ D_{(\ell)}, \! \big[ D_{(\ell)}, D_{(\ell)} \big] \big] \cdots
\big] \big] \, $  ($ (p-\!1) $  brackets).  Then
%
%
%
 \vskip-11pt
  $$  \displaylines{
   \delta^{\ell+1}(\eta^p) \; = \; \big( \id_{\scriptscriptstyle H} \!
- \epsilon \big)^{\otimes (\ell+1)} \big( \Delta^{\ell+1}(\eta^p) \big)
\; \in \;  {\delta_{\ell+1}(\eta)}^p  \; +   \hfill  \cr
%
%
%
   \hfill   \; +  {\textstyle \sum\limits_{\Sb  \hskip-1pt  e_1, \dots,
e_p < p \\  \hskip-1pt  e_1 + \cdots e_p = p  \endSb}}  \hskip-3pt
{\textstyle \Big( \! {p \atop {e_1, \dots, e_p}} \!\Big)} \,
{\textstyle \sum\limits_{\cup_k \Lambda_k = \{1,\dots,\ell+1\}}}
\; {\textstyle \prod\limits_{k=1}^p} \, {\delta_{\Lambda_k}\!
(\eta)}^{e_k} \,  +  \, \big( \id_{\scriptscriptstyle H} \!
- \epsilon \big)^{\otimes (\ell+1)} \Big( \hskip-2pt  {\big(
\hbox{ad}_{[\ ,\ ]}(D_{(\ell)}) \big)}^{p-1} \big( D_{(\ell)}
\big) \hskip-2pt  \Big) \, .  \cr }  $$
 \vskip-6pt
   Now,  $ \, {\delta^{\ell+1}(\eta)}^p = 0 \, $  by construction and
$ \, \Big(\! {p \atop {e_1, \dots, e_p}} \!\Big) \, $  is a multiple
of  $ p \, $,  hence it is zero because  $ \, p = \Char(\Bbbk) \, $;  \,
therefore  $ \; \delta_{\ell+1}(\eta) \in \big( \id_{\scriptscriptstyle
H} \! - \epsilon \big)^{\otimes (\ell+1)} \Big( \hskip-2pt {\big(
\hbox{ad}_{[\ ,\ ]} (D_{(\ell)}) \big)}^{p-1} \big( D_{(\ell)} \big)
\hskip-2pt \Big) \, $.  By Lemma 3.2--3,  $ \, D_{s_i} \cdot D_{s_j}
\subseteq D_{s_i + s_j} \, $  and  $ \, \big[ D_{s_i}, D_{s_j} \big]
\subseteq D_{(s_i + s_j) - 1} \, $  by the commutativity of
$ \widetilde{B} $.  This together with Leibniz' rule implies
$ \, {\big( \hbox{ad}_{[\ ,\ ]}(D_{(\ell)}) \big)}^{p-1} \big(
D_{(\ell)} \big) \subseteq \hskip-7pt \sum\limits_{\sum_t r_t =
p\,\ell+1-p} \hskip-4pt \otimes_{t=1}^{\ell+1} D_{r_t} \, $;  \,
moreover  $ \, \big(\id_{\scriptscriptstyle H} \! - \epsilon
\big)^{\otimes (\ell+1)} \Big( \hskip-2pt {\big( \hbox{ad}_{[\ ,\ ]}
(D_{(\ell)}) \big)}^{p-1} \big( D_{(\ell)} \big) \hskip-2pt \Big)
  \subseteq  \hskip0pt {\textstyle \sum_{\,
\Sb  \sum_t r_t = p\,\ell+1-p  \\
r_1, \dots, r_{\ell+1} > 0  \endSb}}
\hskip-7pt  \otimes_{t=1}^{\ell+1} D_{r_t} \, $
because  $ \, D_0 = \hbox{\sl Ker}\,(\delta_1) = \hbox{\sl
Ker}\big(\id_{\scriptscriptstyle H} \! - \epsilon\big) \, $.
In particular, in the last term above we have  $ \, D_{r_1} \subseteq
D_{(p-1)\ell+1-p} := \hbox{\sl Ker}\,\big(\delta_{(p-1)\ell+2-p}\big)
\subseteq \hbox{\sl Ker}\,\big(\delta_{(p-1)\ell}\big) \, $:
\, thus by coassociativity of the  $ \delta_n $'s  
 \vskip-15pt
  $$  \delta_{p\,\ell}(\eta) = \big( \big( \delta_{(p-1)\ell} \otimes
\id^\ell \,\big) \circ \delta_{\ell+1} \big)(\eta)  \hskip3pt  \subseteq
{\textstyle \sum\limits_{\Sb  \sum_t r_t = p\,\ell-1  \\  r_1, \dots,
r_{\ell+1} > 0  \endSb}}  \hskip-3pt  \delta_{(p-1)\ell}(D_{r_1})
\otimes D_{r_2} \otimes \cdots \otimes D_{r_{\ell+1}}  \hskip3pt
=  \hskip3pt  0  $$
 \vskip-4pt
\noindent   i.e.~$ \, \delta_{p\,\ell}(\eta) = 0 \, $.  This means  $ \,
\eta \in D_{p\,\ell-1} \, $,  \, whereas, on the other hand,  $ \, \eta^p
\in D_\ell^{\;p} \subseteq D_{p\,\ell} \, $:  \, then  $ \, \bar{\eta}^p
:= \overline{\eta^p} = \bar{0} \in D_{p\,\ell} \Big/ D_{p\,\ell - 1}
\subseteq \widetilde{B} \, $,  \, by the definition of the product
in  $ \, \widetilde{B} \, $.  Finally, if  $ \Bbbk $  is perfect by
general theory since  $ G_+ $  has dimension 0 and height 1, then 
$ \, F[G_+] \cong \widetilde{B} \, $  is  {\sl truncated polynomial}, 
namely  $ \, F[G_+] \cong \Bbbk \big[ {\{x_i\}}_{i \in \Cal{I}} \big]
\Big/ \big( {\{x_i^{\,p}\}}_{i \in \Cal{I}} \big) \, $  for some
index set  $ \Cal{I} \, $.   \qed
\enddemo

\vskip7pt

{\bf 3.6 Crystal functors on bialgebras and Hopf algebras.} \,
The analysis in the present section shows that, when restricted
to $ \Bbbk $--bialgebras,  
             the output of the previous\break   
 \eject   
\noindent
 functors are
objects of Poisson-geometric type: Lie bialgebras and Poisson groups.
Therefore  $ \, B \mapsto \widehat{B} \, $  and  $ \, B \mapsto
\widetilde{B} \, $  (for  $ \, B \in \calB \, $)  are ``geometrification
functors'', in that they associate to  $ B $  some  {\sl geometrical\/}
symmetries.   
                                          \par
   It is worthwhile pointing out that, by construction, either crystal
functor forgets some information about the initial object, yet still
saves something.  So  $ \widehat{B} $  tells nothing about the coalgebra
structure of  $ B^\vee $  (for all enveloping algebras   --- like
$ \widehat{B} $  ---   roughly look the same from the coalgebra point
of view), yet it grasps some information on its algebra structure. 
Conversely,  $ \widetilde{B} $  gives no information about the algebra
structure of  $ B' $  (in that the latter is simply a polynomial algebra),
but tells something non-trivial about its coalgebra structure.
                                          \par
   We finish this section with the Hopf duality relationship
between these functors:

\vskip7pt

\proclaim{Theorem 3.7} \, Let  $ \, B, P \in \calB \, $  and let
$ \, \pi \, \colon B \times P \loongrightarrow \Bbbk \, $  be
a bialgebra pairing.  Then  $ \pi $  induces  {\sl filtered}
bialgebra pairings  $ \, \pi_f \, \colon B^\vee \times P'
\loongrightarrow \Bbbk \, $,  $ \, \pi^f \colon B' \times P^\vee
\loongrightarrow \Bbbk \, $,  and  {\sl graded}  bialgebra pairings
$ \, \pi_{{}_G} \, \colon \, \widehat{B} \times \widetilde{P}
\loongrightarrow \Bbbk \, $,  $ \, \pi^{{}_G} \, \colon \,
\widetilde{B} \times \widehat{P} \loongrightarrow \Bbbk \, $;  \,
$ \pi_f $  and  $ \pi_{{}_G} $  are perfect on the right, $ \pi^f
$  and  $ \pi^{{}_G} $  on the left.  If in addition  $ \, \pi $
is perfect then all these induced pairings are perfect as well.
If in particular  $ \, B, P \in \HA \, $  are Hopf algebras
and  $ \pi $  is a Hopf algebra pairing, then all the induced
pairings are (filtered or graded) Hopf algebra pairings.
\endproclaim

\demo{Proof} The notion of bialgebra pairing is the ``left-right
symmetrization'' of the notion of augmented pairing when both the
augmented and the coaugmented coalgebra involved are bialgebras. 
Thus the claim above for bialgebras follows simply by a twofold
application of Theorem 2.14.  As for Hopf algebras, it is enough
to remark in addition that the antipode preserves on both sides
of the pairing the filtrations involved (by Lemma 3.3).  Hence,
if the initial pairing is a Hopf one, 
     \hbox{the induced pairings will
clearly be Hopf pairings as well.   \qed}   
\enddemo

\vskip7pt

   {\it $ \underline{\hbox{{\it Remark}}} $:} \, the
$ \delta_\bullet $--filtration  has been independently
introduced and studied by L.~Foissy in the case of graded
Hopf algebras (see [Fo1--2]).

 \vskip1,3truecm

\centerline {\bf \S \; 4 \  Deformations I --- Rees algebras, Rees
coalgebras, etc. }

\vskip10pt

{\bf 4.1 Filtrations and ``Rees objects''.} \, Let  $ V $  be a
vector space over  $ \Bbbk $,  and let  $ \; {\{F_z\}}_{z \in \Z}
:= \underline{F} : \big( \{0\} \subseteq \,\big) \cdots \subseteq
F_{-m} \subseteq \cdots \subseteq F_{-1} \subseteq F_0 \subseteq
F_1 \subseteq \cdots \subseteq F_n \subseteq \cdots \, \big(\,
\subseteq V \,\big) \; $  be a filtration of  $ V $  by vector
subspaces  $ F_z $  ($ z \in \Z $)\,.
                                     \par
   First we define the associated  {\sl blowing space\/}  as the
$ \Bbbk $--subspace  $ \Cal{B}_{\underline{F}}(V) $  of  $ V \big[
t, t^{-1} \big] $  (with  $ t $  an indeterminate) given by  $ \,
\Cal{B}_{\underline{F}}(V) := \sum_{z \in \Z} t^z F_z \, $;  \,
this is isomorphic to  $ \, \bigoplus_{z \in \Z} F_z \, $. 
Second, we define the associated  {\sl Rees module\/}  as the 
$ \Bbbk[t] $--submodule  $ \Cal{R}^t_{\underline{F}}(V) $  of 
$ V \big[ t, t^{-1} \big] $  generated by  $ \Cal{B}_{\underline{F}}
(V) \, $.  Straightforward computations give  $ \Bbbk $--vector 
space isomorphisms
 \vskip-17pt
  $$  \Cal{R}^t_{\underline{F}}(V) \Big/ (t-1) \,
\Cal{R}^t_{\underline{F}}(V) \; \cong \; {\textstyle
\bigcup_{z \in \Z}} F_z =: V^{\underline{F}} \;\, ,  \quad
\Cal{R}^t_{\underline{F}}(V) \Big/ t \, \Cal{R}^t_{\underline{F}}(V)
\; \cong \, {\textstyle \bigoplus\limits_{z \in Z}} F_z \big/ F_{z-1}
=: G_{\underline{F}}(V)  $$
 \vskip-5pt
\noindent
In other words,  $ \Cal{R}^t_{\underline{F}}(V) $  is a
$ \Bbbk[t] $--module  which specializes to  $ \, V^{\underline{F}}
\, $  for  $ \, t = 1 \, $  and specializes to  $ \, G_{\underline{F}}(V)
\, $  for  $ \, t = 0 \, $.  Therefore the  $ \Bbbk $--vector  spaces
$ \, V^{\underline{F}} \, $  and  $ \, G_{\underline{F}}(V) \, $  can
be seen as 1-parameter (polynomial) deformations of each other via
the 1-parameter family of  $ \Bbbk $--vector  spaces  given by
$ \Cal{R}^t_{\underline{F}}(V) $,  in short  $ \;\; V^{\underline{F}}
\hskip3pt \underset{\Cal{R}^t_{\underline{F}}(V)}  \to
{\overset{1 \,\leftarrow\, t \,\rightarrow\, 0}
\to{\longleftarrow\joinrel\relbar\joinrel%
\relbar\joinrel\relbar\joinrel\llongrightarrow}}
\hskip2pt  G_{\underline{F}}(V) \; $.
                                        \par
   We can repeat this construction within the category of algebras,
coalgebras, bialgebras or Hopf algebras over  $ \Bbbk $  with
a filtration in the proper sense: then we'll end up with
corresponding objects  $ \Cal{B}_{\underline{F}}(V) $, 
$ \Cal{R}^t_{\underline{F}}(V) $,  etc.{} of the
same type (algebras, coalgebras, etc.).   

\vskip7pt

  {\bf 4.2 Connecting functors and Rees modules.} \, Let  $ \,
A \in \calA^+ \, $  be an augmented algebra.  By Lemma 2.2 the
$ \epsilon $--filtration  $ \underline{J} $  of  $ A $  is an
algebra filtration, hence we have the Rees  {\sl algebra}  $ \,
\Cal{R}^t_{\underline{J}}(A) \, $.  By the previous analysis,
this yields a 1-parameter deformation
$ \; A  \cong  \Cal{R}^t_{\underline{J}}(A){\Big\vert}_{t=1}
\hskip-1pt  \underset{ {\Cal{R}^t_{\underline{J}}(A)}}
\to  {\overset{1 \,\leftarrow\, t \,\rightarrow\, 0}
\to {\longleftarrow\joinrel\relbar\joinrel%
\relbar\joinrel\relbar\joinrel\llongrightarrow}}
\hskip1pt  \Cal{R}^t_{\underline{J}}(A){\Big\vert}_{t=0}
\hskip-1pt  \cong  G_{\underline{J}}(A) =: \widehat{A} \; $,
\; where  $ \; M{\Big\vert}_{t=c} \!\! := M \! \Big/ (t-c)
\, M \; $  for any  $ \Bbbk[t] $--module  $ M $  and any  $ \,
c \in \Bbbk \, $.  All fibers in this deformation are pairwise
isomorphic as vector spaces  {\sl but\/}  perhaps for the fiber
at  $ \, t = 0 \, $,  i.e.~exactly  $ \widehat{A} $,  for there
the subspace  $ \, J^\infty \, $  is ``shrunk to zero''.  This is
settled passing from  $ A $  to  $ A^\vee $,  for which the same
deformation is  {\sl regular},  i.e.~{\sl all\/}  fibers in it
are pairwise isomorphic (as vector spaces): the scheme is
 \vskip-13pt
  $$  A^\vee := A \Big/ J^\infty  \, \cong \,
\Cal{R}^t_{\underline{J}}\big(A^\vee\big){\Big\vert}_{t=1}
\hskip-1pt  \underset{ {\Cal{R}^t_{\underline{J}}(A^\vee)}}
\to {\overset{1 \,\leftarrow\, t \,\rightarrow\, 0}
\to{\longleftarrow\joinrel\relbar\joinrel\relbar\joinrel%
\relbar\joinrel\relbar\joinrel\relbar\joinrel\llongrightarrow}}
\hskip2pt  \Cal{R}^t_{\underline{J}}(A){\Big\vert}_{t=0}
\hskip-1pt  \, \cong \,  G_{\underline{J}}(A) =: \widehat{A}
\eqno (4.1)  $$
 \vskip-6pt
\noindent
where we implicitly used the identities  $ \, 
\Cal{R}^t_{\underline{J}}(A^\vee){\Big\vert}_{t=0}
= \, \widehat{A^\vee} \, = \, \widehat{A} \, = \,
\Cal{R}^t_{\underline{J}}(A) {\Big\vert}_{t=0} \, $.
                                              \par
   The situation is (dually!) similar for the connecting functor
on coaugmented coalgebras.  Indeed, let  $ \, C \in \calC^+ \, $:
by Lemma 2.8 the  $ \delta_\bullet $--filtration  $ \underline{D} $
of  $ C $  is a coalgebra filtration, so we have the Rees  {\sl
coalgebra}  $ \, \Cal{R}^t_{\underline{D}}(C) \, $.  The analysis
above yields a 1-parameter deformation
 \vskip-13pt
  $$  C' := {\textstyle \bigcup_{n \in \N}} \, D_n  \,
\cong \,  \Cal{R}^t_{\underline{D}}(C){\Big\vert}_{t=1}
\hskip-1pt  \underset{ {\Cal{R}^t_{\underline{D}}(C')}}
\to  {\overset{1 \,\leftarrow\, t \,\rightarrow\, 0}
\to {\longleftarrow\joinrel\relbar\joinrel\relbar%
\joinrel\relbar\joinrel\llongrightarrow}}  \hskip2pt
\Cal{R}^t_{\underline{D}}(C){\Big\vert}_{t=0}  \hskip-1pt
\, \cong \,  G_{\underline{D}}(C) =: \widetilde{C}
\eqno (4.2)  $$
 \vskip-6pt
\noindent
which is also  {\sl regular},  i.e.~all its fibers are pairwise
isomorphic as vector spaces.

\vskip7pt

  {\bf 4.3 The bialgebra and Hopf algebra case.} \, We now consider
a bialgebra  $ \, B \in \calB \, $.  In this case, the results of
\S 3 ensure that  $ B^\vee $  is a bialgebra,  $ \widehat{B} $  is a
(graded, etc.) Hopf algebra, and  $ \Cal{R}^t_{\underline{J}}(B^\vee) $
is a  $ \Bbbk[t] $--bialgebra,  because  $ \underline{J} $  is a
bialgebra filtration.  Using Theorem 3.5, formula (4.1) becomes
$ \; B^\vee \underset{ {\Cal{R}^t_{\underline{J}}(B^\vee)}}
\to {\overset{1 \,\leftarrow\, t \,\rightarrow\, 0}
\to{\longleftarrow\joinrel\relbar\joinrel\relbar\joinrel%
\relbar\joinrel\relbar\joinrel\relbar\joinrel\llongrightarrow}}
\hskip1pt  \widehat{B} \, \cong \, \U(\gerg_-) \; $  for some restricted
Lie bialgebra  $ \gerg_- $  as in  Theorem 3.5{\it (a)}.  Similarly,
by \S 3 we know that  $ B' $  is a Hopf algebra,  $ \widetilde{B} $ is a (graded, etc.) Hopf algebra and  $ \Cal{R}^t_{\underline{D}}(B') $ 
is a Hopf  $ \Bbbk[t] $--algebra,  because  $ \, \underline{D}(B') =
\underline{D}(B) \, $  is a Hopf algebra filtration of  $ B \, $. 
Thus, again by Theorem 3.5, formula (4.2) becomes  $ \; B' \hskip0pt
\underset{ {\Cal{R}^t_{\underline{D}}(B')}} \to {\overset{1 \,
\leftarrow \, t \, \rightarrow \, 0}
\to {\longleftarrow\joinrel\relbar\joinrel\relbar%
\joinrel\relbar\joinrel\llongrightarrow}}  \hskip1pt
\widetilde{B} \cong F[G_+] \; $  for some Poisson algebraic group
$ G_+ \, $,  as in  Theorem 3.5{\it (b)}.  The overall outcome is
  $$  F[G_+]  \, \cong \,  \widetilde{B}  \hskip2pt
\underset{ {\Cal{R}^t_{\underline{D}}(B')}}  \to
{\overset{0 \,\leftarrow\, t \,\rightarrow\, 1}
\to {\longleftarrow\joinrel\relbar\joinrel\relbar%
\joinrel\relbar\joinrel\llongrightarrow}}  \hskip2pt
B'  \,{\lhook\joinrel\relbar\joinrel\rightarrow}\,  B
\relbar\joinrel\twoheadrightarrow  B^\vee  \hskip0pt
\underset{ {\Cal{R}^t_{\underline{J}}(B^\vee)}}
\to {\overset{1 \,\leftarrow\, t \,\rightarrow\, 0}
\to{\longleftarrow\joinrel\relbar\joinrel\relbar\joinrel%
\relbar\joinrel\relbar\joinrel\relbar\joinrel\longrightarrow}}
\hskip2pt  \widehat{B}  \; \cong \, \U(\gerg_-)   \eqno (4.3)  $$
This drawing shows how the bialgebra  $ B $  gives rise to two Hopf
algebras of Poisson geometrical type, namely  $ F[G_+] $  on the
left-hand side and  $ \U(\gerg_-) $  on the right-hand side, through
bialgebra morphisms and regular bialgebra deformations.  Indeed, in
both cases one has first a ``reduction step'',  i.e.~$ \, B \mapsto
B' \, $  or  $ \, B \mapsto B^\vee $,  (yielding ``connected'' objects,
cf.~Remark 2.12{\it (a)\/}),  then a regular 1-parameter deformation
via Rees bialgebras.
                                         \par
   Finally, if  $ \, H \in \HA \, $  is a Hopf algebra then all
objects in (4.3) are Hopf algebras too, i.e.~also  $ H^\vee $
(over  $ \Bbbk $)  and  $ \Cal{R}^t_{\underline{J}}(B^\vee) $
(over  $ \Bbbk[t] \, $).  Therefore (4.3) reads
  $$  F[G_+]  \, \cong \,  \widetilde{H}  \hskip1pt
\underset{ {\Cal{R}^t_{\underline{D}}(H')}}  \to
{\overset{0 \,\leftarrow\, t \,\rightarrow\, 1} \to
{\longleftarrow\joinrel\relbar\joinrel\relbar\joinrel%
\relbar\joinrel\llongrightarrow}}  \hskip1pt  H'
\,{\lhook\joinrel\relbar\joinrel\rightarrow}\,  H
\relbar\joinrel\twoheadrightarrow  H^\vee  \hskip0pt
\underset{ {\Cal{R}^t_{\underline{J}}(H^\vee)}}
\to {\overset{1 \,\leftarrow\, t \,\rightarrow\, 0}
\to{\longleftarrow\joinrel\relbar\joinrel\relbar\joinrel%
\relbar\joinrel\relbar\joinrel\relbar\joinrel\longrightarrow}}
\hskip1pt   \widehat{H}  \; \cong \, \U(\gerg_-)   \eqno (4.4)  $$
with the one-way arrows being morphisms  {\sl of Hopf algebras\/}
and the two-ways arrows being 1-parameter regular deformations
{\sl of Hopf algebras}.  When in addition  $ H $  is connected,
that is  $ \, H = H' \, $,  and ``coconnected'', that is  $ \,
H = H^\vee $,  \, the scheme (4.4) takes the simpler form
  $$  F[G_+]  \, \cong \,  \widetilde{H}  \hskip1pt
\underset{ {\Cal{R}^t_{\underline{D}}(H)}}  \to
{\overset{0 \,\leftarrow\, t \,\rightarrow\, 1}
\to {\longleftarrow\joinrel\relbar\joinrel\relbar\joinrel%
\relbar\joinrel\relbar\joinrel\llongrightarrow}}  \hskip1pt
H  \underset{ {\Cal{R}^t_{\underline{J}}(H)}} \to
{\overset{1 \,\leftarrow\, t \,\rightarrow\, 0}
\to {\longleftarrow\joinrel\relbar\joinrel\relbar\joinrel%
\relbar\joinrel\relbar\joinrel\llongrightarrow}}  \hskip1pt
\widehat{H}  \; \cong \, \U(\gerg_-)   \eqno (4.5)  $$
which means we can (regularly) deform  $ H $  itself to ``Poisson
geometrical'' Hopf algebras.

\vskip7pt

   {\it $ \underline{\hbox{{\it Remarks}}} $:} \, {\it (a)} \, There
is no simple relationship,  {\it a priori},  between the Poisson group
$ G_+ $  and the Lie bialgebra  $ \gerg_- $  in (4.3) or (4.4), or
even (4.5): examples do show that.  In particular, either  $ G_+ $
or  $ \gerg_- $  may be trivial while the other is not (see \S 8).
                                             \par
   {\it (b)} \, The Hopf duality relationship between connecting
functors of the two types explained in \S\S 2.11--14 extends to
deformations via Rees modules.  Indeed, by construction there is
a neat category-theoretical duality between the definition of
$ \Cal{R}^t_{\underline{J}}(A) $  and of
$ \Cal{R}^t_{\underline{D}}(C) $  (for  $ \,
A \in \calA^+ \, $  and  $ \, C \in \calC^+ $)\,.
Even more, Theorem 2.14 extends to the following:

\vskip7pt

\proclaim{Theorem 4.4} \, Let  $ \, A \in \calA^+ $,  $ C \in
\calC^+ $  and let  $ \, \pi \colon \, A \times C \loongrightarrow
\Bbbk \, $  be an augmented pairing.  Then  $ \pi $  induces an
augmented pairing  $ \, \pi_{\Cal{R}} \, \colon
\Cal{R}^t_{\underline{J}}(A) \times \Cal{R}^t_{\underline{D}}(C)
\loongrightarrow \Bbbk[t] \, $  which is perfect on the right.  If
in addition  $ \, \pi $  is perfect then  $ \, \pi_{\Cal{R}} $  is
perfect as well, and
  $$  \hbox{ $ \matrix
   \Cal{R}^t_{\underline{J}}(A)  &  \hskip-5pt  = \, \Big\{\, \eta \in
A(t) \,\Big|\, \pi_t\big(\eta,\kappa\big) \in \Bbbk[t] \, , \; \forall
\, \kappa \in \Cal{R}^t_{\underline{D}}(C) \Big\}  &  \hskip-5pt  =: \,
\big( \Cal{R}^t_{\underline{D}}(C) \big)^\bullet  \\  
   \Cal{R}^t_{\underline{D}}(C)  &  \hskip-5pt  = \, \Big\{\, \kappa \in
C(t) \,\Big|\, \pi_t\big(\eta,\kappa\big) \in \Bbbk[t] \, , \; \forall
\, \eta \in \Cal{R}^t_{\underline{J}}(A) \Big\}  &  \hskip-5pt  =: \,
\big( \Cal{R}^t_{\underline{J}}(A) \big)^\bullet  \endmatrix $ }
\eqno (4.6)  $$
where  $ \, S(t) := \Bbbk(t) \otimes_\Bbbk S \, $  for
$ \, S \in \{A,C\} \, $  and  $ \, \pi_t \, \colon \, A(t)
\times C(t) \loongrightarrow \Bbbk(t) \, $  is the obvious
$ \Bbbk(t) $--linear  pairing induced by  $ \, \pi \, $.  If
$ A $  and  $ C $  are bialgebras, resp.~Hopf algebras, then
everything holds with bialgebra, resp.~Hopf algebra, pairings
instead of augmented pairings.
\endproclaim

\demo{Proof} The obvious augmented  $ \Bbbk(t) $--linear  pairing
$ \, \pi_t \, \colon \, A(t) \times C(t) \loongrightarrow \Bbbk(t)
\, $  induced by  $ \, \pi \, $  restricts to an augmented pairing
$ \, \pi_{\Cal{R}} \, \colon \Cal{R}^t_{\underline{J}}(A) \times
\Cal{R}^t_{\underline{D}}(C) \loongrightarrow \Bbbk(t) \, $.
By (2.1) both these pairings are perfect on the right, and
$ \pi_{\Cal{R}} $  takes value into  $ \Bbbk[t] $.  In addition,
in the proof of Theorem 2.14 we showed that if  $ \pi $  is perfect
formula (2.1) improves, in that  $ \, {J_{\!{}_A}}^{\! n+1} = \big(
D^{\scriptscriptstyle C}_n \big)^\perp \, $  for all  $ \, n \in \N
\, $.  This implies that both  $ \pi_t $  and  $ \pi_{\Cal{R}} $
are perfect as well, and also that the identities (4.6) do hold,
q.e.d.  The final part about bialgebras or Hopf algebras is clear.
\qed
\enddemo
%
%
%
 \eject

\centerline {\bf \S \; 5 \  Deformations II --- from Rees
bialgebras to quantum groups }

\vskip10pt

{\bf 5.1 From Rees bialgebras to quantum groups via the GQDP.} \,
In this section we show how, for any  $ \Bbbk $--bialgebra  $ B $,
we can get another deformation scheme like (4.3).  In fact, this
will be built upon the latter, applying part of the ``Global
Quantum Duality Principle'' explained in [Ga1--2] in its
stronger form about bialgebras.
                                     \par
   Indeed, the deformations in (4.3) were realized through Rees
bialgebras, namely  $ \Cal{R}^t_{\underline{J}}\big(B^\vee\big) $ 
and  $ \Cal{R}^t_{\underline{D}}\big(B'\big) $.  These are torsion-free
(actually, free) as  $ \Bbbk[t] $--modules,  hence one can apply the
construction made in  [{\it loc.~cit.}]  via the so-called Drinfeld's
functors, to get some new torsion-free  $ \Bbbk[t] $--bialgebras.  The
latters (just like the Rees bialgebras we start from) again specialize
to special bialgebras at  $ \, t = 0 \, $.  In particular, if  $ B $ 
is a Hopf algebra the new bialgebras are Hopf algebras too, and
precisely ``quantum groups'' in the sense of  [{\it loc.~cit.}].
                                     \par
   To begin with, set  $ \, B^\vee_t := \Cal{R}^t_{\underline{J}}
\big(B^\vee\big) \, $:  \, this is a torsion-free, 
$ \Bbbk[t] $--bialgebra.  We define   
  $$  \big( B^\vee_t \big)' \, := \, \Big\{\, b \in B^\vee_t \,\Big|\;
\delta_n(b) \in t^n \big(B^\vee_t\big)^{\otimes n} \, ,  \; \forall
\;\, n \in \N \,\Big\} \; .  $$
   \indent   On the other hand, let  $ \, B'_t :=
\Cal{R}^t_{\underline{D}}(B) = \Cal{R}^t_{\underline{D}}\big(B'\big)
\, $:  \, this is a torsion-free  $ \Bbbk[t] $--bialgebra  as well. 
Letting  $ \, J' := \text{\sl Ker}\,\big( \epsilon \, \colon \, B'_t
\longrightarrow \Bbbk[t] \,\big) \, $  and  $ \, B'(t) := \Bbbk(t)
\otimes_{\Bbbk[t]} B'_t = \Bbbk(t) \otimes_{\Bbbk} B' \, $,  \, set
  $$  \big( B'_t \big)^{\!\vee} \, := \, {\textstyle \sum_{n \geq 0}}
\, t^{-n} \big( J' \,\big)^n \, = \, {\textstyle \sum_{n \geq 0}}
\big( t^{-1} J' \,\big)^n  \qquad  \big(\, \subseteq B'(t) \,\big)
\; .  $$

   \indent   The first important point is the following:

\vskip7pt

\proclaim{Proposition 5.2} Both  $ \big( B^\vee_t \big)' $
and  $ \big( B'_t \big)^{\!\vee} $  are free (hence torsion-free)\/
$ \Bbbk[t] $--bialgebras;  moreover, the mappings  $ \, B \mapsto
\big( B^\vee_t \big)' \, $  and  $ \, B \mapsto \big( B'_t
\big)^{\!\vee} \, $  are functorial.  The analogous results hold
for Hopf\/  $ \Bbbk $--algebras,  just replacing ``bialgebra(s)''
with ``Hopf algebra(s)'' throughout.
\endproclaim

\demo{Proof} This is a special instance of Theorem 3.6 in [Ga1] (or
in [Ga2]), but for the fact that the final objects are free, not only
torsion-free.  This is easy to see: taking a  $ \Bbbk $--basis  of
$ B^\vee $  or  $ B' $  respectively which is ``compatible with
$ \underline{J} \, $'',  resp.~``with  $ \underline{D} \, $'',  in
the obvious sense, one gets at once from that a  $ \Bbbk[t] $--basis
of  $ B^\vee_t $  or of  $ B'_t $,  which then are free.  Now
$ \big(B^\vee_t\big)' $  is a  $ \Bbbk[t] $--submodule  of the
free  $ \Bbbk[t] $--module  $ B^\vee_t $,  hence it is free as well. 
On the other hand, one can rearrange a basis of  $ B'_t $  so to get
another one ``compatible'' with the filtration  $ \, \Big\{ \big( J'
\, \big)^n \Big\}_{n \in \N} \, $.  Then from this new basis one
immediately obtains a  $ \Bbbk[t] $--basis for  $ \big( B'_t
\big)^{\!\vee} $,  via suitable rescaling by negative powers
of  $ t \, $.   
                                           \par
   The Hopf algebra case goes through in the same way, by similar
arguments.   \qed  
\enddemo  

\vskip3pt

   Next we show that  $ \big( B^\vee_t \big)' $  and  $ \big( B'_t
\big)^{\!\vee} $  are regular deformations of  $ B^\vee $  and  $ B' $
respectively:

\vskip5pt

\proclaim{Proposition 5.3} \, Let  $ B $  be a  $ \Bbbk $--bialgebra 
(or a Hopf\/  $ \Bbbk $--algebra).  Then:   
                                   \hfill\break
   \indent   (a)  $ \; \big( B^\vee_t \big)'{\Big|}_{t=1} :=
\, \big( B^\vee_t \big)' \! \Big/ (t-1) \big( B^\vee_t \big)'
\, \cong \, B^\vee \; $  as \,  $ \Bbbk $--bialgebras  (or as
Hopf\/  $ \Bbbk $--algebras);   
                                   \hfill\break
   \indent   (b)  $ \; \big( B'_t \big)^{\!\vee}{\Big|}_{t=1} := \,
\big( B'_t \big)^{\!\vee} \! \Big/ (t-1) \big( B'_t \big)^{\!\vee}
\, \cong \, B' \; $  as \,  $ \Bbbk $--bialgebras  (or as
Hopf\/  $ \Bbbk $--algebras).   
\endproclaim  

\demo{Proof} Claim  {\it (b)\/} follows from the chain of isomorphisms
$ \; \big( B'_t \big)^{\!\vee}{\Big|}_{t=1} \cong \, B'_t{\Big|}_{t=1}
\cong \, B' \; $  of  $ \Bbbk $--bialgebras, which follow directly
from definitions.  As for  {\it (a)},  we can define another
$ \Bbbk $--bialgebra  $ \Big( \! \big( B^\vee_t \big)' \Big)^{\!\vee} $
as we did for  $ \big( B'_t \big)^{\!\vee} $:  the very definition then
gives  $ \; \Big( \! \big( B^\vee_t \big)' \Big)^{\!\vee}{\Big|}_{t=1}
\cong \, \big( B^\vee_t \big)'{\Big|}_{t=1} \, $.   
      \hbox{The claim then follows from  $ \Big(\!
\big( B^\vee_t \big)' \Big)^{\!\!\vee} \!\! = B^\vee_t $,
by  Theorem 3.6{\it (c)\/}  of [Ga1--2].   $ \! \square $}
\enddemo

\vskip4pt

   The key result is then the following

\vskip5pt

\proclaim{Theorem 5.4}
                                   \hfill\break
   \indent   (a)  $ \; \big( B^\vee_t \big)'{\Big|}_{t=0} \! :=
\big( B^\vee_t \big)' \! \Big/ t \, \big( B^\vee_t \big)' \; $ 
is a commutative, i.p.-free Poisson\/  $ \Bbbk $--bialgebra. 
   More-\break
over,
if  $ \, p > 0 \, $  each non-zero element of  $ \, \text{\sl Ker}
\,\Big( \epsilon \, \colon \, \big( B^\vee_t \big)'{\Big|}_{t=0}
\longrightarrow \Bbbk \,\Big) \, $
     \hbox{has nilpotency order  $ p \, $.}
Therefore  $ \big( B^\vee_t \big)'{\Big|}_{t=0} $  is the function
algebra  $ F[M] $  of some connected Poisson algebraic monoid  $ M $,
and if  $ \, p > 0 \, $  then  $ M $  has dimension 0 and height 1.
                                   \hfill\break
   \indent   If in addition  $ \, B = H \in \HA \, $  is a Hopf\/ 
$ \Bbbk $--algebra,  then  $ \big( H^\vee_t \big)'{\Big|}_{t=0} $
is a Poisson Hopf\/  $ \Bbbk $--algebra, and  $ \, K_+ := \text{\it
Spec}\, \Big(\! \big( H^\vee_t \big)'{\Big|}_{t=0} \Big) = M \, $ 
is a (connected) algebraic Poisson group.
                                   \hfill\break
    \indent
   (b)  $ \; \big( B'_t \big)^{\!\vee}{\Big|}_{t=0} := \, \big(
B'_t \big)^{\!\vee} \! \Big/ t \, \big( B'_t \big)^{\!\vee} \; $
is a connected cocommutative Hopf\/  $ \Bbbk $--algebra  generated
by  $ P \Big(\! \big( B'_t \big)^{\!\vee}{\Big|}_{t=0} \Big) $.
Therefore  $ \, \big( B'_t \big)^{\!\vee}{\Big|}_{t=0} \! =
\U(\gerk_-) \, $  for some Lie bialgebra\/  $ \gerk_- \, $.
                                   \hfill\break
   \indent   (c)  If  $ \, p = 0 \, $  and  $ \, B = H \in \HA \, $
is a Hopf\/  $ \Bbbk $--algebra,  let  $ \, \widehat{H} = U(\gerg_-)
\, $  and  $ \, \widetilde{H} = F[G_+] \, $  as in Theorem 3.5.  Then
(notation of (a) and (b))  $ \; K_+ = G_-^\star \; $  and  $ \,\;
\gerk_- = \, \gerg_+^\times \, $,  \; that is  $ \; \text{\sl coLie}
\,(K_+) \, = \, \gerg_- \; $  and  $ \,\; \gerk_- = \, \text{\sl coLie}
\,(G_+) \; $  as Lie bialgebras.
\endproclaim

\demo{Proof} {\it (a)} \,  This follows from Theorem 3.8 in
[Ga1--2], applied to the  $ \Bbbk[t] $--bialgebra  $ B^\vee_t \, $.
                                               \par
   {\it (b)} \,  This follows from Theorem 3.7 in [{\it loc.~cit.}],
applied to the  $ \Bbbk[t] $--bialgebra  $ B'_t \, $.
                                               \par
   {\it (c)} \,  This follows from Theorem 4.8 in  [{\it loc.~cit.}]
applied to the Hopf  $ \Bbbk[t] $--algebra  $ H^\vee_t \, $,  \, because
by construction the latter is a QrUEA (in the terminology therein),
and from Theorem 4.7 in  [{\it loc.~cit.}]  applied to the Hopf
$ \Bbbk[t] $--algebra  $ H'_t \, $,  \, because by construction the
latter is a QFA (in the terminology therein, noting that the condition
$ \; \bigcap\limits_{n \in \N} \big( \text{\sl Ker}\,( H \,{\buildrel
\epsilon \over \longrightarrow}\, \Bbbk[t] \,{\buildrel {t \mapsto 0}
\over \loongrightarrow}\, \Bbbk ) \big)^n = \bigcap\limits_{n \in \N}
t^n H \; $  required in  [{\it loc.~cit.}]  for a QFA is satisfied).  
\qed
\enddemo

\vskip7pt

{\bf 5.5 Deformations through Drinfeld's functors.} \, The outcome of
the previous analysis is that, for each  $ \Bbbk $--bialgebra  $ \,
B \in \calB \, $,  there is a second scheme   --- besides (4.3) ---  
which yields regular 1-parameter deformations, namely
  $$  \U(\gerk_-)  \hskip1pt
\underset{ (B'_t)^\vee}  \to
{\overset{0 \,\leftarrow\, t \,\rightarrow\, 1}
\to {\longleftarrow\joinrel\relbar\joinrel\relbar%
\joinrel\relbar\joinrel\llongrightarrow}}  \hskip1pt  B'
\hskip1pt  \,{\lhook\joinrel\relbar\joinrel\rightarrow}\,
\hskip1pt  B  \relbar\joinrel\relbar\joinrel\twoheadrightarrow
B^\vee  \hskip0pt  \underset{(B^\vee_t)'}
\to  {\overset{1 \,\leftarrow\, t \,\rightarrow\, 0} \to
{\longleftarrow\joinrel\relbar\joinrel\relbar\joinrel\relbar%
\joinrel\relbar\joinrel\relbar\joinrel\longrightarrow}}
\hskip1pt   F[M]   \eqno (5.1)  $$
   \indent   This gives another recipe, besides (4.3), to make
two other bialgebras of Poisson geometrical type, namely  $ F[M] $
and  $ \U(\gerk_-) $,  out of the bialgebra  $ B $,  through bialgebra
morphisms and regular bialgebra deformations.  Like for (4.3), in both
cases there is first the ``reduction step''  $ \, B \mapsto B' \, $  or
$ \, B \mapsto B^\vee $\,  and then regular 1-parameter deformations
via  $ \Bbbk[t] $--bialgebras.  However, this time on right-hand side
we have in general only a bialgebra, not a Hopf algebra.  When  $ \, B
= H \in \HA \, $  is a Hopf  $ \Bbbk $--algebra,  then (5.1) improves:
all objects therein are Hopf algebras too, and morphisms and
deformations are ones of Hopf algebras.  In particular  $ \,
M = K_+ \, $  is a Poisson  {\sl group}, 
     \hbox{not only a monoid: at a glance}
  $$
\U(\gerk_-)  \hskip2pt  \underset{(H'_t)^\vee }  \to
{\overset{0 \,\leftarrow\, t \,\rightarrow\, 1}
\to {\longleftarrow\joinrel\llongrightarrow}}  \hskip3pt  H'
\hskip3pt  {\lhook\joinrel\relbar\joinrel\relbar\joinrel\rightarrow}
\hskip3pt  H  \hskip2pt
{\relbar\joinrel\joinrel\relbar\joinrel\twoheadrightarrow}
\hskip4pt  H^\vee  \hskip-1pt  \underset{(H^\vee_t)'} \to
{\overset{1 \,\leftarrow\, t \,\rightarrow\, 0}
\to {\longleftarrow\joinrel\llongrightarrow}}
\hskip1pt  F[K_+]   \eqno (5.2)  $$
This yields another recipe, besides (4.4), to make two Hopf algebras
of Poisson geometrical type   ---  i.e.~$ F[K_+] $  and  $ \U(\gerk_-) $ 
---   out of  $ H $,  through Hopf algebra morphisms and regular Hopf
algebra deformations.  Again we have first a ``reduction step''  $ \,
H \mapsto H' \, $  or  $ \, H \mapsto H^\vee $,  then regular
1-parameter deformations via Hopf  $ \Bbbk[t] $--algebras.   
                                         \par
   In the special case when  $ H $  is connected, that is  $ \, H
= H' \, $,  and ``coconnected'', that is  $ \, H = H^\vee $,  \,
formula (5.2) takes the simpler form, the analogue of (4.5),
  $$  \U(\gerk_-)  \hskip3pt  \underset{(H'_t)^\vee }  \to
{\overset{0 \,\leftarrow\, t \,\rightarrow\, 1}
\to {\longleftarrow\joinrel\relbar\joinrel\relbar%
\joinrel\relbar\joinrel\llongrightarrow}}  \hskip3pt
H  \hskip3pt  \underset{(H^\vee_t)'} \to
{\overset{1 \,\leftarrow\, t \,\rightarrow\, 0}
\to {\longleftarrow\joinrel\relbar\joinrel\relbar%
\joinrel\relbar\joinrel\llongrightarrow}}
\hskip3pt  F[K_+]   \eqno (5.3)  $$
so we can again (regularly) deform  $ H $  into Poisson geometrical
Hopf algebras.
                             \par
   In particular, when  $ \, H' = H = H^\vee \, $  patching together
$ (4.5) $  and  $ (5.3) $  we find
  $$  \hskip5pt   F[G_+]  \hskip3pt
\underset{\;H'_t}  \to
{\overset{0 \,\leftarrow\, t \,\rightarrow\, 1}
\to{\longleftarrow\joinrel\relbar\joinrel%
\relbar\joinrel\relbar\joinrel\relbar\joinrel\llongrightarrow}}
\hskip2pt {}  H'  \hskip1pt
\underset{\;(H'_t)^\vee}  \to
{\overset{1 \,\leftarrow\, t \,\rightarrow\, 0}
\to{\longleftarrow\joinrel\relbar\joinrel%
\relbar\joinrel\relbar\joinrel\relbar\joinrel\llongrightarrow}}
\hskip3pt  \U(\gerk_-)
\hskip11pt  \Big( = U\big(\gerg_+^{\,\times}\big)  \,\text{\ if \ }
\text{\it Char}\,(\Bbbk) = 0 \Big)  $$
 \vskip-25pt
  $$  \hskip-133pt   \Big|\Big|  $$
 \vskip-20pt
  $$  \hskip-129pt  H_{\phantom{\displaystyle I}}  $$
 \vskip-21pt
  $$  \hskip-133pt   \Big|\Big|  $$
 \vskip-29pt
  $$  \hskip9pt   \U(\gerg_-)  \hskip2pt
\underset{H^\vee_t}  \to
{\overset{0 \,\leftarrow\, t \,\rightarrow\, 1}
\to{\longleftarrow\joinrel\relbar\joinrel%
\relbar\joinrel\relbar\joinrel\relbar\joinrel\llongrightarrow}}
\hskip2pt {}  H^\vee  \hskip1pt
\underset{\;(H^\vee_t)'}  \to
{\overset{1 \,\leftarrow\, t \,\rightarrow\, 0}
\to{\longleftarrow\joinrel\relbar\joinrel%
\relbar\joinrel\relbar\joinrel\relbar\joinrel\llongrightarrow}}
\hskip2pt  F[K_+]
\hskip9pt  \Big( = F\big[G_-^\star\big]  \,\text{\ if \ }
\!\text{\it Char}\,(\Bbbk) = 0 \Big)  $$
 This gives  {\sl four\/}  different regular 1-parameter deformations
from  $ H $  to Hopf algebras encoding geometrical objects of Poisson
type (i.e.~Lie bialgebras or Poisson algebraic groups).

\vskip7pt

{\bf 5.6 Drinfeld-like functors.} \, The constructions in the present
section show that mapping  $ \, B \mapsto \big( B'_t \big)^{\!\vee}
{\big|}_{t=0} \, $  and  $ \, B \mapsto \big( B^\vee_t \big)'{\big|}_{t=0}
\, $  (for all  $ \, B \in \calB \, $)  define two endofunctors of 
$ \calB \, $.  Their output describes objects of Poisson-geometric type,
namely Lie bialgebras and connected Poisson algebraic monoids.  Therefore,
both  $ \, B \mapsto \big( B'_t \big)^{\!\vee}{\big|}_{t=0} \, $  and 
$ \, B \mapsto \big( B^\vee_t \big)'{\big|}_{t=0} \, $  (for  $ \,
B \in \calB \, $)  are  {\sl geometrification functors\/}  on 
$ \Bbbk $--bialgebras,  just like the ones in \S 3.6: we call
them ``Drinfeld-like functors''.  Thus we have four functorial
recipes (our  {\sl geometrification functors\/})  to sort out
of the generalized symmetry  $ B $  some  {\sl geometrical\/}
symmetries.  Now we explain the duality relation between
Drinfeld-like functors and associated deformations:

\vskip7pt

\proclaim{Theorem 5.7} \, Let  $ \, B, P \in \calB \, $,  and let
$ \, \pi \, \colon \, B \times P \loongrightarrow \Bbbk \, $  be
a\/  $ \Bbbk $--bialgebra  pairing.  Then  $ \pi $  induces a\/
$ \Bbbk[t] $--bialgebra  pairing  $ \, \pi'_\vee \, \colon \big(
B^\vee_t \big)' \times \big(P_t^{\,\prime}\big)^{\!\vee} \loongrightarrow
\Bbbk[t] \, $  and a\/  $ \Bbbk $--bialgebra  pairing  $ \, \pi'_\vee
{\big|}_{t=0} \, \colon \big(B^\vee_t\big)'{\Big|}_{t=0} \times \big(
P_t^{\,\prime} \big)^{\!\vee}{\Big|}_{t=0} \! \loongrightarrow \Bbbk
\, $.  If in addition  $ \, \pi $  is perfect and  $ \, \Char(\Bbbk)
= 0 \,$  then both the induced pairings are perfect as well, and
 \vskip-7pt   
  $$  \hbox{ $ \matrix   
   \big(B^\vee_t\big)' &  \hskip-6pt  = \, \Big\{\, \eta \in B(t)
\,\Big|\, \pi_t\big(\eta,\kappa\big) \in \Bbbk[t] \, , \; \forall
\, \kappa \in \big(P'_t\big)^{\!\vee} \Big\}  &  \hskip-5pt  =: \,
\Big(\! \big(P'_t\big)^{\!\vee} \Big)^{\!\bullet}  \\  
   \big(P'_t\big)^{\!\vee} &  \hskip-6pt  = \, \Big\{\, \kappa
\in P(t) \,\Big|\, \pi_t\big(\eta,\kappa\big) \in \Bbbk[t] \, ,
\; \forall \, \eta \in \big(B^\vee_t\big)' \Big\}  &  \hskip-5pt 
=: \, \Big(\! \big(B^\vee_t\big)' \Big)^{\!\bullet}  \endmatrix $ }
\eqno (5.4)  $$
 \vskip-1pt   
\noindent   
 where  $ \, S(t) := \Bbbk(t) \otimes_\Bbbk S \, $  for
$ \, S \in \{B,P\} \, $  and  $ \, \pi_t \, \colon \, B(t)
\times P(t) \loongrightarrow \Bbbk(t) \, $  is the obvious
$ \Bbbk(t) $--linear  pairing induced by  $ \, \pi \, $.  If
$ \, B, P \in \HA \, $  are ~Hopf\/  $ \Bbbk $--algebras  then
everything holds with Hopf\/  $ \Bbbk $--algebra  pairings
instead of bialgebra pairings.
\endproclaim

\demo{Proof} By Theorem 4.4,  $ \pi $  induces a pairing  $ \,
\pi_{\Cal{R}} \, \colon \, B^\vee_t \times P_t^{\,\prime}
\loongrightarrow \Bbbk[t] \, $
of $ \Bbbk[t] $--bialgebras.  Applying to it Proposition 4.4 in [Ga1--2]
we get a pairing  $ \, \pi'_\vee \, \colon \big( B^\vee_t \big)' \times
\big(P_t^{\,\prime}\big)^{\!\vee} \hskip-3pt \loongrightarrow \Bbbk[t]
\, $  as claimed,  and a pairing  $ \, \pi'_\vee{\big|}_{t=0} \, \colon
\big( B^\vee_t\big)'{\Big|}_{t=0} \times \big(P_t^{\,\prime}\big)^{\!\vee}
{\Big|}_{t=0} \hskip-3pt \loongrightarrow \Bbbk \, $  as required,
by specialization.
                                             \par
   When  $ \, \pi $  is perfect, by Theorem 4.4 again
$ \pi_{\Cal{R}} $  is perfect as well, and formulas (4.6) hold.
Then if  $ \, \Char(\Bbbk) = 0 \, $  we can apply Theorem 4.11 of
[Ga1--2] to  $ \pi_{\Cal{R}} $  (switching the positions of  $ F_\hbar $
and  $ U_\hbar $  therein),  since  $ P_t^{\,\prime} $  is a QFA and
$ B^\vee_t $  is a QrUEA (notation of [{\it loc.~cit.}]).  The outcome
is that formulas (5.4) hold, which implies that both  $ \, \pi'_\vee
\, $  and  $ \, \pi'_\vee{\big|}_{t=0} \, $  are perfect.
                                              \par
   The final claim about the Hopf algebra case is clear.   \qed
\enddemo

 \vskip1,3truecm

\centerline {\bf \S \; 6 \  The case of group-related Hopf algebras }

\vskip10pt

  {\bf 6.1 The function algebra case.} \, Let  $ G $  be any
algebraic group over the field  $ \Bbbk $.  Let  $ \, F[G]_t :=
\Bbbk[t] \otimes_\Bbbk F[G] \, $:  \, this is trivially a QFA at
$ t \, $  (in the sense of [Ga1--2]), for  $ \, F[G]_t \big/ t \,
F[G]_t = F[G] \, $,  which induces on  $ G $  the trivial Poisson
structure.  Then the cotangent Lie bialgebra of  $ G $  is simply
$ \, \gerg^\times := J \big/ J^2 \, $,  \, where  $ \, J := J_{F[G]}
\equiv \text{\sl Ker}\, \big( \epsilon_{F[G]} \big) \, $,  \, with
trivial Lie bracket and Lie cobracket dual to the Lie bracket of
$ \gerg \, $.   
                                          \par
   We begin by computing  $ \, \widehat{F[G]} := G_{\underline{J}}
\big(F[G]\big) \, $  and its deformation  $ \, {F[G]}_t^\vee :=
\Cal{R}^t_{\underline{J}}\big(F[G]\big) \, $.
                                          \par
   Let  $ \, {\{j_b\}}_{b \in \Cal{S}} \, (\subseteq J \,) \, $  be a
system of parameters of  $ F[G] $,  i.e.~$ \, {\big\{ y_b := j_b \mod
J^2 \big\}}_{b \in \Cal{S}} \, $  is a  $ \Bbbk $--basis  of  $ \,
\gerg^\times := J \Big/ J^2 \, $.  Then, for all  $ \, n \in \N \, $, 
the set  $ \, \big\{\, j^{\,\underline{e}} \mod J^{n+1} \;\big|\;
\underline{e} \in \N^{\Cal{S}}_f \, , \; |\underline{e}| = n \,\big\}
\, $  spans  $ \, J^n \big/ J^{n+1} \, $  over  $ \Bbbk \, $,  where 
$ \, \N^{\Cal{S}}_f := \big\{\, \underline{e} \in \N^{\Cal{S}} \,\big|\,
\underline{e}(s) \! = \! 0 \, \hbox{\ for almost all} \; s \! \in \!
\Cal{S} \,\big\} \, $,  $ \, |\underline{e}| := \sum_{b \in \Cal{S}}
\underline{e}(b) \, $,  and  $ \, j^{\,\underline{e}} := \prod_{b \in
\Cal{S}} j_b^{\,\underline{e}(b)} \, $  (w.r.t.~some fixed total
order on  $ \Cal{S} \, $).  This implies that   
 \vskip3pt
   \centerline{ $ {F[G]}^\vee_t \, = \, \sum_{\underline{e}
\in \N^{\Cal{S}}_f} \Bbbk[t\,] \cdot t^{-|\underline{e}|}
j^{\,\underline{e}} \bigoplus \Bbbk\big[t,t^{-1}\big] \,
J^\infty \, = \, \sum_{\underline{e} \in \N^{\Cal{S}}_f}
\Bbbk[t\,] \cdot {(j^\vee)}^{\,\underline{e}} \bigoplus
\Bbbk\big[t,t^{-1}\big] \, J^\infty $ }
 \vskip1pt
\noindent   where  $ \, J^\infty := \bigcap_{n \in \N} J^n \, $  and
$ \, j^\vee_s := t^{-1} j_s \, $  for all  $ \, s \in \Cal{S} \, $,
\, and also that  $ \, \widehat{F[G]} = G_{\underline{J}}\big(F[G]\big)
\, $  is  $ \Bbbk $--spanned by  $ \, \big\{\, j^{\,\underline{e}} \mod
J^{|\underline{e}| + 1} \;\big|\; \underline{e} \in \N^{\Cal{S}}_f
\,\big\} \, $,  \, so  $ \, \widehat{F[G]} = G_{\underline{J}}
\big(F[G]\big) \, $  is a quotient of  $ \, S(\gerg^\times) \, $.
                                          \par
   Assume for simplicity that  $ \Bbbk $  be perfect.  Let  $ \,
F[[G]] \, $  be the  $ J $--adic completion of  $ \, H = F[G] \, $.
By standard results on algebraic groups (cf.~[DG]) 
         one can choose the above system
 \eject   
\noindent   
of parameters so that  $ \; F[[G]] \, \cong \, \Bbbk \big[
\big[ {\{x_b\}}_{b \in \Cal{S}} \big] \big] \Big/ \Big( \Big\{ x_b^{\,
p^{n(b)}} \Big\}_{\! b \in \Cal{S}} \Big) \; $  
(the algebra of truncated formal power series), as  $ \Bbbk $--algebras  via  $ \, y_b \mapsto x_b
\, $,  \, where  $ p^{n(b)} $  is the nilpotency order of  $ y_b \, $,
\, for all  $ \, b \in \Cal{S} \, $.  Since  $ \, \widehat{F[G]} :=
G_{\underline{J}}\,\big(F[G]\big) = G_{\underline{J}}\,\big(F[[G]]\big)
\cong \Bbbk \big[ {\{ \overline{x}_b \}}_{b \in \Cal{S}} \big] \cong
S(\gerg^\times) \, $,  \, we deduce that  $ \, \widehat{F[G]} \cong
S(\gerg^\times) \bigg/ \! \Big( \Big\{\,\overline{x}^{\,p^{n(x)}}
\Big\}_{\! x \in \Cal{N}(F[G])} \,\Big) =: {S(\gerg^\times)}_{\text{\it
red}} \; $  as algebras, where  $ \Cal{N}\big(F[G]\big) $  is the
nilradical of  $ F[G] $  and  $ \, p^{n(x)} \, $  is the order of
nilpotency of  $ \, x \in \Cal{N}\big(F[G]\big) $.  In addition,
tracking the very definition of the co-Poisson Hopf structure
onto  $ \widehat{F[G]} $  (after our construction) we see that 
$ \, \widehat{F[G]} \cong {S(\gerg^\times)}_{\text{\it red}} \, $
{\sl as co-Poisson Hopf algebras}.  Here the Hopf structure on
$ {S(\gerg^\times)}_{\text{\it red}} $  is induced from the
standard one on  $ S(\gerg^\times) $  (the ideal to be modded out
being a  {\sl Hopf\/}  ideal), and the co-Poisson structure is the
one induced from the Lie cobracket of  $ \gerg^\times $.  Indeed,
for all  $ \, y \in J \, $  and  $ \, \overline{y} := y \mod J^2
\in J \big/ J^2 \subseteq \widehat{F[G]} \, $  we have  $ \, \Delta
\big( \overline{y} \big) \equiv \overline{y} \otimes 1 + 1 \otimes
\overline{y} \in {\widehat{F[G]}}^{\otimes 2} $  (cf.~the proof
of  Theorem 3.5{\it (a)\/}),  so  $ \, \widehat{F[G]} \cong
{S(\gerg^\times)}_{\text{\it red}} \, $  as  {\sl Hopf\/}
algebras too.  The co-Poisson structure on both sides is
uniquely determined by the Lie cobracket of  $ \gerg^\times $,
hence  $ \; \widehat{F[G]} \cong {S(\gerg^\times)}_{\text{\it red}}
\; $  as  {\sl co-Poisson\/}  Hopf algebras too.  Finally, if
$ G $  is a  {\sl Poisson\/}  group then both  $ \widehat{F[G]} $
and  $ S(\gerg^\times) $  are canonically  {\sl bi-Poisson\/}
Hopf algebras, with  $ \, \{\,\ ,\ \}\Big|_{\gerg^\times} =
\, {[\,\ ,\ ]}_{\gerg^\times} \, $,  \; isomorphic to each
other through the above isomorphism.
                                             \par
   Now, for any Lie algebra  $ \gerh $  we can consider  $ \,
\gerh^{{[p\hskip0,7pt]}^\infty} := \Big\{\, x^{{[p\hskip0,7pt]}^n}
\! := x^{p^n} \,\Big\vert\, x \in \gerh \, , n \in \N \,\Big\} \, $,
\, the  {\sl restricted Lie algebra generated by\/  $ \gerh $}  inside
$ U(\gerh) $,  with the  $ p $--operation  given by  $ \, x^{[p
\hskip0,7pt]} := x^p \, $;  \, thus one always has  $ \, U(\gerh)
= \u\big( \gerh^{{[p\hskip0.7pt]}^\infty} \big) \, $.  Then the
previous analysis gives  $ \, {S(\gerg^\times)}_{\text{\it red}}
= \u\left(\gerg^\times_{\text{\it res}}\right) \, $  (for
$ \gerg^\times $  is Abelian),  hence finally  $ \; \widehat{F[G]} \cong
\u\left(\gerg^\times_{\text{\it res}}\right) \, $  as co-Poisson Hopf
algebras, and even as bi-Poisson Hopf algebras whenever  $ G $  is a
Poisson group.
                                             \par
   If  $ \Bbbk $  is not perfect the same analysis applies, but
modifying a bit the previous arguments.
                                             \par
   As for  $ F[G]^\vee $,  it is known (by [Ab], Lemma 4.6.4) that
$ \, F[G]^\vee =F[G] \, $  whenever  $ G $  is finite dimensional
and there exists no  $ \, f \in F[G] \setminus \Bbbk \, $  which
is separable algebraic over  $ \Bbbk \, $.
                                             \par
   It is also interesting to consider  $ {\big( {F[G]}^\vee_t
\big)}' $.  If  $ \, \text{\it Char}\,(\Bbbk) = 0 \, $,  Proposition
4.3 in [Ga1--2] gives  $ \, {\big( {F[G]}^\vee_t \big)}' = F[G]_t
\, $.  If instead  $ \, \text{\it Char}\,(\Bbbk) = p > 0 \, $,  then
the situation might change drastically.  Indeed, if  $ G $  has eight
1   --- i.e.,  $ \, F[G] = \Bbbk\big[ {\{x_i\}}_{i \in \Cal{I}} \big]
\big] \Big/ \big( \big\{ x_i^p \,\big|\, i \in \Cal{I} \,\big\} \big)
\, $  as a  $ \Bbbk $--algebra  ---   then the same analysis as in
the characteristic zero case applies, with a few minor changes,
whence one gets again  $ \, {\big( {F[G]}^\vee_t \big)}' = {F[G]}_t
\, $.  Otherwise, let  $ \, y \in J \setminus \{0\} \, $  be primitive
and such that  $ \, y^p \not= 0 \, $  (for instance, this occurs for
$ \, G \cong \Bbb{G}_a \, $).  Then  $ \, y^p \, $  is primitive as
well, hence  $ \, \delta_n(y^p) = 0 \, $  for each  $ \, n > 1 \, $.
It follows that  $ \, 0 \not= t \, {(y^\vee)}^p \in {\big(
{F[G]}^\vee_t \big)}' \, $,  \, whereas  $ \, t \, {(y^\vee)}^p
\not\in {F[G]}_t \, $,  as follows from our previous description
of  $ {F[G]}^\vee_t $.  Thus  $ \, {\big( {F[G]}^\vee_t \big)}'
\supsetneqq {F[G]}^\vee_t \, $,  \, a  {\sl counterexample to
Proposition 4.3 of [Ga1--2]}.
                                             \par
   What for  $ \widetilde{F[G]} $,  $ {F[G]}' $  and  $ {F[G]}'_t
\, $?  This heavily depends on the group  $ G $  under consideration.
We provide two simple examples, both ``extreme'', and opposite to
each other.
                                             \par
   First, let  $ \, G := \Bbb{G}_a = \hbox{\it Spec}\big(\Bbbk[x]\big)
\, $,  \, so  $ \, F[G] = F[\Bbb{G}_a] = \Bbbk[x] \, $  and  $ \,
{F[\Bbb{G}_a]}_t := \Bbbk[t] \otimes_\Bbbk \Bbbk[x] = \big(\Bbbk[t]
\big)[x] \, $.  Then since  $ \, \Delta(x) := x \otimes 1 + 1
\otimes x \, $  and  $ \, \epsilon(x) = 0 \, $  we find that  $ \,
{F[\Bbb{G}_a]}'_t = \big(\Bbbk[t]\big)[t{}x] \, $  (like in \S 6.3
below), which at  $ \, t = 1 \, $  gives  $ \, {F[\Bbb{G}_a]}' =
F[\Bbb{G}_a] \, $.  Second, let  $ \, G := \Bbb{G}_m = \hbox{\it
Spec} \, \Big( \Bbbk \big[ z^{+1}, z^{-1} \big] \Big) \, $,  \, that
is  $ \, F[G] = F[\Bbb{G}_m] = \Bbbk \big[ z^{+1}, z^{-1} \big] \, $ 
so that  $ \, {F[\Bbb{G}_m]}_t := \Bbbk[t] \otimes_\Bbbk \Bbbk \big[
z^{+1}, z^{-1} \big] = \Bbbk[t] \big[ z^{+1}, z^{-1} \big] \, $. 
Then since  $ \, \Delta \big( z^{\pm 1} \big) := z^{\pm 1} \otimes
z^{\pm 1} \, $  and  $ \, \epsilon \big( z^{\pm 1} \big) = 1 \, $ 
we find  $ \, \delta_n \big( z^{\pm 1} \big) = {\big( z^{\pm 1}
- 1 \big)}^{\otimes n} \, $  for all  $ \, n \in \N \, $,  \,
whence  $ \, {F[\Bbb{G}_m]}'_t = \Bbbk[t] \cdot 1 \, $  and 
$ \, {F[\Bbb{G}_m]}' = \Bbbk \cdot 1 \, $.
                                     \par
   We resume the main facts above in the following statement:

\vskip7pt

\proclaim{Theorem 6.2} \, Let  $ \, H = F[G] \, $  be the function
algebra of an algebraic Poisson group.  Then  $ \, F[G]^\vee =F[G]
\, $  whenever  $ G $  is finite dimensional and there exists no
$ \, f \in F[G] \setminus \Bbbk \, $  which is separable algebraic
over  $ \Bbbk \, $.  In any case,  $ \widehat{F[G]} $  is a
bi-Poisson Hopf algebra, namely
 \vskip-8pt
  $$  \widehat{F[G]} := G_{\underline{J}}\,\big(F[G]\big) \; \cong
\; S(\gerg^\times) \bigg/ \! \Big( \Big\{\, \overline{x}^{\,p^{n(x)}}
\,\Big|\, x \in \Cal{N}\big(F[G]\big) \Big\} \Big) \; \cong \;
\u\left(\gerg^\times_{\text{\it res}}\right)  $$
 \vskip-3pt
\noindent
where  $ \Cal{N}\big(F[G]\big) $  is the nilradical of  $ F[G] \, $, 
$ \, p^{n(x)} \, $  is the order of nilpotency of  $ \, x \in
\Cal{N}\big(F[G]\big) $  and the bi-Poisson Hopf structure of  $ \;
S(\gerg^\times) \bigg/ \! \Big( \Big\{\, \overline{x}^{\,p^{n(x)}}
\,\Big|\, x \in \Cal{N}\big(F[G]\big) \Big\} \Big) \; $  is the
quotient one from  $ \, S(\gerg^\times) \, $.  In particular,
if  $ \, G $  is smooth then  $ \; \widehat{F[G]} \, \cong
\, S(\gerg^\times) \, = \, U(\gerg^\times) \; $.   \qed
\endproclaim

\vskip7pt

  {\bf 6.3 The enveloping algebra case.} \, Let  $ \gerg $  be a
Lie algebra over  $ \Bbbk \, $.  Assume  $ \, p = 0 \, $,  and set
$ \, {U(\gerg)}_t := \Bbbk[t] \otimes_\Bbbk U(\gerg) \, $.  Then
$ {U(\gerg)}_t $  is trivially a QrUEA at  $ t $  (cf.~[Ga1--2]),
for  $ \, {U(\gerg)}_t \big/ t \, {U(\gerg)}_t = U(\gerg) \, $, 
inducing on  $ \gerg $  the trivial Lie cobracket.  The dual
Poisson group is  $ \gerg^\star $,  the topological dual of  $ \gerg $ 
w.r.t.~the weak topology, w.r.t.~addition, with  $ \gerg $  as cotangent
Lie bialgebra and  $ \, F[\gerg^\star] = S(\gerg) \, $.  The Hopf
structure is the standard one, and the Poisson structure is the
Kostant-Kirillov one, induced by  $ \, \{x,y\} := [x,y] \, $
for  $ \, x $,  $ y \in \gerg \, $,  as in \S 1.3.
                                          \par
   Similarly, if  $ \, p > 0 \, $  and  $ \gerg $  is restricted, the
Hopf  $ \Bbbk[t] $--algebra  $ \, {\u(\gerg)}_t := \Bbbk[t] \otimes_\Bbbk
\u(\gerg) \, $  is a QrUEA at  $ t $,  because  $ \, {\u(\gerg)}_t
\big/ t \, {\u(\gerg)}_t = \u(\gerg) \, $,  inducing on  $ \gerg $
the trivial Lie cobracket.  Then the dual Poisson group is again
$ \gerg^\star $,  with cotangent Lie bialgebra  $ \gerg $  and
with  $ \, F[\gerg^\star] = S(\gerg) \, $,  \, the Poisson Hopf
structure being as above.  Recall also that  $ \, U(\gerg) =
\u \big( \gerg^{{[p\hskip0.7pt]}^\infty} \big) \, $  (cf.~\S 1.3).
                                          \par
   First we describe  $ \, {\u(\gerg)}'_t := \Cal{R}^t_{\underline{D}}
\, $,  $ \, {\u(\gerg)}' \, $  and  $ \, \widetilde{\u(\gerg)} \, $,
computing the
        \hbox{$ \delta_\bullet $--filtration  $ \underline{D} \, $
of  $ \u(\gerg) $.}
                                          \par
   By the PBW theorem, once an ordered basis  $ B $  of  $ \gerg $
is fixed  $ \u(\gerg) $  admits as basis the set of ordered monomials
in the elements of  $ B $  whose degree, w.r.t.~each element of  $ B $,
is less than  $ p \, $.  This yields a Hopf algebra filtration of
$ \u(\gerg) $  by the total degree, the so-called  {\sl standard
filtration}.  A straightforward calculation shows that  $ \underline{D} $
coincides with the standard filtration.  This implies  $ \, {\u(\gerg)}'
= \u(\gerg) \, $  and  $ \, {\u(\gerg)}'_t = \langle \tilde{\gerg} \rangle
= \langle t \, \gerg \rangle \, $,  \, where hereafter  $ \, \tilde{\gerg}
:= t \, \gerg \, $,  \, and similarly  $ \, \tilde{x} := t \, x \, $
for all  $ \, x \in \gerg \, $.  Therefore the
           presentation\footnote"${}^\dagger$"{\, Hereafter,  $ \, T_R(M)
\, $,  resp.~$ \, S_R(M) \, $,  \, is the tensor, resp.~symmetric
algebra of an  $ R $--module  $ M $.}
  $ \; {\u(\gerg)}_t \, = \, T_{\Bbbk[t]}(\gerg) \Big/ \big( \big\{\,
x \, y - y \, x - [x,y] \, , \, z^p - z^{[p\hskip0,7pt]}
\;\big\vert\; x, y, z \in \gerg \,\big\} \big) \; $  implies
 \vskip-13pt
  $$  \eqalign{
   {\u(\gerg)}'_t \;  &  = \; T_R\big(\tilde{\gerg}\big) \Big/
\big( \big\{\,\tilde{x} \, \tilde{y} - \tilde{y} \, \tilde{x}
- t \, \widetilde{[x,y]} \, , \, \tilde{z}^{\,p} - t^{\,p-1}
\widetilde{z^{[p\hskip0,7pt]}} \;\big\vert\; x, y, z \in
\gerg \,\big\} \big)  \cr
   \hbox{and}  \qquad  \widetilde{\u(\gerg)} \;  &  = \;
T_\Bbbk(\tilde{\gerg}) \Big/ \big( \big\{\, \tilde{x} \, \tilde{y}
- \tilde{y} \, \tilde{x} \, , \, \tilde{z}^{\,p} \;\big\vert\;
\tilde{x}, \tilde{y}, \tilde{z} \in \tilde{\gerg} \,\big\} \big)
%
%
\; = \; F[\gerg^\star] \Big/ \! \big( \big\{\, z^p
\,\big\vert\; z \in \gerg \,\big\} \big)  \cr }  $$
as  {\sl Hopf\/}  algebras and as  {\sl Poisson Hopf\/}  algebras
respectively.
%
 \eject   
     {\it Therefore  $ \widetilde{\u(\gerg)} $  is the function
algebra of, and  $ {\u(\gerg)}'_t $  is a QFA for, a connected
algebraic Poisson group of dimension 0 and height 1 with cotangent
Lie bialgebra  $ \gerg \, $,  \hbox{\it hence dual to  $ \, \gerg
\, $}}.
                                              \par
   Since  $ \, U(\gerg) = \u\big( \gerg^{{[p\hskip0,7pt]}^\infty}
\big) \, $  the previous analysis yields  $ \, {U(\gerg)}' =
U(\gerg) \, $  and a description of  $ {U(\gerg)}'_t $  and
$ \widetilde{U(\gerg)} $,  \, in particular  $ \; \widetilde{U(\gerg)}
\, = \, F \! \left[ {\big( \gerg^{{[p\hskip0,7pt]}^\infty} \big)}^{\!
\star} \right] \hskip-2pt \bigg/ \hskip-2pt \Big( {\big\{\, z^p \,
\big\}}_{z \in \gerg^{{[p\hskip0,7pt]}^\infty}} \Big) \; $.
                                          \par
   Furthermore,  $ \, {\u(\gerg)}'_t = \langle
\tilde{\gerg} \rangle \, $  implies that  $ \, J_{{\u(\gerg)}'_t}
\, $  is generated (as an ideal) by  $ \, \tilde{\gerg} \, $,  \,
so  $ \, t^{-1} J_{{\u(\gerg)}'_t} \, $ is generated by  $ \, \gerg
\, $,  \, thus  $ \, \big( {\u(\gerg)}'_t \big)^{\!\vee} := \sum_{n
\geq 0} {\big( t^{-1} J_{{\u(\gerg)}'_t} \big)}^n = \sum_{n \geq 0}
\gerg^n = {\u(\gerg)}_t \; $.
                                          \par
   When  $ \, \text{\it Char}\,(\Bbbk) = 0 \, $  and we look at
$ U(\gerg) $,  the like argument applies:  $ \underline{D} $
coincides with the standard filtration of  $ U(\gerg) $  given
by total degree, via the PBW theorem.  This and definitions imply
$ \, {U(\gerg)}' = U(\gerg) \, $  and  $ \, {U(\gerg)}'_t = \langle
\tilde{\gerg} \rangle = \langle t \, \gerg \rangle \, $,  so that
the presentation  $ \; {U(\gerg)_t} = T_{\Bbbk[t]}(\gerg) \! \Big/
\! \big( {\big\{\, x \, y - y \, x - [x,y] \,\big\}}_{x, y, z \in
\gerg\,} \big) \; $  yields
 \vskip-6pt
  $$  \eqalign{
   {U_t(\gerg)}' \;  &  = \; T_{\Bbbk[t]}(\tilde{\gerg}) \Big/
\big( {\big\{ \tilde{x} \, \tilde{y} - \tilde{y} \, \tilde{x}
- t \, \widetilde{[x,y]} \,\big|\, \tilde{x}, \tilde{y} \in
\tilde{\gerg}\,} \,\big\} \big)  \cr
   \text{and}  \qquad  \widetilde{U(\gerg)} \;  &  = \;
T_\Bbbk(\tilde{\gerg}) \Big/ \big( \big\{\, \tilde{x} \,
\tilde{y} - \tilde{y} \, \tilde{x} \;\big\vert\; \tilde{x},
\tilde{y} \in \tilde{\gerg} \,\big\} \big) \, = \, S(\gerg)
\, = \, F[\gerg^\star]  \cr }  $$
 \vskip-4pt
\noindent   as  {\sl Hopf\/}  and as  {\sl Poisson Hopf\/}  algebras
(as predicted by  Theorem 2.2{\it (c)\/}  in [Ga1-2]).  Furthermore,
$ \, J_{{U(\gerg)}'_t} \, $  is generated by  $ \, \tilde{\gerg}
\, $:  \, therefore  $ \, t^{-1} J_{{U(\gerg)}'_t} \, $  is
generated by  $ \, \gerg \, $,  \, whence eventually  $ \, {\big(
{U(\gerg)}'_t \big)}^{\!\vee} := \sum_{n \geq 0} {\big(\, t^{-1}
J_{{U(\gerg)}'_t} \,\big)}^n = \sum_{n \geq 0} \gerg^n =
{U(\gerg)}_t \; $.
                                          \par
   What for the functor  $ \, {(\ )}^\vee \, $?  This heavily
depends on the  $ \gerg $  we start from.   
                                          \par
   First assume  $ \, \hbox{\it Char}\,(\Bbbk) = 0 \, $.  Let
$ \, \gerg_{(1)} := \gerg \, $,  $ \, \gerg_{(k)} := \big[ \gerg,
\gerg_{(k-1)} \big] \, $  ($ k \in \N_+ $),  be the  {\sl lower
central series\/}  of  $ \gerg \, $.  Set  $ \, \gerg_{\langle
k \rangle} := \gerg_{(k)} \big/ \gerg_{(k+1)} \, $  for all  $ \,
k \in \N \, $,  \, and  $ \, \gerg_{\text{\it gr}} := \oplus_{k \in \N}
\gerg_{\langle k \rangle} \, $,  \, the latter being a  {\sl graded\/}
Lie algebra in a natural way; set also  $ \, \gerg_{(\infty)} := \cap_{k
\in \N} \gerg_{(k)} \, $.  Recall that  $ \, {U(\gerg)}^\vee_t :=
\sum\nolimits_{n \geq 0} t^{-n} J^n \, $  for  $ \, J := J_{U_t(\gerg)}
\, $.  Pick subsets  $ \, B_k \, \big( \! \subset \! \gerg_{(k)} \big)
\, $  such that  $ \, B_k \mod \gerg_{(k+1)} \, $  be a  $ \Bbbk $--basis 
of  $ \, \gerg_{\langle k \rangle} \, $  (for all  $ \, k \in \N \, $), 
give any total order to  $ \, B := \cup_{k \in \N} B_k \, $  and set 
$ \, \partial(b) := k \, $  iff  $ \, b \in B_k \, $.  Applying the
PBW theorem to this ordered basis of  $ \gerg \, $,  we get that  $ \,
J^n \big/ J^{n+1} \, $  has basis  $ \, \big\{\, b_1^{e_1} b_2^{e_2}
\cdots b_s^{e_s} \mod J^{n+1} \;\big|\; s \in \N \, , b_r \in B \, ,
\, \sum_{r=1}^s b_r \, \partial(b_r) = n \,\big\} \, $.  Then one
finds that  $ \, {U(\gerg)}^\vee_t \, $  is generated by  $ \,
\big\{\, t^{-1} b_r \;\big|\; b_r \in B_1 \,\big\} \cup \big(
\cup_{n \in \N} t^{-n} \gerg_{(\infty)} \big) \, $,  \, while 
$ \; {U(\gerg)}^\vee \, = \, U(\gerg) \Big/ \gerg_{(\infty)}
\, U(\gerg) \; $  and  $ \; \widehat{U(\gerg)} \, \cong \,
U(\gerg_{\text{\it gr}}) \; $  as  {\sl graded\/}  Hopf algebras.
                                          \par
   Now assume  $ \, \hbox{\it Char}\,(\Bbbk) = p > 0 \, $.
Let  $ \, \gerg_n := \Big\langle \bigcup_{m \, p^k \geq n}
{(\gerg_{(m)})}^{[p^k]} \Big\rangle \, $  for all  $ \, n
\in \N_+ \, $  (where  $ \langle X \rangle $  denotes the
Lie subalgebra of  $ \gerg $  generated by  $ X \, $):  \,
then  $ \, {\big\{ \gerg_n \big\}}_{n \in \N_+} \, $  is
{\sl the  $ p $--lower  central series of\/}  $ \gerg \, $,
which is a  {\sl strongly central series\/}  of  $ \gerg \, $. 
Set  $ \, \gerg_{[k]} := \gerg_k \big/ \gerg_{k+1} \, $  for
all  $ \, k \in \N \, $,  \, and  $ \, \gerg_{\text{\it p-gr}} :=
\oplus_{k \in \N} \, \gerg_{[k]} \, $,  \, the latter being a  {\sl
graded\/}  restricted Lie algebra in a natural way; set also  $ \,
\gerg_\infty := \cap_{k \in \N} \, \gerg_k \, $.  Now definitions give
$ \, \gerg_n \subseteq J^n \, $  for all  $ \, n \in \N \, $,  where
$ \, J := J_{\u(\gerg)} \, $.  More precisely, we can proceed as above,
taking suitable lifts  $ B_k $  of bases of each  $ \gerg_{[k]} \, $. 
Then the (restricted) PBW theorem for  $ \, \u(\gerg) \, $  implies
that  $ \, J^n \big/ J^{n+1} \, $  has  $ \Bbbk $--basis  the set of
ordered monomials  $ \, x_{i_1}^{e_1} x_{i_2}^{e_2} \cdots x_{i_s}^{e_s}
\, $  (with the  $ x_i $'s  in the union of the  $ B_k $'s)  such that 
$ \, \sum_{r=1}^s e_r \partial(x_{i_r}) = n \, $,  where  $ \, \partial
(x_{i_r}) \in \N \, $  is uniquely determined by the condition  $ \,
x_{i_r} \in \gerg_{\partial(x_{i_r})} \setminus \gerg_{\partial
(x_{i_r}) + 1} \, $.  This yields a description of  $ \underline{J}
\, $,  hence of  $ {\u(\gerg)}^\vee $,  $ {\u_{\,t}(\gerg)}^\vee $,  and
$ \widehat{\u(\gerg)} \, $:  \, in particular,  $ \; {\u(\gerg)}^\vee
\, = \, \u(\gerg) \Big/ \gerg_\infty \, \u(\gerg) \; $  and  $ \,
\widehat{\u(\gerg)} \, \cong \, \u(\gerg_{\text{\it p-gr}}) \, $ 
as  {\sl graded\/}  Hopf algebras.
                                               \par
   We collect the main results above in the following statement:

\vskip7pt

\proclaim{Theorem 6.4}
                                      \hfill\break
   \indent   (a) \, Let  $ \, \Char(\Bbbk) = 0 \, $,  and  $ \,
\gerg \, $  be a Lie bialgebra.  Then  $ \, {U(\gerg)}' = U(\gerg)
\, $,  \, and  $ \, \widetilde{U(\gerg)} := G_{\underline{D}}\,\big(
U(\gerg) \big) \, $  is a bi-Poisson Hopf algebra, namely  $ \;
\widetilde{U(\gerg)} \, \cong \, S(\gerg) \, = \, F[\gerg^\star]
\; $  (with notation of\/ \S 1), where the bi-Poisson Hopf structure on
$ S(\gerg^\star) $  is the canonical one (see\/ \S 1.3).  On the other
hand,  $ \; {U(\gerg)}^\vee \, = \, U(\gerg) \Big/ \gerg_{(\infty)}
\, U(\gerg) \; $  and  $ \; \widehat{U(\gerg)} \, \cong \,
U(\gerg_{\text{\it gr}}) \; $  as  {\sl graded\/}  Hopf algebras.
                                      \hfill\break
   \indent   (b) \, Let  $ \, \Char(\Bbbk) = p > 0 \, $,  and
$ \, \gerg \, $  be a restricted Lie bialgebra.  Then  $ \,
{\u(\gerg)}' = \u(\gerg) \, $,  \, and  $ \, \widetilde{\u(\gerg)}
:= G_{\underline{D}}\,\big(\u(\gerg)\big) \, $  is a bi-Poisson
Hopf algebra, namely  $ \; \widetilde{\u(\gerg)} \, \cong \, S(\gerg)
\Big/ \big( \big\{ x^p \,\big|\, x \in \gerg \big\} \big) \, = \,
F[G^\star] \; $  (with notation of\/ \S 1).  Here the bi-Poisson Hopf
structure on  $ \, S(\gerg) \Big/ \big( \big\{ x^p \,\big|\, x \in \gerg
\big\} \big) \, $  is induced by the canonical one on  $ S(\gerg) $ 
(see\/ \S 1.3) and  $ \, G^\star $  is a connected algebraic Poisson
group of dimension 0 and height 1 whose cotangent Lie bialgebra is 
$ \gerg \, $.  On the other hand,  $ \; {\u(\gerg)}^\vee \, = \,
\u(\gerg) \Big/ \gerg_\infty \, \u(\gerg) \; $  and  $ \;
\widehat{\u(\gerg)} \, \cong \, \u(\gerg_{\text{\it p-gr}}) \; $
as  {\sl graded\/}  Hopf algebras.   \qed
\endproclaim

\vskip3pt

{\it $ \underline{\hbox{{\it Remark}}} $:\/}  For any given Lie algebra 
$ \gerg \, $,  the group-scheme theoretic version of Lie's Third Theorem
claims the existence of a connected algebraic group-scheme  {\sl of
height 1\/}  having  $ \gerg $  as tangent Lie algebra.  Part  {\it
(b)\/}  of Theorem 6.4 gives a Poisson-dual result:  $ \, G^\star =
\gerg^\star \, $  is an algebraic Poisson group-scheme of height 1
   \hbox{having  $ \gerg $  as  {\sl cotangent\/}  Lie algebra.}   

\vskip7pt

  {\bf 6.5 The hyperalgebra case.} \, Let  $ G $  be an algebraic
group, which for simplicity we assume to be finite-dimensional.
The hyperalgebra associated to  $ G $  is  $ \, \hyp(G) := {\big(
{F[G]}^\circ \big)}_\epsilon = \big\{\, \phi \in {F[G]}^\circ
\,\big|\, \phi({\germ_e}^{\!n}) = 0 \, , \, \forall \; n \gg 0
\,\big\} \, $,  \, that is the irreducible component of the  {\sl
dual\/}  Hopf algebra  $ \, {F[G]}^\circ \, $  containing  $ \,
\epsilon = \epsilon_{\!{}_{F[G]}} \, $.  This is a Hopf subalgebra
of  $ {F[G]}^\circ $,  connected and cocommutative.  There is a
Hopf algebra morphism  $ \, \Phi : U(\gerg) \longrightarrow \hyp(G)
\, $;  \, if  $ \, p = 0 \, $  then  $ \Phi $  is an isomorphism,
so  $ \hyp(G) $  identifies to  $ U(\gerg) $;  \, if  $ \, p >
0 \, $  then  $ \Phi $  factors through  $ \u(\gerg) $  and the
induced morphism  $ \, \overline{\Phi} : \u(\gerg) \longrightarrow
\hyp(G) \, $  is injective, so that  $ \u(\gerg) $  identifies with
a Hopf subalgebra of  $ \hyp(G) \, $.   
     \hbox{Now we study  $ {\hyp(G)}' $,  $ {\hyp(G)}^\vee $,
$ \widetilde{\hyp(G)} \, $,  $ \widehat{\hyp(G)} \, $.}
                                               \par
   As  $ \hyp(G) $  is connected, we have  $ \; \hyp(G) = {\hyp(G)}'
\; $.  Now,  Theorem 3.5{\it (b)\/}  gives  $ \, \widetilde{\hyp(G)}
= F[\varGamma\,] \, $  for some connected algebraic Poisson group
$ \varGamma \, $.  Theorem 6.2 yields  $ \, \widehat{F[G]} \cong
S(\gerg^*) \bigg/ \! \Big( \Big\{\, \overline{x}^{\,p^{n(x)}}
\Big\}_{x \in \Cal{N}_{F[G]}} \Big) = \u \bigg(\! P \bigg( S(\gerg^*)
\bigg/ \! \Big( \Big\{\, \overline{x}^{\,p^{n(x)}} \Big\}_{x \in
\Cal{N}_{F[G]}} \Big) \!\bigg) \!\bigg) \! = \u \Big(\! \big(
\gerg^* \big)^{p^\infty} \Big) \, $,  \, with  $ \, \big( \gerg^*
\big)^{p^\infty} \! := \text{\sl Span}\, \Big( \Big\{\, x^{p^n}
\,\Big\vert\; x \in \gerg^* \, , n \in \N \,\Big\} \Big) \subseteq
\widehat{F[G]} \, $,  \, and noting that  $ \, \gerg^\times \! =
\gerg^* \, $.  On the other hand, exactly like for  $ U(\gerg) $
and  $ \u(\gerg) $,  respectively in case  $ \, \Char(\Bbbk) = 0 \, $
and  $ \, \Char(\Bbbk) > 0 \, $,  \, the filtration  $ \underline{D} $
of  $ \hyp(G) $  is the natural filtration given by the order of
differential operators.  This implies
     \hbox{$ \; \hyp(G)^{\,\prime}_t := \, \Cal{R}_t\big(\hyp(G)\big)
= \big\langle \big\{\, t^n x^{(n)} \,\big|\; x \in \gerg \, , n \in
\N \,\big\} \big\rangle \, $,}
\, where hereafter  $ x^{(n)} $  is the  $ n $--th  divided power
or the  $ n $--th  binomial coefficient of  $ \, x \in \gerg \, $
(such  $ x^{(n)} $'s  generate  $ \hyp(G) $  as an algebra).  Then
one easily checks that the graded Hopf pairing between  $ \, \hyp(G)^{\,
\prime}_t \Big/ t \, \hyp(G)^{\,\prime}_t = \widetilde{\hyp(G)} =
F[\varGamma] \, $  and  $ \, \widehat{F[G]} \, $  in Theorem 3.7
is perfect, from which we deduce that      
     \hbox{$ \varGamma $  has cotangent   
Lie bialgebra isomorphic to  $ \, \Big( \hskip-2pt \big( \gerg^*
\big)^{p^\infty} \Big)^{\!*} \, $.}   

\vskip1,3truecm

   \centerline {\bf \S \; 7 The Crystal duality Principle }

\vskip10pt

{\bf 7.1 The Crystal Duality Principle.} \, We can finally motivate
the expression ``Crystal Duality Principle''.  In short, we gave
functorial recipes to get, out from  {\sl any\/}  Hopf algebra  $ H $, 
four Hopf algebras of Poisson-geometrical type arranged in couples,
namely  $ \, \big(\, \widehat{H}, \widetilde{H} \,\big) = \big(
\U(\gerg_-), F[G_+] \big) \, $  and  $ \, \Big(\! \big(H^\vee_t\big)'
{\big|}_{t=0}, \big(H'_t \big)^{\!\vee}{\big|}_{t=0} \Big) = \big(
F[K_+], \U(\gerk_-) \big) \,  $,  \, hence four Poisson-geometrical
symmetries  $ G_+ \, $,  $ \gerg_- \, $,  $ K_+ \, $  and  $ \gerk_-
\, $:  this is the ``Principle''.  The word ``Crystal'' reminds the
fact that the first couple, namely  $ \big( \widehat{H}, \widetilde{H}
\,\big) $,  is obtained via a  {\sl crystallization process\/} 
(cf.~\S\S 2.15 and 3.6).  Finally, ``Duality'' witnesses
that if  $ \, \Char(\Bbbk) = 0 \, $  then the link between the two
couples of special Hopf algebras is Poisson duality (see  Theorem
5.4{\it (c)\/}),  in that  $ \, K_+ = G_-^\star \, $  and  $ \,
\gerk_- = \gerg_+^{\,\times} \, $.  Moreover, in  {\sl any\/} 
characteristic, when  $ H $  is a Hopf algebra of (Poisson-)geometrical
type applying the crystal functor leading to Hopf algebras
of dual type the result is ruled   
      \hbox{by Poisson duality (see Theorems 6.2 and 6.4).}

\vskip7pt

   {\bf 7.2 The CDP as corollary of the GQDP.} \, The construction
of Drinfeld-like functors passes through the application of the Global
Quantum Duality Principle (=GQDP in the sequel).  In this section we
briefly outline how  {\sl the whole CDP can be obtained as a corollary
of the GQDP\/}  (but for some minor details); see also [Ga1--2], \S 5.
                                                \par
   For any  $ \, H \in \HA \, $,  \, let  $ \, H_t := H[t] \equiv
\Bbbk[t] \otimes_\Bbbk H \, $.  Then  $ H_t $  is a torsionless
Hopf algebra over  $ \Bbbk[t] $,  hence one of those to which the
constructions in [Ga1--2] can be applied.  In particular, we can
act on it with Drinfeld's functors considered therein: these give
quantum groups, namely a quantized (restricted) universal enveloping
algebra (=QrUEA) and a quantized function algebra (=QFA).
Straightforward computations show that the QrUEA is just
$ \, H^\vee_t := \Cal{R}^t_{\underline{J}}(H) \, $,  \, and
the QFA is  $ \, H'_t := \Cal{R}^t_{\underline{D}}(H) \, $,  \,
with  $ \, \widehat{H} \cong H^\vee_t{\big|}_{t=0} \, $  and  $ \,
\widetilde{H} \cong H'_t{\big|}_{t=0} \; $.  It follows that all
properties of  $ \widehat{H} $  and  $ \widetilde{H} $  spring out as
special cases of the results proved in [Ga1--2] for Drinfeld's functors,
but for their being graded.  Similarly, the fact that  $ H' $  be a
Hopf subalgebra of  $ H $  follows from the fact that  $ H'_t $  itself
is a Hopf algebra and  $ \, H' = H'_t{\big|}_{t=1} \, $.  Instead,
$ H^\vee $  is a quotient Hopf algebra of  $ H $  because  $ H^\vee_t $
is a Hopf algebra, hence  $ \, \overline{H^\vee_t} := H^\vee_t \Big/
\bigcap_{n \in \N} t^n H^\vee_t \, $  is a Hopf algebra, and finally
$ \, H^\vee = \overline{H^\vee_t}{\big|}_{t=1} \, $.  The fact that
$ H'_t $ and  $ H^\vee_t $  be regular 1-parameter deformations of
$ H' $  and  $ H^\vee $  is then clear by construction.  Finally,
the parts of the CDP dealing with Poisson duality are direct
consequences of the like items in the GQDP applied to  $ H'_t $
and to  $ H^\vee_t $  (but for  Theorem 6.4{\it (b)\/}).  The
     \hbox{cases of (co)augmented (co)algebras or bialgebras can
be easily treated the same.}   

\vskip1,3truecm

\centerline {\bf \S \; 8 \  The Crystal Duality Principle on group
algebras }

\vskip10pt

  {\bf 8.1 Group-related algebras.} \, In this section,  $ G $  is
any abstract group.  For any commutative unital ring  $ \A \, $,
by  $ \, \A[G] \, $  we mean the group algebra of  $ G $  over
$ \A \, $  and, when  $ G $  is  {\sl finite},  we denote by  $ \,
A_\A(G) := {\A[G]}^* \, $  (the linear dual of  $ \A[G] \, $)  the
function algebra of  $ G $  over  $ \A \, $.  Hereafter  $ \Bbbk $
will be a field with  $ \, p := \Char(\Bbbk) \, $,  \, and
$ \, R := \Bbbk[t] \, $.   
                                          \par
   Recall that  $ \, H := \A[G] \, $  admits  $ G $  itself
as a special basis, with Hopf algebra structure given by
$ \, g \cdot_{{}_H} \gamma := g \cdot_{{}_G} \gamma \, $,  $ \,
1_{{}_H} := 1_{{}_G} \, $,  $ \, \Delta(g) := g \otimes g \, $,
$ \, \epsilon(g) := 1 \, $,  $ \, S(g) := g^{-1} \, $,  \, for all
$ \, g, \gamma \in G \, $.  Dually,  $ \, H := A_\A(G) \, $  has
basis  $ \, \big\{ \varphi_g \,\big|\, g \! \in \! G \big\} \, $  dual
to the basis  $ G $  of  $ \A[G] \, $,  \, with  $ \, \varphi_g(\gamma)
:= \delta_{g,\gamma} \, $  for all  $ \, g, \gamma \in G \, $.  Its Hopf
algebra structure is given by  $ \, \varphi_g \cdot \varphi_\gamma :=
\delta_{g,\gamma} \varphi_g \, $,  $ \, 1_{{}_H} := \sum_{g \in G}
\varphi_g \, $,  $ \, \Delta(\varphi_g) := \sum_{\gamma \cdot \ell = g}
\varphi_\gamma \otimes \varphi_\ell \, $,  $ \, \epsilon(\varphi_g) :=
\delta_{g,1_G} \, $,  $ \, S(\varphi_g) := \varphi_{g^{-1}} \, $,  \,
            for all  $ \, g $,\break
\noindent
$ \gamma \in G \, $.
  \hbox{Thus  $ \, {\Bbbk[G]}_t = R \otimes_\Bbbk \Bbbk[G] = R[G] \, $,
$ \, {A_\Bbbk[G]}_t = R \otimes_\Bbbk A_\Bbbk[G] = A_R[G] \, $.
First we have}

\vskip7pt

\proclaim{Theorem 8.2} $ \, {{\Bbbk[G]}_t}^{\!\prime} = R \cdot 1 \, $,
$ \; {\Bbbk[G]}' = \Bbbk \cdot 1 \; $  and  $ \; \widetilde{\Bbbk[G]}
= \Bbbk \cdot 1 = F\big[\{*\}\big] \, $.
\endproclaim

\demo{Proof} The claim follows easily from the formula  $ \, \delta_n(g)
= {(g \! - \! 1)}^{\otimes n} $,  for  $ \, g \in G $,  $ n \in \N \, $.
\qed
\enddemo

\vskip7pt

  {\bf 8.3  $ {\Bbbk\boldkey{[}\boldkey{G}\,
\boldkey{]}}_{\boldkey{t}}^{\,\boldsymbol\vee} $,  $ \,
{\Bbbk\boldkey{[}\boldkey{G}\,\boldkey{]}}^{\boldsymbol\vee} $ 
and  $ \widehat{\Bbbk\boldkey{[}\boldkey{G}\,\boldkey{]}} $  and
the dimension subgroup problem.}  \, In contrast with the triviality
result in Theorem 8.2, things are more interesting for  $ \, {R[G]}^\vee
\! = {{\Bbbk[G]}_t}^{\!\!\vee} \, $,  $ \, {\Bbbk[G]}^\vee $  and  $ \,
\widehat{\Bbbk[G]} \, $.  Note, however, that since  $ \, \Bbbk[G] \, $ 
is cocommutative the induced Poisson cobracket on  $ \, \widehat{\Bbbk[G]}
\, $  is trivial, and the Lie cobracket of $ \, \gerk_G := P \Big(
\widehat{\Bbbk[G]} \Big) \, $  is trivial as well.
                                          \par
   Studying  $ {\Bbbk[G]}^\vee $  and  $ \widehat{\Bbbk[G]} $
amounts to study the filtration  $ {\big\{ J^n \big\}}_{n \in \N}
\, $,  \, with  $ J := \text{\sl Ker}\,(\epsilon_{{}_{\Bbbk[G]}}) $,
\, which is a classical topic.  Indeed, for  $ \, n \! \in \! \N \, $
let  $ \, \Cal{D}_n(G) := \big\{\, g \in G \,\big|\, (g \! - \! 1) \in
J^n \,\big\} \, $: \, this is a characteristic subgroup of  $ G $,  called
{\sl the  $ n^{\text{th}} $  dimension subgroup of  $ G \, $}.  All these
form a filtration inside  $ G \, $:  \, characterizing it in terms of
$ G $  is the  {\sl dimension subgroup problem}.  For group algebras
over fields, it is completely solved (see [Pa], Ch.~11, \S 1, and
[HB], and references therein); this also gives a description of  $ \,
\big\{ J^n \big\}_{n \in \N_+} \, $.  Thus we find ourselves within
the domain of classical group theory.  Now we use the results which
solve the dimension subgroup problem to deduce a description of 
$ {\Bbbk[G]}^\vee $,  $ \widehat{\Bbbk[G]} $  and  $ {\Bbbk[G]_t}^{\!
\vee} $.  Later on we'll get from this a description
      \hbox{of  $ \big( {\Bbbk[G]_t}^{\!\vee} \big)' $  and its
semiclassical limit too.}
                                          \par
   By construction,  $ J \, $  has  $ \Bbbk $--basis  $ \, \big\{ \eta_g
\,\big|\; g \! \in \! G \setminus \{1_{{}_G}\} \big\} \, $,  \, where
$ \, \eta_g := (g-1) \, $.  Then  $ \, {\Bbbk[G]}^\vee \, $  is
generated by  $ \, \big\{\, \eta_g \! \mod J^\infty \,\big|\; g \in
G \setminus \{1_{{}_G}\} \big\} \, $,  \, and  $ \, \widehat{\Bbbk[G]}
\, $  by  $ \, \big\{\, \overline{\,\eta_g} \;\big|\; g \! \in \! G
\setminus \{1_{{}_G}\} \big\} \, $.  Hereafter  $ \, \overline{x}
:= x \mod J^{n+1} \, $  for all  $ \, x \in J^n \, $,  \, that is
$ \overline{x} $  is the element in  $ \widehat{\Bbbk[G]} $  which
corresponds to  $ \, x \in \Bbbk[G] \, $.  Moreover,  $ \, \overline{g}
= \overline{\,1 + \eta_g} = \overline{1} \, $  for all  $ \, g \in G \, $;
\, also,  $ \; \Delta \big(\overline{\,\eta_g}\big) = \overline{\,\eta_g}
\otimes \overline{g} + 1 \otimes \overline{\,\eta_g} = \overline{\,\eta_g}
\otimes 1 + 1 \otimes \overline{\,\eta_g} \; $:  \; thus  $ \overline{\,
\eta_g} $  is primitive, so  $ \, \big\{\,\overline{\,\eta_g} \;\big|\;
g \! \in \! G \setminus \{1_{{}_G}\} \big\} \, $  generates  $ \,
\gerk_G := P \Big( \widehat{\Bbbk[G]} \Big) \, $.

\vskip7pt

  {\bf 8.4 The Jennings-Hall theorem.} \, The description of 
$ \Cal{D}_n(G) $  is given by the Jennings-Hall theorem, which
we now recall.  The involved construction strongly depends on
whether  $ \, p := \Char(\Bbbk) \, $  is zero or not, so we
shall distinguish these two cases.
                                          \par
   First assume  $ \, p = 0 \, $.  Let  $ \, G_{(1)} := G \, $,  $ \,
G_{(k)} := (G,G_{(k-1)}) \, $  ($ k \in \N_+ $),  form the  {\sl lower
central series\/}  of  $ G \, $;  hereafter  $ (X,Y) $  is the
commutator subgroup of  $ G $  generated by the set of commutators
$ \, \big\{ (x,y) := x \, y \, x^{-1} y^{-1} \,\big|\, x \in X, y
\in Y \big\} \, $.  Then let  $ \, \sqrt{G_{(n)}} := \big\{ x \in
G \,\big|\, \exists \, s \in \N_+ : x^s \in G_{(n)} \big\} \, $
for all  $ n \in \N_+ \, $:  these form a descending series of
characteristic subgroups in  $ G $,  such that each composition
factor  $ \, A^G_{(n)} := \sqrt{G_{(n)}} \Big/ \! \sqrt{G_{{(n+1)}}}
\, $  is torsion-free Abelian.  Therefore  $ \, \L_0(G) :=
\bigoplus_{n \in \N_+} A^G_{(n)} \, $  is a graded Lie ring, with
Lie bracket  $ \, \big[ \overline{g}, \overline{\ell} \,\big] :=
\overline{(g,\ell\,)} \, $  for all  {\sl homogeneous}  $ \overline{g} $,
$ \overline{\ell} \in \L_0(G) \, $,  \, with obvious notation.  It is
easy to see that the map  $ \; \Bbbk \otimes_\Z \L_0(G) \longrightarrow
\gerk_G \, $,  $ \, \overline{g} \mapsto \overline{\eta_g} \, $,  \;
is an epimorphism  {\sl of graded Lie rings\/}:  \, therefore  {\sl
the Lie algebra  $ \, \gerk_G \, $  is a quotient of  $ \; \Bbbk
\otimes_\Z \L_0(G) \, $};  in fact, the above is an isomorphism
(see below).  We shall use notation  $ \, \partial(g) := n \, $  for
all  $ \, g \in \sqrt{G_{(n)}} \, \setminus \sqrt{G_{(n+1)}} \; $.
                                     \par
   For each $ \, k \in \N_+ \, $  pick in  $ A^G_{(k)} $  a subset
$ \overline{B}_k $  which is a  $ \Bbb{Q} $--basis  of  $ \, \Bbb{Q}
\otimes_\Z A^G_{(k)} \, $.  For each  $ \, \overline{b} \in
\overline{B}_k \, $,  \, choose a fixed  $ \, b \in \sqrt{G_{(k)}}
\, $  such that its coset in  $ A^G_{(k)} $  be  $ \overline{b} $,  \,
and denote by  $ \, B_k \, $  the set of all such elements  $ b \, $.
Let  $ \, B := \bigcup_{k \in \N_+} B_k \, $:  \, we call such a set
{\sl t.f.l.c.s.-net\/}  (\,=\,``torsion-free-lower-central-series-net'')
on  $ G $.  Clearly  $ \, B_k = \Big( B \cap \sqrt{G_{(k)}} \,\Big)
\setminus \Big( B \cap \sqrt{G_{(k+1)}} \,\Big) \, $  for all
$ k \, $.  By an  {\sl ordered t.f.l.c.s.-net\/}  we mean a
t.f.l.c.s.-net  $ B $  which is totally ordered in such a way that:
{\it (i)\/}  if  $ \, a \in B_m \, $,  $ \, b \in B_n \, $,  $ \, m
< n \, $,  \, then  $ \, a \preceq b \, $;  \, {\it (ii)\/}  for each
$ k $,  every non-empty subset of  $ B_k $  has a greatest element.
An ordered t.f.l.c.s.-net always exists.
                                     \par
   Now assume instead  $ \, p > 0 \, $.  Starting from the lower
central series  $ \, {\big\{G_{(k)}\big\}}_{k \in \N_+} $,  define
$ \, G_{[n]} := \prod_{k p^\ell \geq n} {(G_{(k)})}^{p^\ell} \; $  for
all  $ \, n \in \N_+ \, $  (hereafter, for any group  $ \varGamma $ 
we denote  $ \varGamma^{p^e} $  the subgroup generated by  $ \, \big\{
\gamma^{p^e} \,\big|\, \gamma \! \in \! \varGamma \,\big\} \, $).  This
gives another strongly central series  $ \, {\big\{G_{[n]}\big\}}_{n
\in \N_+} $  in  $ G $,  \, with the additional property that  $ \,
{(G_{[n]})}^p \leq G_{[n+1]} \, $  for all  $ n \, $,  \, called
{\sl the  $ p $--lower  central series of\/}  $ G \, $.  Then  $ \,
\Cal{L}_p(G) := \oplus_{n \in \N_+} G_{[n]} \big/ G_{[n+1]} \, $  is
a graded restricted Lie algebra over  $ \, \Z_p := \Z \big/ p \, \Z \, $,
\, with operations  $ \, \overline{g} + \overline{\ell} := \overline{g
\cdot \ell} \, $,  $ \, \big[\overline{g},\overline{\ell}\,\big] :=
\overline{(g,\ell\,)} \, $,  $ \, \overline{g}^{\,[p\,]} := \overline{g^p}
\, $,  \, for all  $ \, g $,  $ \ell \in G \, $  (cf.~[HB], Ch.~VIII,
\S 9).  Like before, we consider the map  $ \; \Bbbk \otimes_{\Z_p}
\Cal{L}_p(G) \longrightarrow \gerk_G \, $,  $ \, \overline{g} \mapsto
\overline{\eta_g} \, $,  \; which now is an epimorphism  {\sl of graded
restricted Lie  $ \Z_p $--algebras},  whose image spans  $ \gerk_G $
over  $ \Bbbk \, $.  Thus  {\sl  $ \, \gerk_G \, $  is a quotient of 
$ \; \Bbbk \otimes_{\Z_p} \Cal{L}_p(G) \, $};  \, in fact, the above
is an iso\-morphism (see below).  Finally, we introduce notation  $ \,
d(g) := n \, $  for all  $ \, g \in G_{[n]} \setminus G_{[n+1]} \, $.
                                      \par
   For each  $ \, k \in \N_+ \, $,  choose a  $ \Z_p $--basis
$ \overline{B}_k $  of the  $ \Z_p $--vector  space  $ \, G_{[k]} \big/
G_{[k+1]} \, $.  For each  $ \, \overline{b} \in \overline{B}_k \, $,
\, fix  $ \, b \in G_{[k]} \, $  such that  $ \, \overline{b} = b \,
G_{[k+1]} \, $,  \, and let  $ \, B_k \, $  be the set of all such
elements  $ b \, $.  Let  $ \, B := \bigcup_{k \in \N_+} B_k \, $:
\, such a set will be called a  {\sl  $ p $-l.c.s.-net\/}  (=
``$ p $-lower-central-series-net''; the terminology in [HB] is
``$ \kappa $-net'') on  $ G $.  Of course  $ \, B_k = \big( B \cap
G_{[k]} \big) \setminus \big( B \cap G_{[k+1]} \big) \, $  for all
$ k \, $.  By an  {\sl ordered  $ p $-l.c.s.-net\/}  we mean a
$ p $-l.c.s.-net  $ B $  which is totally ordered in such a way
that:  {\it (i)\/}  if  $ \, a \in B_m \, $,  $ \, b \in B_n \, $,
$ \, m < n \, $,  \, then  $ \, a \preceq b \, $;  \, {\it (ii)\/}
for each  $ k $,  every non-empty subset of  $ B_k $  has a greatest
element (like for  $ \, p = 0 \, $).  Again,  $ p $-l.c.s.-nets  do
exist.
                                      \par
   We can now describe each  $ \Cal{D}_n(G) $,  hence also each graded
summand  $ \, J^n \big/ J^{n+1} \, $  of  $ \, \widehat{\Bbbk[G]} \, $, 
%
%
in terms of a fixed ordered t.f.l.c.s.-net or 
$ p $-l.c.s.-net.  To unify notations, set  $ \, G_n := G_{(n)} \, $, 
$ \, \theta(g) := \partial(g) \, $  if  $ \, p \! = \! 0 \, $,  \, and 
$ \, G_n := G_{[n]} \, $,  $ \, \theta(g) := d(g) \, $  if  $ \, p >
0 \, $,  \, set  $ \, G_\infty := \bigcap_{n \in N_+} \! G_n \, $, 
\, let  $ \, B := \bigcup_{k \in \N_+} B_k \, $  be an ordered
t.f.l.c.s.-net or  $ p $-l.c.s.-net according to whether  $ \,
p = 0 \, $  or  $ \, p > 0 \, $,  \, and set  $ \, \ell(0) :=
+ \infty \, $  and  $ \, \ell(p) := p \, $  for  $ \, p > 0
\, $.  The key result we need is

\vskip5pt

\noindent   {\it  $ \underline{\hbox{\sl Jennings-Hall theorem}} $
(cf.~[HB], [Pa] and references therein).  Let  $ \, p:= \Char(\Bbbk)
\, $.
                                       \par
   (a) \, For all  $ \, g \in G \, $,  $ \; \eta_g
\in J^n \Longleftrightarrow g \in \! G_n \, $.  Therefore
$ \, \Cal{D}_n(G) = G_n \; $  for all  $ \, n \in \N_+ \, $.
                                       \par
   (b) \, For any  $ \, n \in \N_+ \, $,  the set of ordered monomials
 \vskip-19pt
  $$  \Bbb{B}_n \, := \, \Big\{\, {\overline{\,\eta_{b_1}}}^{\;e_1}
\cdots {\overline{\,\eta_{b_r}}}^{\;e_r} \;\Big|\; b_i \in B_{d_i}
\, , \; e_i \in \N_+ \, , \; e_i < \ell(p) \, , \; b_1 \precneqq
\cdots \precneqq b_r \, , \; {\textstyle \sum}_{i=1}^r e_i \, d_i
= n \,\Big\}  $$
 \vskip-12pt
\noindent   is a\/  $ \Bbbk $--basis  of  $ \, J^n \big/ J^{n+1} \, $,
\, and  $ \; \Bbb{B} \, := \, \{1\} \cup \bigcup_{n \in \N} \Bbb{B}_n
\; $  is a\/  $ \Bbbk $--basis  of  $ \; \widehat{\Bbbk[G]} \, $.
                                          \par
   (c) \;  $ \big\{\, \overline{\,\eta_b} \,\;\big|\;
b \in B_n \big\} \; $  is a  $ \Bbbk $--basis  of the  $ n $--th
graded summand  $ \, \gerk_G \cap \big( J^n \big/ J^{n+1} \big) \, $
of the graded restricted Lie algebra\/  $ \gerk_G \, $,  \, and
$ \, \big\{\, \overline{\,\eta_b} \,\;\big|\; b \in B \,\big\} \; $
is a\/  $ \Bbbk $--basis  of\/  $ \gerk_G \, $.
                                          \par
   (d) \;  $ \big\{\, \overline{\,\eta_b} \,\;\big|\;
b \in B_1 \big\} \; $  is a minimal set of generators of
the (restricted) Lie algebra\/  $ \gerk_G \, $.
                                          \par
   (e) \; The map  $ \; \Bbbk \otimes_\Z \L_p(G) \longrightarrow
\gerk_G \, $,  $ \, \overline{g} \mapsto \overline{\,\eta_g} \, $,
\, is an isomorphism of graded restricted Lie algebras.  Therefore
$ \; \widehat{\Bbbk[G]} \, \cong \, \U \big( \Bbbk \otimes_\Z \L_p(G)
\big) \; $  as Hopf algebras (notation of \S 1.3).   
                                          \par
   (f) \;  $ J^\infty \, = \, \hbox{\sl Span} \big( \big\{\,
\eta_g \,\big|\, g \in G_\infty \,\big\} \big) \, $,  \, whence\/
$ \; {\Bbbk[G]}^\vee \cong \, \bigoplus_{\overline{g} \in G/G_\infty}
\! \Bbbk \cdot \overline{g} \; \cong \, \Bbbk \big[ G \big/ G_\infty
\big] \; $.   \qed}

\vskip5pt

   Recall that  $ A\big[x,x^{-1}\big] \, $,  for any  $ A \, $,  has
$ A $--basis  $ \, \big\{ {(x \!-\! 1)}^n x^{-[n/2]} \,\big|\, n \in
\N \big\} \, $,  \, where  $ [q] $  is the integer part of  $ \, q
\in \Bbb{Q} \, $.  Then from Jennings-Hall theorem and definitions
we deduce   

\vskip7pt

\proclaim{Proposition 8.5} \, Let  $ \; \chi_g := t^{-\theta(g)} \eta_g
\, $,  \, for all  $ \, g \in \{G\} \setminus \{1\} \, $.  Then
 \vskip-16pt
  $$  \displaylines{
   {\Bbbk[G]}_{\,t}^{\,\vee}  = \,  \Big( {\textstyle
\bigoplus_{\Sb  b_i \in B, \; 0 < e_i < \ell(p)  \\
                r \in \N, \; b_1 \precneqq \cdots \precneqq b_r  \endSb}}
R \cdot \chi_{b_1}^{\;e_1} \, b_1^{\,-[e_1\!/2]} \cdots \chi_{b_r}^{\;
e_r} \, b_r^{\,-[e_r/2]} \Big) \,{\textstyle
\bigoplus}\, R\big[t^{-1}\big] \cdot J^\infty  \; =   \hfill  \cr
   {} \;  = \,  \Big( {\textstyle
\bigoplus_{\Sb  b_i \in B, \; 0 < e_i < \ell(p)  \\
                r \in \N, \; b_1 \precneqq \cdots \precneqq b_r  \endSb}}
R \cdot \chi_{b_1}^{\,e_1} \, b_1^{\,-[e_1\!/2]} \cdots \chi_{b_r}^{\,
e_r} \, b_r^{\,-[e_r/2]} \Big) \,{\textstyle \bigoplus}\,
\Big( \hskip1pt {\textstyle \sum_{\gamma \in G_\infty}}
\hskip-0pt R\big[t^{-1}\big] \cdot \eta_\gamma \Big) \; ;
\hskip-2,5pt  \cr }  $$
 \vskip-7pt
\noindent
If  $ \, J^\infty \! = \! J^n $  for some  $ \, n \! \in \! \N $
(iff  $ \, G_\infty \! = G_n $)  we can drop the factors  $ \,
b_1^{-[e_1\!/2]}, \dots, b_r^{-[e_r/2]} \, . \, \square $
\endproclaim

\vskip7pt

   {\bf 8.6 Poisson groups from  $ \Bbbk\boldkey{[}\boldkey{G}\,
\boldkey{]}$.} \, The previous discussion attached to  $ G $  and
$ \Bbbk $  the (maybe restricted) Lie algebra  $ \gerk_G \, $.  By
Jennings-Hall theorem, this is just the scalar extension of the Lie
ring  $ \L_p $  associated to  $ G $  via the central series of the
$ G_n $'s.  In particular, the functor  $ \, G \mapsto \gerk_G \, $ 
is one considered since a long time in group theory.  
                                           \par   
   By  Theorem 5.4{\it (a)\/}  we know that  $ \, {\big(
{\Bbbk[G]}_{\,t}^{\,\vee} \big)}'{\Big|}_{t=0} = F \big[
\varGamma_G \big] \, $  for some connected Poisson group 
$ \varGamma_G \, $,  \, of dimension zero and height 1 if 
$ \, p > 0 \, $.  This defines a functor  $ \, G \mapsto
\varGamma_G \, $:  in particular,  $ \varGamma_G $  {\sl
is a new invariant for abstract groups}.
                                            \par
   The description of  $ {\Bbbk[G]}_{\,t}^{\,\vee} $  in
Proposition 8.5 leads us to an explicit description of
$ {\big( {\Bbbk[G]}_{\,t}^{\,\vee} \big)}' $,  \, hence of
$ \, {\big( {\Bbbk[G]}_{\,t}^{\,\vee} \big)}'{\Big|}_{t=0} \! =
F\big[\varGamma_G\big] \, $  and  $ \varGamma_G $  too.  Indeed, direct
inspection gives  $ \, \delta_n\big(\chi_g\big) = t^{(n-1) \, \theta(g)}
\chi_g^{\;\otimes n} \, $,  \, so  $ \, \psi_g := t \, \chi_g = t^{1 -
\theta(g)} \eta_g \in {\big( {\Bbbk[G]}_{\,t}^{\,\vee} \big)}' \setminus
t \, {\big( {\Bbbk[G]}_{\,t}^{\,\vee} \big)}' \, $  for each  $ \, g \in
G \setminus G_\infty \, $.  Instead for  $ \, \gamma \in G_\infty $
we have  $ \, \eta_\gamma \in J^\infty \, $,  which implies also  $ \,
\eta_\gamma \in {\big( {\Bbbk[G]}_{\,t}^{\,\vee} \big)}' \, $,  \,
and even  $ \, \eta_\gamma \in \bigcap_{n \in \N} t^n {\big(
{\Bbbk[G]}_{\,t}^{\,\vee} \big)}' \, $.  Therefore  $ {\big(
{\Bbbk[G]}_{\,t}^{\,\vee} \big)}' $  is generated by  $ \, \big\{\,
\psi_g \;\big|\; g \in G \setminus \{1\} \big\} \cup \big\{\, \eta_\gamma
\,\big|\, \gamma \in G_\infty \big\} \, $.  Moreover,  $ \, g = 1 +
t^{\theta(g)-1} \psi_g \in {\big( {\Bbbk[G]}_{\,t}^{\,\vee} \big)}' \, $ 
for every  $ \, g \in G \setminus G_\infty \, $,  \, and  $ \, \gamma =
1 + (\gamma - 1) \in 1 + J^\infty \subseteq {\big( {\Bbbk[G]}_{\,t}^{\,
\vee} \big)}' \, $  for  $ \, \gamma \in G_\infty \, $.  This and the
previous analysis, along with Proposition 8.5, prove next result,
which in turn is the basis for Theorem 8.8 below.

\vskip7pt

\proclaim{Proposition 8.7}
  $$  \eqalign{
   {\big( {\Bbbk[G]}_{\,t}^{\,\vee} \big)}' \;  &  = \;
\Big( {\textstyle \bigoplus_{\Sb  b_i \in B, \; 0 < e_i < \ell(p)  \\
                r \in \N, \; b_1 \precneqq \cdots \precneqq b_r  \endSb}}
R \cdot \psi_{b_1}^{\;e_1} \, b_1^{\,-[e_1\!/2]} \cdots \psi_{b_r}^{\;
e_r} \, b_r^{\,-[e_r/2]} \Big) \,{\textstyle \bigoplus}\,
R\big[t^{-1}\big] \cdot J^\infty  \; =   \hfill  \cr
   {}  &  = \;  \Big( {\textstyle
\bigoplus_{\Sb  b_i \in B, \; 0 < e_i < \ell(p)  \\
                r \in \N, \; b_1 \precneqq \cdots \precneqq b_r  \endSb}}
R \cdot \psi_{b_1}^{\,e_1} \, b_1^{\,-[e_1\!/2]} \cdots \psi_{b_r}^{\,
e_r} \, b_r^{\,-[e_r/2]} \Big) \,{\textstyle \bigoplus}\,
\Big( \hskip1pt {\textstyle \sum_{\gamma \in G_\infty}}
\hskip-0pt R\big[t^{-1}\big] \cdot \eta_\gamma \Big) \; .  \cr }  $$
In particular,  $ \; {\big( {\Bbbk[G]}_{\,t}^{\,\vee} \big)}' =
{\Bbbk[G]}_{\,t} \; $  {\sl if and only if  $ \; G_2 = \{1\} =
G_\infty \; $.}  If in addition  $ \, J^\infty \! = \! J^n \, $ 
for some  $ \, n \in \N $  (iff  $ \, G_\infty = G_n $),  then
we can drop the factors $ \, b_1^{\, -[e_1\!/2]}, \dots,
b_r^{\,-[e_r/2]} \, $.   \qed
\endproclaim

\vskip7pt

\proclaim{Theorem 8.8}\, Let  $ \; x_g := \psi_g \mod t \; {\big(
{\Bbbk[G]}_{\,t}^{\,\vee} \big)}' \, $,  $ \, z_g := g \mod t \;
{\big( {\Bbbk[G]}_{\,t}^{\,\vee} \big)}' \, $  for all  $ \, g \not=
1 \, $,  \, and  $ \, B_1 := \big\{\, b \in \! B \,\big\vert\, \theta(b)
= 1 \big\} \, $,  $ \, B_> := \big\{\, b \in \! B \,\big\vert\, \theta(b)
> 1 \big\} \, $.
                                        \hfill\break
   \indent   (a) \, If  $ \, p = 0 \, $,  \, then  $ \, F \big[
\varGamma_G \big] = {\big( {\Bbbk[G]}_{\,t}^{\,\vee} \big)}'
{\Big|}_{t=0} $  is a  {\sl polynomial/Laurent polynomial algebra},
namely  $ \, F \big[ \varGamma_G \big] = \Bbbk \big[ {\big\{
{z_b}^{\!\pm 1} \big\}}_{b \in B_1} \!\! \cup {\{x_b\}}_{b \in B_>}
\big] \, $,  \; the  $ z_b $'s  being group-like and the  $ x_b $'s 
being primitive.  In particular  $ \, \varGamma_G \cong \big(
\Bbb{G}_m^{\,\times B_1} \big) \times \big( \Bbb{G}_a^{\,
\times B_>} \big) \, $ as algebraic groups.
                                        \hfill\break
   \indent   (b) \, If  $ \, p > 0 \, $,  \, then  $ \, F \big[
\varGamma_G \big] = {\big( {\Bbbk[G]}_{\,t}^{\,\vee} \big)}'
{\Big|}_{t=0} $  is a  {\sl truncated polynomial/Laurent polynomial
algebra}, namely  $ \, F \big[ \varGamma_G \big] = \, \Bbbk \big[
{\big\{ {z_b}^{\!\pm 1} \big\}}_{b \in B_1} \!\! \cup {\{x_b\}}_{b
\in B_>} \big] \Big/ \! \big( \{ z_b^{\,p} - 1 \}_{b \in B_1} \cup
\{ x_b^{\,p} \}_{b \in B_>} \big) \, $,  \, the  $ z_b $'s  being
group-like and the  $ x_b $'s  being primitive.  In particular 
$ \, \varGamma_G \cong \big( {{\boldsymbol\mu}_p}^{\! \times B_1}
\big) \times \big( {{\boldsymbol\alpha}_p}^{\!\times B_>} \big)
\, $  as algebraic groups of dimension zero and height 1.   
                                        \hfill\break
   \indent   (c) \, The Poisson group  $ \varGamma_G $  has cotangent
Lie bialgebra  $ \gerk_G \, $,  that is  $ \, \text{\sl coLie}\,
(\varGamma_G) = \gerk_G \, $.
\endproclaim

\demo{Proof}  {\it (a)} \, Definitions give  $ \, \partial(g\,\ell\,)
\geq \partial(g) + \partial(\ell\,) \, $  for all  $ \, g, \ell \in G
\, $,  \, so that  $ \; [\psi_g,\psi_\ell\big] =
%
%
t^{1 - \partial(g)
- \partial(\ell) + \partial((g,\ell))} \, \psi_{(g,\ell)} \, g \, \ell
\in t \cdot {\big( {\Bbbk[G]}_{\,t}^{\,\vee} \big)}' \, $,  \; which
proves (directly) that  $ \, {\big( {\Bbbk[G]}_{\,t}^{\,\vee} \big)}'
{\Big|}_{t=0} \, $  is commutative.  Moreover, the relation  $ \, 1 =
g^{-1} \, g = g^{-1} \, \big(1 + t^{\partial(g)-1} \psi_g \big) \, $
(for any  $ \, g \in G \, $)  yields  $ \, z_{g^{-1}} = {z_g}^{\!-1}
\, $,  \, and  $ \, z_{g^{-1}} = 1 \, $  iff  $ \, \partial(g) > 1
\, $.  Note also that  $ \, J^\infty \equiv 0 \!\! \mod t \,
{\big( {\Bbbk[G]}_{\,t}^{\,\vee} \big)}' \, $  and  $ \, g
= 1 + t^{\partial(g) - 1} \psi_g \equiv 1 \mod t \, {\big(
{\Bbbk[G]}_{\,t}^{\,\vee} \big)}' \, $  for  $ \, g \in G \setminus
G_\infty \, $,  \, and also  $ \, \gamma = 1 + (\gamma - 1) \in 1 +
J^\infty \equiv 1 \mod t \, {\big( {\Bbbk[G]}_{\,t}^{\,\vee} \big)}'
\, $  for  $ \, \gamma \in G_\infty \, $.  Then Proposition 8.7 gives
  $$  F\big[\varGamma_G\big]  \; = \; {\big( {\Bbbk[G]}_{\,t}^{\,\vee}
\big)}'{\Big|}_{t=0} = \,  \Big( {\textstyle
\bigoplus_{\hskip-3pt   \Sb  b_i \in B_1, \; a_i \in \Z  \\
                s \in \N, \; b_1 \precneqq \cdots \precneqq b_s  \endSb}}
\hskip-2pt  \Bbbk \cdot z_{b_1}^{\,a_1} \cdots z_{b_s}^{\,a_s} \Big)
\,{\textstyle \bigotimes}\; \Big( {\textstyle
\bigoplus_{\hskip-1pt   \Sb  b_i \in B_>, \; e_i \in \N_+  \\
                r \in \N, \; b_1 \precneqq \cdots \precneqq b_r  \endSb}}
\hskip-2pt  \Bbbk \cdot x_{b_1}^{\,e_1} \cdots x_{b_r}^{\,e_r} \Big)  $$  
which means that  $ F\big[\varGamma_G\big] $  is a polynomial/Laurent
polynomial algebra, as claimed.  
                                                \par    
   Again definitions imply   $ \, \Delta(z_g) = z_g \otimes z_g \, $ 
for all  $ \, g \in G \, $  and  $ \, \Delta(x_g) = x_g \otimes 1
+ 1 \otimes x_g \, $  iff  $ \, \partial(g) > 1 \, $.  Thus the 
$ z_b $'s  are group-like and the  $ x_b $'s  are primitive, as claimed.  
                                          \par  
   {\it (b)} \, The definition of  $ d $  implies  $ \, d(g\,\ell\,)
\geq d(g) + d(\ell\,) \, $  ($ g, \ell \in G $),  \, whence we get
$ \; [\psi_g,\psi_\ell] \, = \, t^{1 - d(g) - d(\ell) + d((g,\ell))}
\, \psi_{(g,\ell)} \, g \, \ell \, \in \, t \cdot {\big(
{\Bbbk[G]}_{\,t}^{\,\vee} \big)}' \, $,  \, proving that
$ \, {\big( {\Bbbk[G]}_{\,t}^{\,\vee} \big)}'{\Big|}_{t=0} \, $
is commutative.  In addition  $ \; d(g^p) \geq p \; d(g) \, $,  \,
so  $ \; \psi_g^{\;p} = t^{\, p \, (1 - d(g))} \, \eta_g^{\;p} =
t^{\, p - 1 + d(g^p) - p\,d(g)} \, \psi_{g^p} \in t \cdot {\big(
{\Bbbk[G]}_{\,t}^{\,\vee} \big)}' \, $,  \, whence  $ \, {\big(
\psi_g^{\;p}{\big|}_{t=0} \big)}^p = 0 \, $  inside  $ \, {\big(
{\Bbbk[G]}_{\,t}^{\,\vee} \big)}'{\Big|}_{t=0} \! = F \big[
\varGamma_G \big] \, $,  \, which proves that  $ \varGamma_G $
has dimension 0 and height 1.  Finally,  $ \; b^p = {(1 + \psi_b)}^p
= 1 + \psi_b^{\,p} \equiv 1 \mod t \, {\big( {\Bbbk[G]}_{\,t}^{\,\vee}
\big)}' \, $  for all  $ \, b \in B_1 \, $,  \, so  $ \, b^{-1} \equiv
b^{p-1} \!\! \mod t \, {\big( {\Bbbk[G]}_{\,t}^{\,\vee} \big)}' \, $.
   \hbox{Thus, letting  $ \, x_g := \psi_g \!\! \mod t \; {\big(
{\Bbbk[G]}_{\,t}^{\,\vee} \big)}' $  (for  $ \, g \! \not= \! 1 $), 
we get}
  $$  F\big[\varGamma_G\big] = {\big( {\Bbbk[G]}_{\,t}^{\,\vee} \big)}'
{\Big|}_{t=0} \!  =  \Big( {\textstyle \bigoplus_{\Sb  b_i \in B_1, \;
0 < e_i < p  \\   
                s \in \N, \; b_1 \precneqq \cdots \precneqq b_s  \endSb}}
\Bbbk \cdot z_{b_1}^{\,e_1} \cdots z_{b_s}^{\,e_s} \Big)
 {\textstyle \bigotimes}\,
     \Big( {\textstyle \bigoplus_{\Sb  b_i \in B_>, \; 0 < e_i < p  \\
                r \in \N, \; b_1 \precneqq \cdots \precneqq b_r  \endSb}}
\Bbbk \cdot x_{b_1}^{\,e_1} \cdots x_{b_r}^{\,e_r} \Big)  $$   
just like for  {\it (a)},  and also taking care that  $ \, z_b =
x_b + 1 \, $  and  $ \, z_b^{\,p} = 1 \, $  for  $ \, b \in B_1 \, $.
Therefore  $ \, {\big( {\Bbbk[G]}_{\,t}^{\,\vee} \big)}'{\Big|}_{t=0}
\, $  is a truncated polynomial/Laurent polynomial algebra as claimed.
The properties of the  $ x_b $'s  and the  $ z_b $'s  w.r.t.~the Hopf
     \hbox{structure are then proved like for \!  {\it (a)\/}  again.}
                                                 \par
   {\it (c)} \, The augmentation ideal  $ \, \germ_e \, $  of  $ \,
{\big( {\Bbbk[G]}_{\,t}^{\,\vee} \big)}'{\Big|}_{t=0} = F \big[
\varGamma_G \big] \, $  is generated by  $ \, {\{ x_b \}}_{b \in B}
\, $.  Then  $ \; t^{-1} \, [\psi_g,\psi_\ell\big] \, = \, t^{\,
\theta((g,\ell)) - \theta(g) - \theta(\ell)} \, \psi_{(g,\ell)}
\, \big( 1 + t^{\, \theta(g) - 1} \psi_g \,\big) \, \big( 1 + t^{\,
\theta(\ell) - 1} \psi_\ell \,\big) \, $  by the previous computation,
whence at  $ \, t = 0 \, $  one has  $ \; \big\{ x_g \, , x_\ell \big\}
\, \equiv \, x_{(g,\ell)} \mod \, \germ_e^{\,2} \, $  if  $ \, \theta
\big( (g,\ell\,) \big) = \theta(g) + \theta(\ell\,) \, $,  and  $ \;
\big\{ x_g \, , x_\ell \big\} \, \equiv \, 0 \mod \, \germ_e^{\,2} \, $ 
if  $ \, \theta \big( (g,\ell\,) \big) > \theta(g) + \theta(\ell\,)
\, $.  This means that the cotangent Lie bialgebra  $ \, \germ_e
\Big/ \germ_e^{\,2} \, $  of  $ \varGamma_G \, $  is isomorphic
to  $ \gerk_G \, $,  \, as claimed.   \qed   
\enddemo   

\vskip2pt

   {\bf Remark 8.9:}  Theorem 8.8 yields functorial recipes to attach
to each abstract group  $ G $  and each field  $ \Bbbk $  a connected
Abelian algebraic Poisson group over\/  $ \Bbbk $,  namely  $ \; G
\mapsto \varGamma_G \equiv K_G^\star \, $,  \; with  $ \, \text{\sl
coLie}\,(\varGamma_G) = \gerk_G \, $.  {\sl Every such  $ \varGamma_G $ 
(for given  $ \Bbbk $)  is then an invariant of  $ G \, $},  a new one
to the author's knowledge.  Yet this invariant is perfectly  {\sl
equivalent\/}  to the well-known invariant  $ \gerk_G $  (over
the same  $ \Bbbk $).  In fact,  $ \, \varGamma_{G_1} \! \cong
\varGamma_{G_2} \, $  implies  $ \, \gerk_{G_1} \! \cong \gerk_{G_2}
\, $,  \, whereas  $ \, \gerk_{G_1} \! \cong \gerk_{G_2} \, $ 
implies that  $ \varGamma_{G_1} $  and  $ \varGamma_{G_2} $  are
isomorphic as algebraic groups, by  Theorem 8.8{\it (a--b)},  and bear
isomorphic Poisson structures, by Theorem 8.8{\it (c)},  so  $ \,
\varGamma_{G_1} \! \cong \varGamma_{G_2} \, $  as Poisson groups.   

\vskip7pt

   {\bf 8.10 The case of  $ \boldkey{A}_\Bbbk\boldkey{(}\boldkey{G}
\boldkey{)} \, $.} \, Let's now  $ G $  be a  {\sl finite\/}  group,
$ \A $  any commutative unital ring, and  $ \Bbbk \, $,  $ \, R :=
\Bbbk[t] \, $  be as before.  By definition  $ \, A_\A(G) = {\A[G]}^*
\, $,  \, hence  $ \, \A[G] = {A_\A(G)}^* \, $,  \, and we have
a natural perfect Hopf pairing  $ \, A_\A(G) \times \A[G]
\longrightarrow \A \, $.  Our first result is

\vskip7pt

\proclaim{Theorem 8.11}  $ \, {A_\Bbbk(G)}_t^\vee \! = R \cdot
1 \oplus R \big[ t^{-1} \big] \, J = {\big( {A_\Bbbk(G)}_t^\vee
\big)}' \, $,  $ \, {A_\Bbbk(G)}^\vee = \Bbbk \!\cdot\! 1 \, $, 
$ \, \widehat{A_\Bbbk(G)} = {A_\Bbbk(G)}_t^\vee{\Big|}_{t=0} \!
= \, \Bbbk \cdot 1 = \, \U(\boldkey{0}) \; $  and  $ \; {\big(
{A_\Bbbk(G)}_t^\vee \big)}'{\Big|}_{t=0} \! = \, \Bbbk \cdot 1
= F\big[\{*\}\big] \; $.
\endproclaim

\demo{Proof} By construction  $ \, J := \text{\sl Ker}\,(\epsilon_{\!
{}_{A_\Bbbk(G\,)}}) \, $  has  $ \Bbbk $--basis  $ \, \big\{ \varphi_g
\big\}_{g \in G \setminus \{1_{{}_G}\}} \cup \big\{ \varphi_{1_G} \! -
1_{\!{}_{A_\Bbbk(G\,)}} \big\} \, $,  \, and since  $ \, \varphi_g =
{\varphi_g}^{\!2} \, $  for all  $ g $  and  $ \, {(\varphi_{1_G} \!
- \! 1)}^2 = -(\varphi_{1_G} \!-\! 1) \, $  we have  $ \, J = J^\infty
\, $,  \, so  $ \, {A_\Bbbk(G)}^\vee = \Bbbk \!\cdot\! 1 \, $  and
$ \, \widehat{A_\Bbbk(G)} = \Bbbk \!\cdot\! 1 \, $.  Similarly, the set
$ \, \big\{ t^{-1} \varphi_g \big\}_{g \in G \setminus \{1_{{}_G}\}}
\cup \big\{ t^{-1} (\varphi_{1_G} \! - \! 1_{\!{}_{A_R(G\,)}}) \big\}
\, $  generates  $ \, {A_\Bbbk(G)}_t^\vee \, $,  \, since  $ \,
{A_\Bbbk(G)}_t = A_R(G) \, $.  Moreover,  $ \, J = J^\infty \, $
implies  $ \, t^{-n} J \subseteq {A_\Bbbk(G)}_{\,t}^{\,\vee} \, $  for
all  $ \, n \, $,  \, so  $ \, {A_\Bbbk(G)}_t^\vee = R \, 1 \oplus
R[t^{-1}] J \, $.  Then  $ \, J_{{A_\Bbbk(G)}_t^\vee} = R \big[ t^{-1}
\big] J \subseteq t \, {A_\Bbbk(G)}_t^\vee \, $,  \, which implies 
$ \, {\big( {A_\Bbbk(G)}_t^\vee \big)}' = {A_\Bbbk(G)}_t^\vee \, $: 
\, in particular  $ \, {\big( {A_\Bbbk(G)}_t^\vee \big)}'{\Big|}_{t=0}
\! = {A_\Bbbk(G)}_t^\vee{\Big|}_{t=0} \! = \Bbbk \cdot 1 \, $,  \, as
claimed.   \qed
\enddemo

%
%
 \vskip-1pt

   {\bf 8.12 Poisson groups from  $ \boldkey{A}_\Bbbk\boldkey{(}
\boldkey{G}\boldkey{)} \, $.} \, Now we look at  $ {A_\Bbbk(G)}_t^{\,
\prime} \, $,  $ {A_\Bbbk(G)}' $  and  $ \widetilde{A_\Bbbk(G)} \, $.  By
construction  $ {A_\Bbbk(G)}_t $  and  $ {\Bbbk[G]}_t $  are in perfect
Hopf pairing, and are free  $ R $--modules  of  {\sl finite rank}.  In
this case Theorem 4.4 yields  $ \, {A_\Bbbk(G)}_t^{\,\prime} = {\big(
{\Bbbk[G]}_t^\vee \big)}^\bullet \, $  (see (4.6)), hence we have  $ \,
{A_\Bbbk(G)}_t^{\,\prime} = {\big( {\Bbbk[G]}_t^\vee \big)}^\bullet
= {\big( {\Bbbk[G]}_t^\vee \big)}^* \, $:  thus  $ \, {A_\Bbbk(G)}_t^{\,
\prime} \, $  is the  {\sl dual\/}  Hopf algebra to  $ {\Bbbk[G]}_t^\vee
\, $.  Then from Proposition 8.5 we can deduce an explicit description of 
$ {A_\Bbbk(G)}_t^{\,\prime} $,  whence also of  $ \big( {A_\Bbbk(G)}_t^{\,
\prime} \big)^{\!\vee} $.  By Theorem 3.7 there are a perfect filtered Hopf
pairing  $ \, {\Bbbk[G]}^\vee \times {A_\Bbbk(G)}' \longrightarrow \Bbbk
\, $  and a perfect graded Hopf pairing  $ \, \widehat{\Bbbk[G]} \times
\widetilde{A_\Bbbk(G)} \! \longrightarrow \Bbbk \, $:  \, thus  $ \,
{A_\Bbbk(G)}' \! \cong \! {\big( {\Bbbk[G]}^\vee \big)}^* \, $  as
filtered Hopf algebras and  $ \, \widetilde{A_\Bbbk(G)} \cong {\big(
\widehat{\Bbbk[G]}\big)}^* \, $  as graded Hopf algebras.
                                         \par
   If  $ p = 0 \, $  then  $ J = J^\infty $,  \, as each  $ \, g \in
G \, $  has finite order and  $ \, g^n = 1 \, $  implies  $ \, g \in
G_\infty \, $.  Then  $ \, {\Bbbk[G]}^\vee \! = \Bbbk \cdot 1 =
\widehat{\Bbbk[G]} \, $,  \, so  $ {A_\Bbbk(G)}' = \Bbbk \cdot 1 =
\widetilde{A_\Bbbk(G)} \, $.  If  $ \, p > 0 \, $  instead, this
analysis gives  $ \, \widetilde{A_\Bbbk(G)} = {\big( \widehat{\Bbbk[G]}
\big)}^* = {\big( \u(\gerk_G) \big)}^* = F[K_G] \, $,  \, where  $ \,
K_G  \, $  is a connected Poisson group of dimension 0, height 1 and
tangent Lie bialgebra  $ \gerk_G \, $.  Thus we get:   

\vskip7pt

\proclaim{Theorem 8.13}
                                                \par
   (a) \, There is a second functorial recipe to attach to each 
{\sl finite}  abstract group a connected algebraic Poisson group
of dimension zero and height 1 over any field\/  $ \Bbbk $  with
$ \, \Char(\Bbbk) > 0 \, $,  namely  $ \; G \mapsto K_G :=
\text{\it Spec}\,\Big( \widetilde{A_\Bbbk(G)} \Big) \, $,
\; with  $ \; \text{\it Lie}\,(K_G) = \, \gerk_G =
\text{\sl coLie}\,(\varGamma_G) \; $.
                                                \par
   (b) \, If  $ \, p:= \Char(\Bbbk) > 0 \, $,  then  $ \;
\big({A_\Bbbk(G)}_t^{\,\prime}\big)^{\!\vee}\Big|_{t=0}
\! = \, \u \big( \gerk_G^{\,\times} \big) = S\big(\gerk_G^{\,
\times}\big) \Big/ \big( \big\{ x^p \,\big|\, x \in \gerk_G^{\,
\times} \big\} \big) \; $.
\endproclaim

\demo{Proof} Claim  {\it (a)\/}  is the outcome of the discussion
above.  Part  {\it (b)\/}  instead requires an explicit description
of  $ \, \big({A_\Bbbk(G)}_t^{\,\prime}\big)^{\!\vee} $.  Since
$ \, {A_\Bbbk(G)}_t^{\,\prime} \cong \big( {\Bbbk[G]}_{\,t}^{\,\vee}
\big)^* \, $,  \, from Proposition 8.5 we get  $ \; {A_\Bbbk(G)}_t^{\,
\prime} = \, \Big(\! \bigoplus_{\Sb  b_i \in B, \; 0 < e_i < p  \\
          r \in \N, \; b_1 \precneqq \cdots \precneqq b_r  \endSb}
R \cdot \rho_{b_1, \dots, b_r}^{e_1, \dots, e_r} \Big) \; $
where each  $ \, \rho_{b_1, \dots, b_r}^{e_1, \dots, e_r} \, $  is
defined by
                      \hfill\break
 \vskip1pt
   \centerline{$ \Big\langle \, \rho_{b_1, \dots, b_r}^{e_1, \dots,
e_r} \; , \; \chi_{\beta_1}^{\;\varepsilon_1} \, \beta_1^{\,
-[\varepsilon_1\!/2]} \cdots \chi_{\beta_s}^{\;\varepsilon_s} \,
\beta_s^{\,-[\varepsilon_s/2]} \,\Big\rangle \; = \; \delta_{r,s}
\, \prod_{i=1}^r \delta_{b_i,\beta_i} \delta_{e_i, \varepsilon_i} $}
 \vskip4pt
\noindent
(for all  $ \, b_i, \beta_j \in B \, $  and  $ \, 0 < e_i,
\varepsilon_j < p \, $).  Now, let  $ \, K := {\Bbbk[G]}_{\,
t}^{\,\vee} \, $,  $ \, H := {A_\Bbbk(G)}_t^{\,\prime} \, $;
by the previous description of  $ H $  and  $ K \, $,  the
natural pairing  $ \, H \times K \longrightarrow R \, $, 
\, which is perfect on the left, has  $ \, K_\infty :=
R\big[t^{-1}\big] \, J^\infty \, $  as right kernel: so
it induces a  {\sl perfect\/}  Hopf pairing  $ \, H \times
\overline{K} \! \longrightarrow R \, $,  \, where  $ \,
\overline{K} := K \big/ K_\infty \, $.  By construction
the latter specializes at  $ \, t = 0 \, $  to the natural
pairing  $ \, F[K_G] \times \u(\gerk_G) \longrightarrow \Bbbk
\, $,  \, which is perfect too.  Then we can apply  Proposition
4.4{\it (c)\/}  in [Ga1--2] (with  $ \overline{K} $  playing the
r\^{o}le of  $ K $  therein) which yields  $ \, \overline{K}^{\;
\prime} \! = \big( H^\vee \big)^\bullet \! = \Big(\! \big(
{A_\Bbbk(G)}_t^{\,\prime} \big)^{\!\vee} \Big)^{\!\bullet} \, $.
By definitions one sees that  $ \, R\big[t^{-1}\big] \, J^\infty
\subseteq \big({\Bbbk[G]}_{\,t}^\vee\big)' \, $,  \, whence  $ \,
\overline{K}^{\;\prime} \! = \big({\Bbbk[G]}_{\,t}^\vee\big)' \!
\Big/ \! \big( R\big[t^{-1}\big] \, J^\infty \big) \, $  follows
at once; Proposition 8.7 describes the latter space as
$ \; \overline{K}^{\;\prime} = \, \Big(\! \bigoplus_{\Sb  b_i \in B,
                                                  \; 0 < e_i < p  \\
            r \in \N, \; b_1 \precneqq \cdots \precneqq b_r  \endSb}
\hskip-3pt  R \cdot \overline{\psi}_{b_1}^{\;e_1} \cdots
\overline{\psi}_{b_r}^{\;e_r} \Big) \, $,
\; where  $ \, \overline{\psi}_{b_i} := \psi_{b_i} \mod
R\big[t^{-1}\big] \, J^\infty \, $  for all  $ i \, $.  Since
we saw that  $ \, \overline{K}^{\;\prime} \! = \Big(\! \big(
{A_\Bbbk(G)}_t^{\,\prime} \big)^{\!\vee} \Big)^{\!\bullet} \, $, 
\, and  $ \, \psi_g = t^{+1} \chi_g \, $,  \, this analysis yields
 \vskip-4pt
  $$  \big({A_\Bbbk(G)}_t^{\,\prime}\big)^{\!\vee} = \, \Big(\!
        {\textstyle \bigoplus_{\Sb  b_i \in B, \; 0 < e_i < p  \\
        r \in \N, \; b_1 \precneqq \cdots \precneqq b_r  \endSb}}
R \cdot t^{- \sum_i e_i} \rho_{b_1, \dots, b_r}^{e_1,\dots,e_r} \Big)
\; \cong \; \big(\overline{K}^{\;\prime}\big)^*  $$
 \vskip-3pt
\noindent
whence  $ \; \big({A_\Bbbk(G)}_t^{\,\prime}\big)^{\!\vee}\Big|_{t=0}
\! \cong \big(\overline{K}^{\;\prime}\big)^*\Big|_{t=0} \! = \big(
K_t^{\,\prime}\big|_{t=0} \big)^* \! = \Big(\!\big({\Bbbk[G]}_{\,t}^{\,
\vee}\big)'\big|_{t=0}\Big)^* \! \cong {F\big[\varGamma_G\big]}^* \!
= \u \big( \gerk_G^{\,\times} \big) = S\big(\gerk_G^{\,\times}\big)
\Big/ \big( \big\{ x^p \,\big|\, x \in \gerk_G^{\,\times} \big\} \big)
\; $  as claimed, the latter identity being trivial (as  $ \gerk_G^{\,
\times} $  is Abelian).   \qed
\enddemo

\vskip7pt

   {\bf 8.14 Remark:}  for each field  $ \Bbbk $  of positive
characteristic,  {\sl the functor  $ \, G \mapsto K_G \, $  is an
invariant for  $ G \, $},  a new one's to the author knowledge,
but again equivalent to  $ \gerk_G \, $.

\vskip7pt

   {\bf 8.15 Examples.}  {\sl (1) \, Finite Abelian  $ p \, $--groups.}
\, Let  $ p \, $  be a prime number,  $ \Bbbk $  a field with  $ \,
\Char(\Bbbk) = p \, $,  and  $ \, G := \Z_{p^{e_1}} \times \cdots
\times \Z_{p^{e_k}} \, $  ($ k, e_1, \ldots, e_k \in \N \, $), 
\, with  $ \, e_1 \geq \cdots \geq e_k \, $.
                                              \par
   First,  $ \, \gerk_G \, $  is Abelian, because  $ G $  is.
Let  $ g_i $  be a generator of  $ \, \Z_{p^{e_i}} \, $  (for
all  $ i \, $),  identified with its image in  $ G \, $.  Since
$ G $  is Abelian we have  $ \, G_{[n]} = G^{p^n} \, $  (for all
$ n \, $),  and an ordered  $ p $-l.c.s.-net  is  $ \, B :=
\bigcup_{r \in \N_+} B_r \, $  with  $ \, B_r := \Big\{\, g_1^{\,p^r},
\, g_2^{\,p^r}, \, \ldots , \, g_{j_r}^{\,p^r} \Big\} \, $  where
$ j_r $  is uniquely defined by  $ \; e_{j_r} > r \, $,  $ \; e_{j_r
+ 1} \leq r \, $.  Then  $ \, \gerk_G \, $  has  $ \Bbbk $--basis
$ \, {\big\{\, \overline{\,\eta_{g_i^{p^{s_i}}}} \,\big\}}_{1 \leq i
\leq k; \; 0 \leq s_i < e_i} \, $,  \, and minimal set of generators
(as a restricted Lie algebra)  $ \, \big\{\,\overline{\,\eta_{g_1}}
\, , \, \overline{\,\eta_{g_2}} \, , \, \ldots, \, \overline{\,
\eta_{g_k}} \,\big\} \, $.  In fact, the  $ p $--operation 
of  $ \gerk_G $  is  $ \, {\big(\overline{\,\eta_{g_i^{p^s}}}
\,\big)}^{[p\hskip0,5pt]} = \overline{\,\eta_{g_i^{p^{s+1}}}} \, $,
\, and the order of nilpotency of each  $ \, \overline{\,\eta_{g_i}}
\, $  is exactly  $ p^{e_i} $,  i.e.~the order of  $ g_i \, $.  In
addition  $ \, J^\infty = \{0\} \, $  so  $ \, {\Bbbk[G]}^\vee \!
= \Bbbk[G] \, $.  The outcome is  $ \; {\Bbbk[G]}^\vee \! = \,
\Bbbk[G] \; $  and
 \vskip11pt
   \centerline{ $  \widehat{\Bbbk[G]}  \, = \,  \u(\gerk_G)  \; = \;
U(\gerk_G) \bigg/ \Big( {\Big\{ {\big(\overline{\,\eta_{g_i^{p^s}}}
\,\big)}^p - \overline{\,\eta_{g_i^{p^{s+1}}}} \,\Big\}}_{1 \leq i
\leq k}^{0 \leq s < e_i} \bigcup \; {\Big\{ {\big(\overline{\,
\eta_{g_i^{p^{e_i-1}}}}\,\big)}^p \,\Big\}}_{1 \leq i \leq k}
\, \Big) $ }
 \vskip1pt
\noindent
whence  $ \; \widehat{\Bbbk[G]} \, \cong \, \Bbbk[x_1,\dots,x_k] \bigg/
\! \Big( \Big\{\, x_i^{p^{e_i}} \;\Big|\; 1 \leq i \leq k \,\Big\} \Big)
\, $,  \; via  $ \; \overline{\,\eta_{g_i^{p^s}}} \mapsto x_i^{\,p^s}
\, $  (for all  $ i $,  $ s \, $).
                                              \par
   As for  $ {\Bbbk[G]}_{\,t}^{\,\vee} $,  for all  $ \, r < e_i
\, $  we have  $ \, d\big(g_i^{p^r}\big) = p^r \, $  and so  $ \,
\chi_{g_i^{p^r}} = t^{-p^r} \big( g_i^{p^r} \!-\! 1 \big) \, $
and  $ \, \psi_{g_i^{p^r}} = t^{1-p^r} \big( g_i^{p^r} \!-\!
1 \big) \, $.  Now  $ \, G_{[\infty]} = \{1\} \, $  (or,
equivalently,  $ \, J^\infty = \{0\} \, $)  and everything
is Abelian, so from the general theory we conclude that both
$ {\Bbbk[G]}_{\,t}^{\,\vee} $  and  $ {\big( {\Bbbk[G]}_{\,
t}^{\,\vee} \big)}' $  are truncated-polynomial algebras, in
the  $ \chi_{g_i^{p^r}} $'s  and in the  $ \psi_{g_i^{p^r}} $'s
respectively, namely
 \vskip5pt
   \centerline{ $ \eqalign{
   {\Bbbk[G]}_{\,t}^{\,\vee}  &  = \; \Bbbk[t] \Big[ \big\{\,
\chi_{g_i^{p^s}} \big\}_{1 \leq i \leq k \, ; \; 0 \leq s < e_i} \Big]
\; \cong \;  \big(\Bbbk[t]\big)\,[\,y_1,\dots,y_k] \bigg/ \! \Big(
\Big\{\,y_i^{p^{e_i}} \;\Big|\; 1 \leq i \leq k \,\Big\} \Big)  \cr
   {\big({\Bbbk[G]}_{\,t}^{\,\vee}\big)}'  &  = \, \Bbbk[t] \Big[
\big\{\, \psi_{g_i^{p^s}} \big\}_{1 \leq i \leq k \, ; \; 0 \leq
s < e_i} \Big] \, \cong \, \big(\Bbbk[t]\big) \, \Big[ \big\{\,
z_{i,s} \big\}_{1 \leq i \leq k \, ; \; 0 \leq s < e_i} \Big]
\bigg/ \! \Big( \Big\{\, {z_{i,s}}^{\!p} \,\Big\}_{1 \leq i
\leq k} \,\Big)  \cr } $ }
 \vskip3pt
\noindent
via the isomorphisms given by  $ \, \overline{\,\chi_{g_i^{p^s}}}
\mapsto y_i^{\,p^s} \, $,  $ \, \overline{\,\psi_{g_i^{p^s}}}
\mapsto z_{i,s} \, $.  When  $ \, e_1 > 1 \, $  this implies 
$ \, {\big( {\Bbbk[G]}_{\,t}^{\,\vee} \big)}' \supsetneqq
{\Bbbk[G]}_{\,t} \, $.  Setting  $ \; \overline{\psi_{g_i^{p^s}}}
:= \psi_{g_i^{p^s}} \mod t \, {\big( {\Bbbk[G]}_{\,t}^{\,\vee}
\big)}' \; $  ($ \, 1 \! \leq \! i \! \leq \! k \, $,  $ \, 0
\! \leq \! s \! < \! e_i \, $)  we have
                              \hfill\break
   \centerline{ $ F\big[\varGamma_G\big] \, = \, {\big(
{\Bbbk[G]}_{\,t}^{\,\vee} \big)}'{\Big|}_{t=0} = \, \Bbbk \Big[
\big\{ \overline{\psi_{g_i^{p^s}}} \,\big\}_{1 \leq i \leq k}^{0 \leq
s < e_i} \Big] \, \cong \, \Bbbk\,\Big[ \big\{\, w_{i,s} \big\}_{1
\leq i \leq k}^{0 \leq s < e_i} \Big] \bigg/ \! \Big( \Big\{\,
w_{i,s}^{\,p} \Big|\; 1 \!\leq\! i \!\leq\! k \,\Big\} \Big) $ }
(via  $ \, \overline{\psi_{g_i^{p^s}}} \mapsto w_{i,s} \, $)
as a  $ \Bbbk $--algebra.  The Poisson bracket is trivial, and
the  $ w_{i,s} $'s  are primitive for  $ \, s > 1 \, $  and
$ \, \Delta(w_{i,1}) = w_{i,1} \otimes 1 + 1 \otimes w_{i,1}
+ w_{i,1} \otimes w_{i,1} \, $  for all  $ 1 \leq i \leq k \, $.
If instead  $ \; e_1 = \cdots = e_k = 1 \, $,  \, then  $ \, {\big(
{\Bbbk[G]}_{\,t}^{\,\vee} \big)}' = {\Bbbk[G]}_{\,t} \, $.  This
is an analogue of Theorem 2.2{\it (b)\/}  in [Ga1--2], although
now  $ \, \Char(\Bbbk) > 0 \, $,  \, in that in this case  $ \,
{\Bbbk[G]}_{\,t} \, $  is a QFA, with  $ \, {\Bbbk[G]}_{\,t}
{\Big|}_{t=0} \! = \Bbbk[G] = F\big[\widehat{G} \,\big] \, $
where  $ \, \widehat{G} \, $  is the  {\sl group of characters\/}
of  $ G \, $.  But then  $ \, F\big[\widehat{G}\,\big] = \Bbbk[G] =
{\Bbbk[G]}_{\,t}{\Big|}_{t=0} \! = {\big( {\Bbbk[G]}_{\,t}^{\,\vee}
\,\big)}'{\Big|}_{t=0} \! = F\big[\varGamma_G\big] \, $  (by the
general analysis), which means that  $ \widehat{G} $  is a finite
connected dimension 0 height 1 
   \hbox{group-scheme dual to  $ \gerk_G \, $,
namely  $ K_G^\star = \varGamma_G \, $.}   
                                                    \par
   Finally, a direct easy calculation shows that   --- letting  $ \,
\chi^*_g := t^{\,d(g)} \, (\varphi_g - \varphi_1) \in {A_\Bbbk(G)}'_t
\, $  and  $ \, \psi^*_g := t^{\,d(g)-1} \, (\varphi_g - \varphi_1)
\in {\big({A_\Bbbk(G)}_t^{\,\prime}\,\big)}^{\!\vee} \, $,  \, for
all  $ \, g \in G \setminus \{1\} \, $  ---   we have also
 \vskip7pt
   \centerline{ $  \eqalign{
   {A_\Bbbk(G)}_t^{\,\prime} \,  &  = \; \Bbbk[t] \Big[ \big\{\,
\chi^*_{g_i^{p^s}} \big\}_{1 \leq i \leq k}^{0 \leq s < e_i}
\Big] \, \cong \; \Bbbk[t] \Big[ \big\{ Y_{i,j} \big\}_{1
\leq i \leq k}^{0 \leq s < e_i} \Big] \bigg/ \! \Big( \big\{
Y_{i,j}^{\;p} \big\}_{1 \leq i \leq k}^{0 \leq s < e_i} \Big)  \cr
   {\big({A_\Bbbk(G)}_t^{\,\prime}\,\big)}^{\!\vee}  &  = \; \Bbbk[t]
\Big[ \big\{\, \psi^*_{g_i^{p^s}} \big\}_{1 \leq i \leq k}^{0 \leq
s < e_i} \Big] \, \cong \; \Bbbk[t] \Big[ \big\{ Z_{i,s} \big\}_{1
\leq i \leq k}^{0 \leq s < e_i} \Big] \bigg/ \! \Big( \big\{
Z_{i,s}^{\;p} - Z_{i,s} \big\}_{1 \leq i \leq k}^{0 \leq s
< e_i} \Big)  \cr } $ }
 \vskip3pt
\noindent
via the isomorphisms given by  $ \, \chi^*_{g_i^{p^s}} \mapsto
Y_{i,s} \, $  and  $ \, \psi^*_{g_i^{p^s}} \mapsto Z_{i,s} \, $,
\, from which one also gets the analogous descriptions of  $ \,
{A_\Bbbk(G)}_t^{\,\prime}{\Big|}_{t=0} \! = \widetilde{A_\Bbbk(G)}
= F[K_G] \, $  and of  $ \, {\big({A_\Bbbk(G)}_t^{\,\prime}\,
\big)}^{\!\vee}{\Big|}_{t=0} \! = \u(\gerk_G^\times) \; $.
%
%
 \eject   
   {\sl (2) \, A non-Abelian  $ p \, $--group.} \, Let  $ p \, $
be a prime number, and let  $ \, \Char(\Bbbk) = p > 0 \, $.
                                              \par
   Let  $ \, G := \Z_p \ltimes \Z_{p^{\,2}} \, $,  \, that is the group
with generators  $ \, \nu $,  $ \tau \, $  and relations  $ \, \nu^p
= 1 \, $,  $ \, \tau^{p^2} = 1 \, $,  $ \, \nu \, \tau \, \nu^{-1}
= \tau^{1+p} \, $.  In this case,  $ \, G_{[2]} = \cdots = G_{[p\,]}
= \big\{ 1, \tau^p \,\big\} \, $,  $ \, G_{[p+1]} = \{1\} \, $,  \,
so we can take  $ \, B_1 = \{\nu \, , \tau \,\} \, $  and  $ \, B_p
= \big\{ \tau^p \,\big\} \, $  to form an ordered  $ p $-l.c.s.-net
$ \, B := B_1 \cup B_p \, $  w.r.t.~the ordering  $ \, \nu \preceq
\tau \preceq \tau^p \, $.  Noting also that  $ \, J^\infty = \{0\}
\, $  (for  $ \, G_{[\infty]} = \{1\} \, $),  we have
 \vskip5pt
   \centerline{ $ {\Bbbk[G]}^\vee_t \; = \; {\textstyle \bigoplus_{a,
b, c = 0}^{p-1}} \, \Bbbk[t] \cdot \chi_\nu^{\,a} \, \chi_\tau^{\,b}
\, \chi_{\tau^p}^{\,c} \; = \; {\textstyle \bigoplus_{a,b,c=0}^{p-1}}
\, \Bbbk[t] \, t^{- a - b - c \, p} \cdot {(\nu-1)}^a \, {(\tau-1)}^b
\, {\big( \tau^p - 1 \big)}^c $ }
 \vskip5pt
\noindent
as  $ \Bbbk[t] $--modules,  since  $ \, d(\nu) = 1 = d(\tau) \, $
and  $ \, d\big(\tau^p\big)) = p \, $,  \, with  $ \, \Delta(\chi_g)
= \chi_g \otimes 1 + 1 \otimes \chi_g + t^{d(g)} \, \chi_g \otimes
\chi_g \, $  for all  $ \, g \in B \, $.  As a direct consequence
we have also
 \vskip5pt
   \centerline{ $ {\textstyle \bigoplus_{a,b,c=0}^{p-1}} \, \Bbbk
\cdot \overline{\chi_\nu}^{\;a} \, \overline{\chi_\tau}^{\;b} \,
\overline{\chi_{\tau^p}}^{\;c} \; = \; {\Bbbk[G]}^\vee_t
{\Big|}_{t=0} \; \cong \; \widehat{\Bbbk[G]} \; = \; {\textstyle
\bigoplus_{a,b,c=0}^{p-1}} \, \Bbbk \cdot \overline{\eta_\nu}^{\;a}
\, \overline{\eta_\tau}^{\;b} \, \overline{\eta_{\tau^p}}^{\;c}
\, . $ }
 \vskip3pt
   The relations  $ \, \nu^p = 1 \, $  and  $ \, \tau^{p^2} = 1 \, $
in  $ G $  yield trivial relations in  $ \Bbbk[G] $  and  $ \big(
\Bbbk[t] \big)[G] \, $.  Instead, the relation  $ \, \nu \, \tau \,
\nu^{-1} = \tau^{1+p} \, $  turns into  $ \, [\,\eta_\nu, \eta_\tau]
= \eta_{\tau^p} \cdot \tau \, \nu \, $,  \, whence  $ \, [\,\chi_\nu,
\chi_\tau] = t^{p-2} \, \chi_{\tau^p} \cdot \tau \, \nu \, $
in  $ {\Bbbk[G]}_{\,t}^{\,\vee} \, $;  \, thus  $ \; [\,
\overline{\chi_\nu} \, , \, \overline{\chi_\tau} \,] =
\delta_{p,2} \, \overline{\chi_{\tau^p}} \, $.  Since  $ \,
[\,\overline{\chi_\tau} \, , \, \overline{\chi_{\tau^p}} \,]
= 0 = [\,\overline{\chi_\nu} \, , \, \overline{\chi_{\tau^p}} \,]
\, $  (for  $ \, \nu \, \tau^p \, \nu^{-1} = {\big( \tau^{1+p} \big)}^p
= \tau^{p + p^2} = \tau^p \, $)  and  $ \, \{ \overline{\chi_\nu} \, ,
\, \overline{\chi_\tau} \, , \, \overline{\chi_{\tau^p}} \,\} \, $  is
a  $ \Bbbk $--basis  of  $ \, \gerk_G = \Cal{L}_p(G) \, $,  \, the
latter has trivial or non-trivial Lie bracket, according to whether
$ \, p \not= 2 \, $  or  $ \, p = 2 \, $.  In addition,  $ \,
\chi_\nu^{\;p} = 0 \, $,  $ \, \chi_{\tau^p}^{\;p} = 0 \, $
and  $ \, \chi_\tau^{\;p} = \chi_{\tau^p} \, $:  \, these give
analogous relations in  $ \, {\Bbbk[G]}_{\,t}^{\,\vee}{\Big|}_{t=0}
\!\!\! = \widehat{\Bbbk[G]} \, $,  \, which read as formulas for
the  $ p $--operation,
   \hbox{namely  $ \; \overline{\chi_\nu}^{\;[p\,]} = 0 \, $,
$ \; \overline{\chi_{\tau^p}}^{\;[p\,]} = 0 \, $,  $ \;
\overline{\chi_\tau}^{\;[p\,]} = \chi_{\tau^p} \, $.}
                                     \par
   Reassuming, we have a complete presentation for
$ {\Bbbk[G]}^\vee_t $  by generators and relations, i.e.   
 \vskip9pt
  \centerline{ $ {\Bbbk[G]}^\vee_t  \;\; \cong \;\;
\Bbbk[t] \, \big\langle x_1, x_2, x_3 \big\rangle \bigg/
\! \bigg( \hskip-11pt
  \hbox{ $ \matrix
 &  \! x_1 \, x_2 - x_2 \, x_1 - t^{p-2} \, x_3
\, (1 + t \, x_\tau) \, (1 + t \, x_\nu) \!  \\
 &  \! x_1 \, x_3 - x_3 \, x_1 \, ,  \quad  x_1^{\,p} \, ,
\quad  x_2^{\,p} - x_3 \, , \quad  x_3^{\,p} \, , \quad
x_2 \, x_3 - x_3 \, x_2 \!
       \endmatrix $ }  \bigg) $ }
 \vskip9pt
\noindent
via  $ \, \chi_\nu \mapsto x_1 \, $,  $ \, \chi_\tau \mapsto x_2
\, $,  $ \, \chi_{\tau^p} \mapsto x_3 \, $.  Similarly (as a
consequence) we have the presentation
 \vskip9pt
  \centerline{ $ \widehat{\Bbbk[G]}  \, = \,  {\Bbbk[G]}^\vee_t
{\Big|}_{t=0} \; \cong \;\;  \Bbbk \, \big\langle y_1,y_2,y_3
\big\rangle \bigg/ \! \bigg(
 \hbox{ $ \matrix
     \! y_1 \, y_2 - y_2 \, y_1 - \delta_{p,2} \, y_3 \, ,
\quad  y_2^{\,p} - y_3 \!  \\
     \! y_1 \, y_3 - y_3 \, y_1 \, ,  \quad  y_1^{\,p} \, ,
\quad  y_3^{\,p} \, ,  \quad  y_2 \, y_3 - y_3 \, y_2 \!
       \endmatrix $ }  \bigg) $ }
 \vskip9pt
\noindent
via  $ \, \overline{\chi_\nu} \mapsto y_1 \, $,  $ \,
\overline{\chi_\tau} \mapsto y_2 \, $,  $ \, \overline{\chi_{\tau^p}}
\mapsto y_3 \, $,  \, with  $ p $--operation  as above and the
           $ y_i $'s  being primitive\break
{\it  $ \underline{\text{Remark}} $:}  \, if  $ \, p \not= 2 \, $
exactly the same result holds for  $ \, G = \Z_p \times \Z_{p^2} \, $,
\, i.e.~$ \; \gerk_{\, \Z_p \hskip-0,5pt \ltimes \Z_{p^2}} =
\gerk_{\, \Z_p \hskip-0,5pt \times \Z_{p^2}} \; $:  \; this
shows that the restricted Lie bialgebra  $ \gerk_G $  may be
     \hbox{not enough to recover the group  $ G \, $.}
                                          \par
   As for  $ \, {\big( {\Bbbk[G]}^\vee_t \big)}' $,  \, it is
generated by  $ \, \psi_\nu = \nu - 1 $,  $ \, \psi_\tau = \tau
- 1 $,  $ \, \psi_{\tau^p} = t^{1-p} \big(\tau^p - 1 \big)
\, $,  \, with relations  $ \; \psi_\nu^{\;p} = 0 \, $,
$ \; \psi_\tau^{\;p} = t^{p-1} \psi_{\tau^p} \, $,  $ \;
\psi_{\tau^p}^{\;\;p} = 0 \, $,  $ \; \psi_\nu \, \psi_\tau
- \psi_\tau \, \psi_\nu = t^{\,p-1} \psi_{\tau^p} \, (1 +
\psi_\tau) \, (1 + \psi_\nu) \, $,  $ \; \psi_\tau \, \psi_{\tau^p}
- \psi_{\tau^p} \, \psi_\tau = 0 \, $,  and  $ \; \psi_\nu \,
\psi_{\tau^p} - \psi_{\tau^p} \, \psi_\nu = 0 \, $.  In particular
$ \; {\big({\Bbbk[G]}^\vee_t\big)}' \supsetneqq \big(\Bbbk[t]\big)[G]
\, $,  \, and
 \vskip9pt
  \centerline{ $ {\big({\Bbbk[G]}^\vee_t\big)}'  \;\; \cong \;\;
\Bbbk[t] \, \big\langle y_1, y_2, y_3 \big\rangle \bigg/ \!
\bigg( \hskip-12pt \hbox{ $ \matrix
   y_\nu \, y_{\tau^p} - y_{\tau^p} \, y_\nu \, ,
\quad  y_\tau^{\;p} - t^{p-1} y_{\tau^p}  \, ,  \quad
y_\tau \, y_{\tau^p} - y_{\tau^p} \, y_\tau  \\
   \quad  y_\nu^{\;p} \, ,  \quad  y_\nu \, y_\tau - y_\tau \,
y_\nu - t^{\,p-1} y_{\tau^p} \, (1 + y_\tau) \, (1 + y_\nu)
\, , \quad  y_{\tau^p}^{\;\;p}
       \endmatrix $ }  \hskip-3pt \bigg) $ }
 \vskip9pt
\noindent
via  $ \, \psi_\nu \mapsto y_1 \, $,  $ \, \psi_\tau \mapsto y_2
\, $,  $ \, \psi_{\tau^p} \mapsto y_3 \, $.  Letting  $ \; z_1 :=
\psi_\nu{\big|}_{t=0} \! + 1 \, $,  $ \; z_2 := \psi_\tau{\big|}_{t=0}
\! + 1 \; $  and  $ \; x_3 := \psi_{\tau^p}{\big|}_{t=0} \; $  this gives
$ \, {\big( {\Bbbk[G]}^\vee_t \big)}'{\Big|}_{t=0} \!\! = \Bbbk \big[ z_1,
z_2, x_3 \big] \Big/ \big( z_1^{\,p} \! - \! 1, z_2^{\,p} \! - \! 1,
x_3^{\,p} \,\big) \, $  as a  $ \Bbbk $--algebra,  with the  $ z_i $'s
group-like,  $ x_3 $  primitive  (cf.~Theorem 8.8{\it(b)\/}),  and
Poisson bracket given by  $ \, \big\{ z_1, z_2 \big\} = \delta_{p,2}
\, z_1 \, z_2 \, x_3 \, $,  $ \, \big\{ z_2, x_3 \big\} = 0 \, $
and  $ \, \big\{z_1, x_3\big\} = 0 \, $.  Thus  $ \, {\big(
{\Bbbk[G]}^\vee_t \big)}'{\Big|}_{t=0} \! = F[\varGamma_G]
\, $  with  $ \, \varGamma_G \cong {\boldsymbol\mu}_p \times
{\boldsymbol\mu}_p \times {\boldsymbol\alpha}_p \, $  as
algebraic groups, with Poisson structure such that  $ \,
\text{\sl coLie}\,(\varGamma_G) \cong \gerk_G \, $.
                                         \par
   Since  $ \, G_\infty = \{1\} \, $  the general theory ensures that
$ \, {A_\Bbbk(G)}' = A_\Bbbk(G) \, $.  We leave to the interested
reader the task of computing the filtration  $ \underline{D} $  of
$ A_\Bbbk(G) $,  and consequently describe  $ \, {A_R(G)}' \, $,
$ \, \big({A_R(G)}'\big)^{\!\vee} \, $,  $ \, \widetilde{A_\Bbbk(G)}
\, $  and the connected Poisson group  $ \, K_G := \text{\it Spec}\,
\big( \widetilde{A_\Bbbk(G)} \big) \, $.

\vskip1,7truecm

\vskip33pt

\enddocument

\bye
\end